\def\longversion{1}

\ifx\longversion\undefined
	\documentclass[final]{siamltex}
\else
	\documentclass[letterpaper]{article}
\fi

\usepackage[margin=1.5in]{geometry}

\usepackage{times}

\usepackage{graphics} % for pdf, bitmapped graphics files
\usepackage{epsfig} % for postscript graphics files
\ifx\longversion\undefined
	
\else
	\usepackage{amsthm}
\fi
\usepackage{amsmath} % assumes amsmath package installed
\usepackage{amssymb}  % assumes amsmath package installed
\usepackage{mathrsfs}
\usepackage{color}
\usepackage{subfigure}
\usepackage{framed}
\usepackage{stackrel}
\usepackage{soul}
\usepackage{cancel}
\usepackage[all]{xy}

\definecolor{MyDarkRed}{rgb}{0.5,0,0.1}
\definecolor{MyDarkBlue}{rgb}{0.1,0.1,0.5}
\definecolor{MyDarkGreen}{rgb}{0.1,0.5,0}
\definecolor{MyRed}{rgb}{1.0,0,0}
\definecolor{MyBlue}{rgb}{0,0,1.0}
\definecolor{MyGreen}{rgb}{0,0.8,0}
\definecolor{lightgray}{rgb}{0.96,0.96,0.96}
\usepackage[pdftex, colorlinks=true, urlcolor=MyDarkBlue, citecolor=MyDarkRed, linkcolor=MyDarkGreen]{hyperref}

\usepackage{helvet}
\usepackage{courier}
\usepackage{xcolor}
\usepackage{chemarr}
\usepackage{extarrows}
\usepackage[small]{caption}

\newcommand{\wipe}[1]{}

\definecolor{MyDarkRed}{rgb}{0.5,0,0.1}
\definecolor{MyRed}{rgb}{1.0,0,0}
\definecolor{MyBlue}{rgb}{0,0,1.0}
\definecolor{MyGreen}{rgb}{0,0.8,0}
\definecolor{MyPurple}{rgb}{0.6,0,0.6}
\definecolor{MyCyan}{rgb}{0,0.5,0.6}

\definecolor{lightgray}{rgb}{0.96,0.96,0.96}
\definecolor{leftBarGray}{rgb}{0.8,0.8,0.8}

\DeclareFixedFont{\picturechar}{OT1}{cmr}{m}{n}{5}

\newcommand{\changedA}[1]{{#1}}
\newcommand{\changedB}[1]{{#1}}
\newcommand{\changedLS}[1]{{#1}}
\newcommand{\changedC}[1]{{#1}}

\ifx\longversion\undefined
	\newcommand{\lo}[1]{}
	\newcommand{\sh}[1]{\changedLS{#1}}
	\newcommand{\ls}[2]{\changedLS{#2}}
\else
	\newcommand{\lo}[1]{}
	\newcommand{\sh}[1]{}
	\newcommand{\ls}[2]{}
	\newcommand{\Lo}[1]{{#1}}
	\newcommand{\Ls}[2]{{#1}}
\fi
\ifx\longversion\undefined
	
	\newtheorem{note}{Note}
	\newtheorem{remark}{Remark}
\else
	
	\newtheorem{theorem}{Theorem}
	\newtheorem{proposition}{Proposition}
	\newtheorem{corollary}{Corollary}
	\theoremstyle{definition}
	\newtheorem{definition}{Definition}
	\newtheorem{remark}{Remark}

\fi

\newcommand{\myurl}[1]{\textbf{\small\url{#1}}}

\newcommand{\indens}{\hspace*{5mm}}
\newcommand{\sqed}{\mbox{\ \ }\rule{5pt}{5pt}\medskip}

\renewcommand{\d}[0]{\textrm{d}}

\newcommand{\craig}[1]{}

\renewcommand{\d}{\textrm{d}}

\renewcommand{\setminus}{ - }

\newenvironment{myleftbar}{%
  \MakeFramed {\advance\hsize-\width \FrameRestore}}%
 {\endMakeFramed}
\newenvironment{quoteproof}{

\noindent\textit{Proof.}~ \begin{myleftbar}\small}{\sqed\end{myleftbar}}

\begin{document}

%\title{Using Linking Number to Identify Homology Classes of Cycles Via Action of Co-cycles}
\title{\changedA{Invariants for Homology Classes with Application to Optimal Search and Planning Problem in Robotics}}
\date{}
\author{Subhrajit Bhattacharya \and David Lipsky \and Robert Ghrist \and Vijay Kumar}
\maketitle

\begin{abstract}
We consider planning problems on a punctured Euclidean spaces,
$\mathbb{R}^D - \widetilde{\mathcal{O}}$, where
$\widetilde{\mathcal{O}}$ is a collection of obstacles. Such spaces
are of frequent occurrence as configuration spaces of robots, where
$\widetilde{\mathcal{O}}$ represent either physical obstacles that the
robots need to avoid (e.g., walls, other robots, etc.) or illegal
states (e.g., all legs off-the-ground). As state-planning is
translated to path-planning on a configuration space, we collate
equivalent plannings via topologically-equivalent paths. This prompts
finding or exploring the different homology classes in such
environments and finding representative optimal trajectories in each
such class.

In this paper we start by considering the problem of finding a
complete set of easily computable homology class invariants for
$(N-1)$-cycles in $(\mathbb{R}^D - \widetilde{\mathcal{O}})$. We
achieve this by finding explicit generators of the $(N-1)^{st}$ de
Rham cohomology group of this punctured Euclidean space, and using
their integrals to define cocycles. The action of those dual cocycles
on $(N-1)$-cycles gives the desired complete set of invariants. We
illustrate the computation through examples.

We further show that, due to the integral approach, this complete set
of invariants is well-suited for efficient search-based planning of
optimal robot trajectories with topological constraints. Finally we
extend this approach to computation of invariants in spaces derived
from $(\mathbb{R}^D - \widetilde{\mathcal{O}})$ by collapsing
subspace, thereby permitting application to a wider class of
non-Euclidean ambient spaces.

\end{abstract}

\section{Introduction}

\subsection{Motivation: Robot Path Planning with Topological Constraints}

\begin{figure*}
  \begin{center}
    \subfigure[In $\mathbb{R}^2$ punctured by obstacles.]{ \label{fig:robotics-motivation-2d}
      \includegraphics[width=0.4\textwidth, trim=0 0 0 0, clip=true]{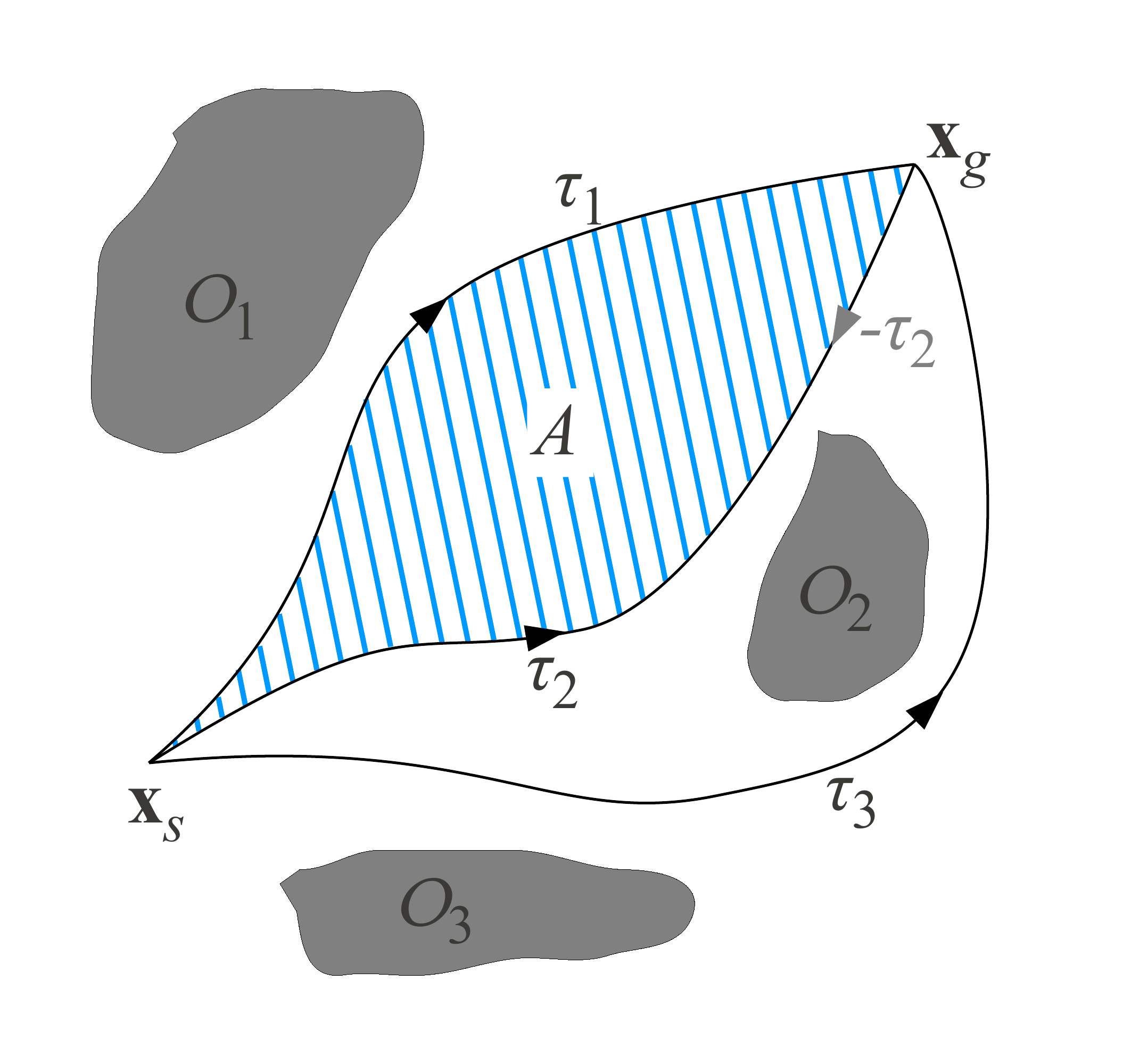}
    }
    \subfigure[In $\mathbb{R}^3$ punctured by obstacles]{ \label{fig:robotics-motivation-3d}
      \includegraphics[width=0.4\textwidth, trim=0 0 0 0, clip=true]{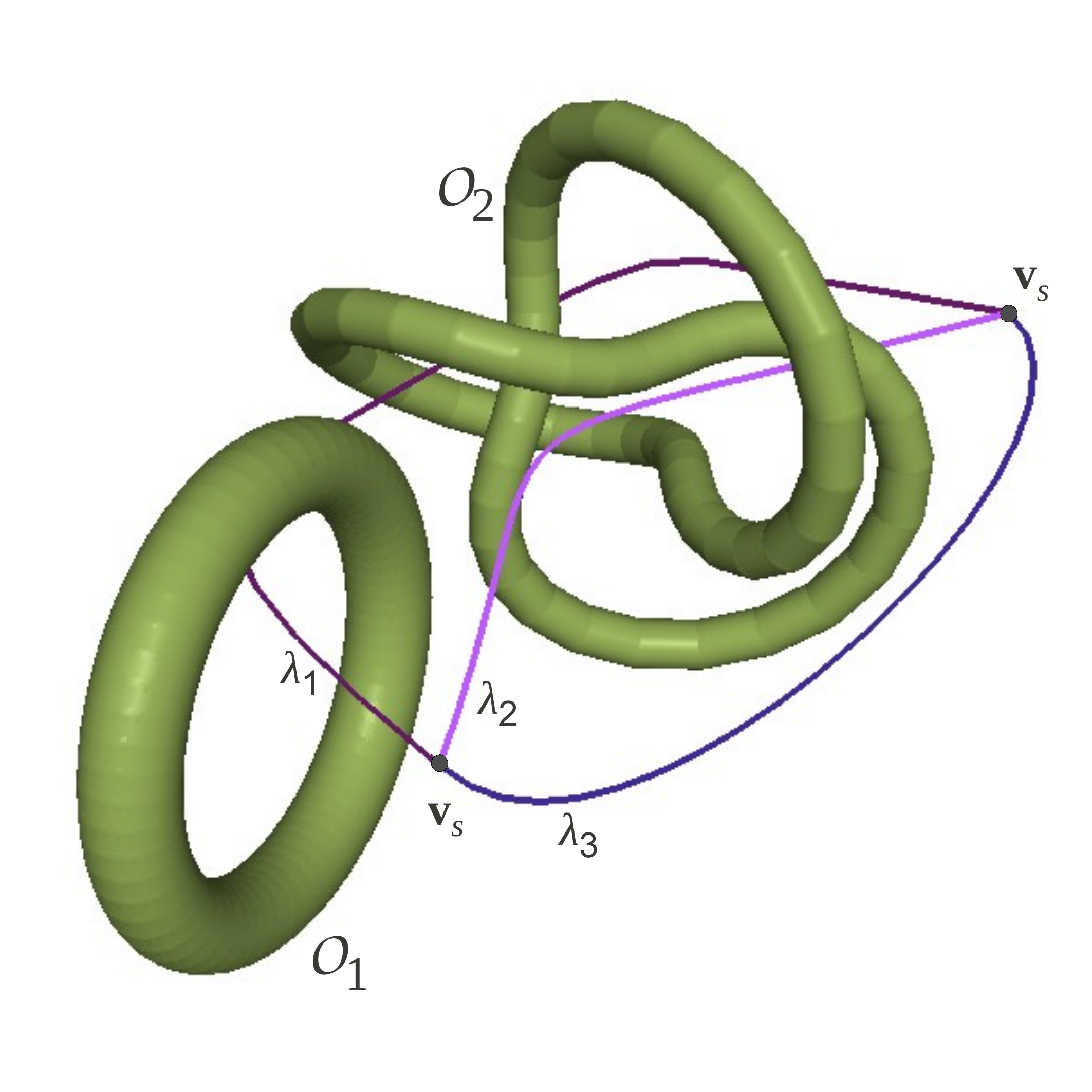}
    }
  \end{center}
  \caption{Homology classes of robot trajectories in Euclidean spaces with obstacles.} \label{fig:robotics-motivation} 
\end{figure*}

In numerous robotics applications, it is important to distinguish
between configuration space paths in different topological classes, as
a means of categorizing continuous families of plans. This motivation
--- connected components of paths relative to endpoints --- leads to
classifying up to homotopy. Examples motivating a classification of
homotopy classes of paths include: (1)
group exploration of an environment
\cite{Bourgault02informationbased}, in which an efficient strategy
involves allocating one agent per homotopy class; (2) visibility,
especially in the tracking of uncertain agents in an environment with
dynamic obstacles \cite{occlusion:yan:06}; and (3) multi-agent
coordination, in which (Pareto-) optimal planning coincides with
homotopy classification \cite{GL:2006}.

%add the following reference to the biblio:
%
%R. Ghrist and S. LaValle (2006) “Nonpositive curvature and Pareto
%optimal motion planning,” SIAM J. Control & Opt., 45(5), 1697-1713.

Although homotopy is a natural topological equivalence relation for
paths, the computational bottlenecks involved, especially in higher
dimensional configuration spaces, present severe challenges in solving
practical problems in robot path planning. Thus we resort to its
computationally-simpler cousin --- homology
(Figure~\ref{fig:robotics-motivation}). We assume a basic familiarity
with first-year algebraic topology, as in \cite{Hatcher:AlgTop} for homology
and \cite{bott1982differential} for differential forms and de Rham
cohomology.

The methods we employ, following
\cite{planning:AURO:12}, construct an explicit
differential $1$-form, the integration of which along
trajectoriesgives complete homology class invariants. Such $1$-forms
are elements of the de Rham cohomology group of the configuration
space,  $H_{dR}^1(\mathbb{R}^D-\widetilde{\mathcal{O}})$. To deal with
the obstacles, we replace ${\mathcal{O}}$ with topologically
equivalent codimesion-2 skeleta (e.g., Figure~\ref{fig:obstacle-eqvs})
and then compute the degrees (or \emph{linking numbers}) of closed
loops with the skeleta.

\subsection{Contributions of this Paper}

We generalize the path-planning problem to higher homology classes and
linking numbers of arbitrary submanifolds (not merely $1$-dimensional
curves representing trajectories). In particular, we will consider
$(N-1)$-dimensional closed manifolds as generalization of
$1$-dimensional curves that constituted the trajectories. Obstacles
will be represented by codimension $N$ closed manifolds (which, in
many cases will be deformation retracts of the original obstacles).

Degree and linking numbers are closely related to homology
\cite{Hatcher:AlgTop,dold1995lectures}. We will in fact show that the
proposed integration along trajectories give homology class invariants
for closed loops (something that was claimed in \cite{planning:AURO:12},
but not proved rigorously).

The primary aim of this paper is two-fold:
\begin{itemize}
 \item[1.] To find certain differential $(N-1)$-forms in the Euclidean
space punctured by obstacles, and show that integration of the forms
along $(N-1)$-dimensional closed manifolds give a \emph{complete set
of invariants} for homology classes of the manifolds in the punctured
space (\emph{i.e.} the value of the integral over two closed manifolds
are equal if and only if the manifolds are homologous),
 \item[2.] To adapt and extend the tools used in
\cite{planning:AURO:12} for robot path planning with topological
reasoning to arbitrary dimensional Euclidean configuration spaces
punctured by obstacles.
\end{itemize}

\subsection{Overview and Organization of this Paper}

\changedB{
The main concept behind the treatment in this paper is to exploit the
pairing $H^{N-1}(\mathbb{R}^D-\widetilde{\mathcal{O}}; \mathbb{G})
\otimes H_{N-1}(\mathbb{R}^D-\widetilde{\mathcal{O}}; \mathbb{G})
\rightarrow \mathbb{G}$, which evaluates $(N-1)$-cocycles over
$(N-1)$-cycles. Given a cycle $\overline{\omega} \in
Z_{N-1}(\mathbb{R}^D-\widetilde{\mathcal{O}}; \mathbb{G})$, and a
large enough set of cocycles, $\mathscr{A} = \{ \alpha_1, \alpha_2,
\cdots, \alpha_m \}, \alpha_i \in
Z^{N-1}(\mathbb{R}^D-\widetilde{\mathcal{O}}; \mathbb{G})$, one can
hope that the set of values $\{ \alpha_1(\overline{\omega}),
\alpha_2(\overline{\omega}), \cdots, \alpha_m(\overline{\omega}) \}
\in \mathbb{G}^m$ will provide some information about the homology
class of $\overline{\omega}$, that is the value of
$[\overline{\omega}] \in H_{N-1}(\mathbb{R}^D-\widetilde{\mathcal{O}};
\mathbb{G})$. In fact choosing the coefficients in $\mathbb{R}$, and
with some assumptions on $\widetilde{\mathcal{O}}$, we will show that
it is sufficient to choose the elements of $\mathscr{A}$ such that
their cohomology classes generate
$H^{N-1}(\mathbb{R}^D-\widetilde{\mathcal{O}}; \mathbb{R})$.

However, the challenge lies in explicitly finding the cochains,
$\alpha_i$, that will serve our purpose and are easy to evaluate on
cycles. Due to De Rham's theorem, the cocycles, $\alpha_i$, can be
represented by some $(N-1)$-form, $\phi_i \in
\Omega^{N-1}(\mathbb{R}^D-\widetilde{\mathcal{O}})$, so that the
evaluation of the cocycle over a cycle is, precisely, integration of
the form over the cycle.
In order to find this form, we exploit the difference map $p:
(\mathbb{R}^D-\widetilde{\mathcal{O}}) \times \widetilde{\mathcal{O}}
\rightarrow (\mathbb{R}^D-\{0\})$. The codomain of this map is the
$D$-dimensional Euclidean space with the origin removed, and is much
simpler and well-studied. Thus, if $\eta_0 \in
\Omega^{D-1}(\mathbb{R}^D-\{0\})$ is a differential $(D-1)$-form in
$(\mathbb{R}^D-\{0\})$, a simple pull-back via $p$ gives the form
$\eta = p^{*} \eta_0 \in
\Omega^{D-1}(\mathbb{R}^D-\widetilde{\mathcal{O}}) \times
\widetilde{\mathcal{O}})$. Upon integration of $\eta$ over some
$(D-N)$-cycle, $\overline{S}$, one may hope to obtain the desired
$(N-1)$-form, $\phi_i = \int_{\overline{S}} p^{*} \eta_0$.
Considering the space $(\mathbb{R}^D-\widetilde{\mathcal{O}}) \times
\widetilde{\mathcal{O}}$ as a fiber bundle over
$(\mathbb{R}^D-\widetilde{\mathcal{O}})$ with
$\widetilde{\mathcal{O}}$ as the fibers, one may be tempted to
integrate $p^{*} \eta_0$ over the fibers. However, the nature of
$\widetilde{\mathcal{O}})$ (its topology, dimensionality) can be quite
arbitrary in general.

Thus we begin by constructing a suitable skeleton
$\widetilde{\mathcal{S}})$ with which to replace
$\widetilde{\mathcal{O}})$, so that the spaces
$(\mathbb{R}^D-\widetilde{\mathcal{O}})$ and
$(\mathbb{R}^D-\widetilde{\mathcal{S}})$ are identical as far as their
$(N-1)^{th}$ homology groups are concerned. However, in that
construction, we will ensure that $\widetilde{\mathcal{S}}$ is a
collection (disjoint union) of codimension-$N$ manifolds, thus
simplifying the problem.

}

\lo{
The rest of the paper is organized as follows:

In Section~\ref{sec:simplify} we first try to simplify the problem by suggesting some replacement of the arbitrary punctures/holes (obstacles in light of robot planning problems) that may be present in an Euclidean space. In particular, we suggest that in a $D$ dimensional punctured Euclidean space, in which we are concerned about computing the homology classes of $(N-1)$-cycles (top-dimensional covers of $(N-1)$-dimensional closed manifolds), we can replace the punctures/holes by some $(D-N)$-dimensional representatives under certain conditions. Thus we obtain the reduced problem definition of Section~\ref{sec:reduced-problem-algtop}.

In Section~\ref{sec:linking-num-pre} we describe some concepts on linking number and their relation to homology classes. This section discusses the topics from a fairly general algebraic topology point of view without specializing to the main problem under consideration. However, in the process, we try to illustrate some of the technical details using simple examples in the $D=3,N=2$ case (which is the case of robot planning problems in $3$-dimensional configuration space). %The main aim of this section is to obtain a general formula (in form of an integration) for computing linking numbers, and showing that linking numbers between two cycles actually is a complete invariants of homology class of one of the cycles in the space formed by removing the other cycle

In Section~\ref{sec:homology-invariant-formula} we specialize some of the results obtained in Section~\ref{sec:linking-num-pre} to fit the reduced problem we described in Section~\ref{sec:reduced-problem-algtop}. We hence obtain an explicit formula for the complete set of invariants for homology class. This is the invariant described in Equation~\eqref{eq:final-psi-eqn}.

In Section~\ref{sec:theo-validation} we plug in some specific values for $D$ and $N$ in the formula we obtained in the previous section, and hence show that the specific formulae we obtain match exactly with some of the well-known results from complex analysis, electromagnetism and electrostatics. Thus we conclude that not only we have unified some seemingly unrelated theories, but have also generalized them to higher dimensions.

In Section~\ref{sec:example-applications} we demonstrate one example with $D=5,N=3$, and show that the proposed formula indeed computes a complete invariant for homology class in that example. Moreover, we illustrate how the proposed formula can be used in search-based robot path planning with topological constraints.

\changedA{Finally, in Section~\ref{sec:non-Eu-extension}, we attempt to achieve certain level of generalization of the discussed approach to some ambient spaces that are not Euclidean. We achieve this by considering subspace, $L$, of the punctured Euclidean space, and trivializing every chain in it. This, in effect, lets us consider homology classes of relative cycles in $(\mathbb{R}^D-\widetilde{\mathcal{O}}, L)$. We demonstrate that the invariants thus computed for such spaces can once again be efficiently used for search-based path planning in robotics.}
}

\sh{
The rest of the paper is organized as follows:
In Section~\ref{sec:simplify} we first try to simplify the problem by suggesting the replacements, $\widetilde{\mathcal{S}}$, of $\widetilde{\mathcal{O}}$. Thus we obtain the reduced problem definition of Section~\ref{sec:reduced-problem-algtop}.
In Section~\ref{sec:linking-num-pre} we describe some concepts on linking number, \changedB{where we introduce the map $p$ for the first time from a fairly general perspective}. We thus establish the relation of linking numbers with homology classes.
In Section~\ref{sec:homology-invariant-formula} we specialize some of the results obtained in Section~\ref{sec:linking-num-pre} to fit the reduced problem we described in Section~\ref{sec:reduced-problem-algtop}. We hence obtain an explicit formula for the complete set of invariants for homology class. This is the invariant described in Equation~\eqref{eq:final-psi-eqn}.
\changedA{We discuss how plugging in some specific low values of $D$ and $N$ in the formula gives us some well-known results from complex analysis, electromagnetism and electrostatics.}
In Section~\ref{sec:example-applications} we demonstrate one example with $D=5,N=3$, and show that the proposed formula indeed computes a complete invariant for homology class in that example. Moreover, we illustrate how the proposed formula can be used in search-based robot path planning with topological constraints.
\changedA{Finally, in Section~\ref{sec:non-Eu-extension}, we attempt to achieve certain level of generalization of the discussed approach to ambient spaces that are not Euclidean. We achieve this by considering subspace, $L$, of the punctured Euclidean space, and trivializing every chain in it. This, in effect, lets us consider homology classes of relative cycles in $(\mathbb{R}^D-\widetilde{\mathcal{O}}, L)$. We demonstrate that the invariants thus computed for such spaces can once again be efficiently used for search-based path planning in robotics.}
}

\vspace{0.1in}
Throughout this paper we consider homology and cohomology with coefficients in the field $\mathbb{R}$. \changedB{As a consequence, all the homology and cohomology groups are freely and finitely generated.}
Also, for simplicity, we will throughout consider $N>1$ to avoid the special treatment of the $0^{th}$ (co)homology groups. \changedA{All topological spaces referred to in this paper are assumed to be Hausdorff.}

\section{On Building Obstacle Equivalents} \label{sec:simplify}

\begin{figure*}[t]
  \begin{center}
    \includegraphics[width=0.65\textwidth, trim=100 90 100 80, clip=true]{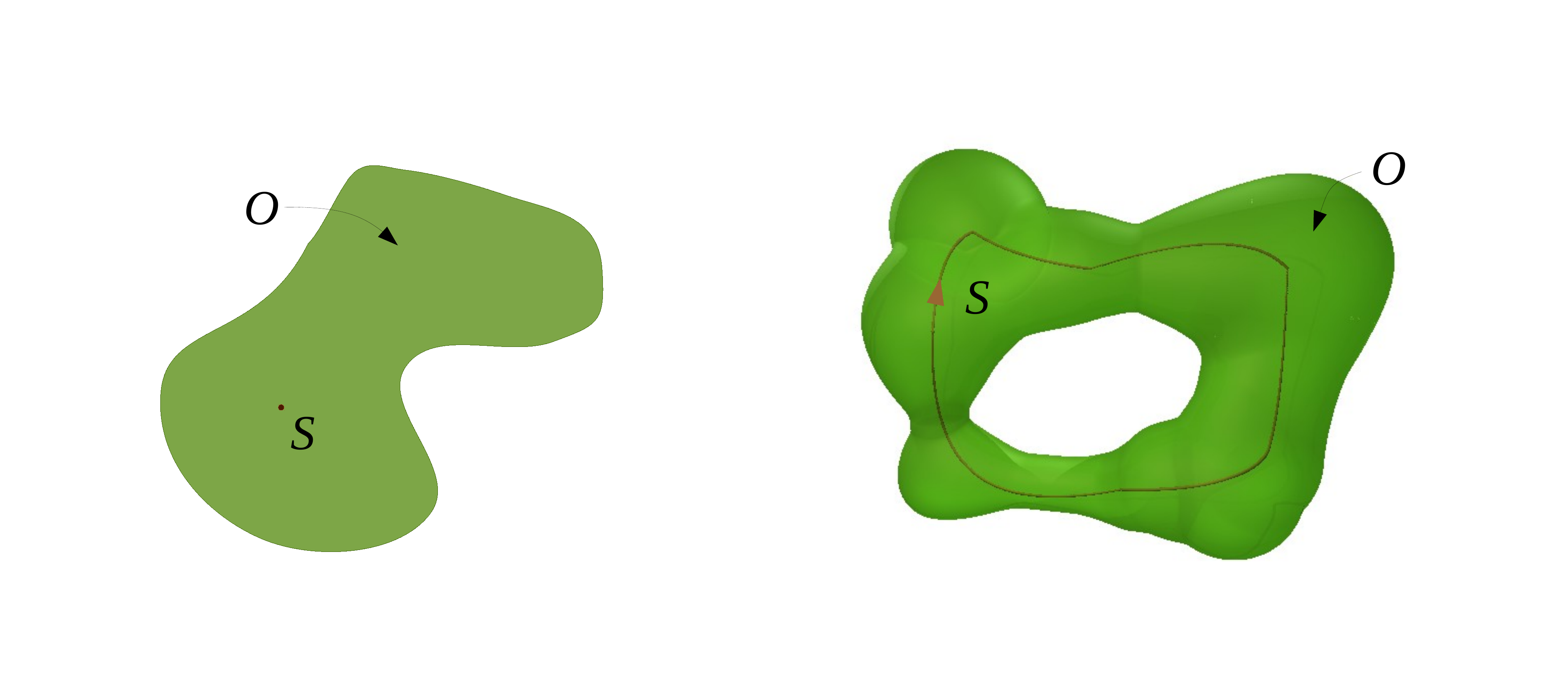}
  \end{center}
  \caption[Obstacles, $O$, can be replaced by equivalents, $S$, without change to $H_{N-1}$ of the complement.]{Obstacles, $O$, can be replaced by equivalents, $S$, without change to $H_{N-1}$ of the complement.
%  In either of these figures, $S$ is a deformation retract of $O$. The justification for this construction is a consequence of Corollary~\ref{cor:deformation-retract-eqv}.
\label{fig:obstacle-eqvs}} 
\end{figure*}

\begin{figure*}
  \centering
  \subfigure%[The hollow torus can be replaced by the image of its generating $1$-cycles, $S$. This replacement does not alter the $(N-1)^{th}$ homology groups of the complement space, neither does it alter the class of a $(N-1)$-cycle in the complement space, $(\mathbb{R}^D-O)$. $O$ and $S$ satisfy the conditions of Proposition~\ref{prop:obstacle-equivalent}.]
  {
%      \label{fig:obstacle-eqvs-genrators-a}
      \includegraphics[width=0.4\textwidth, trim=95 50 80 50, clip=true]{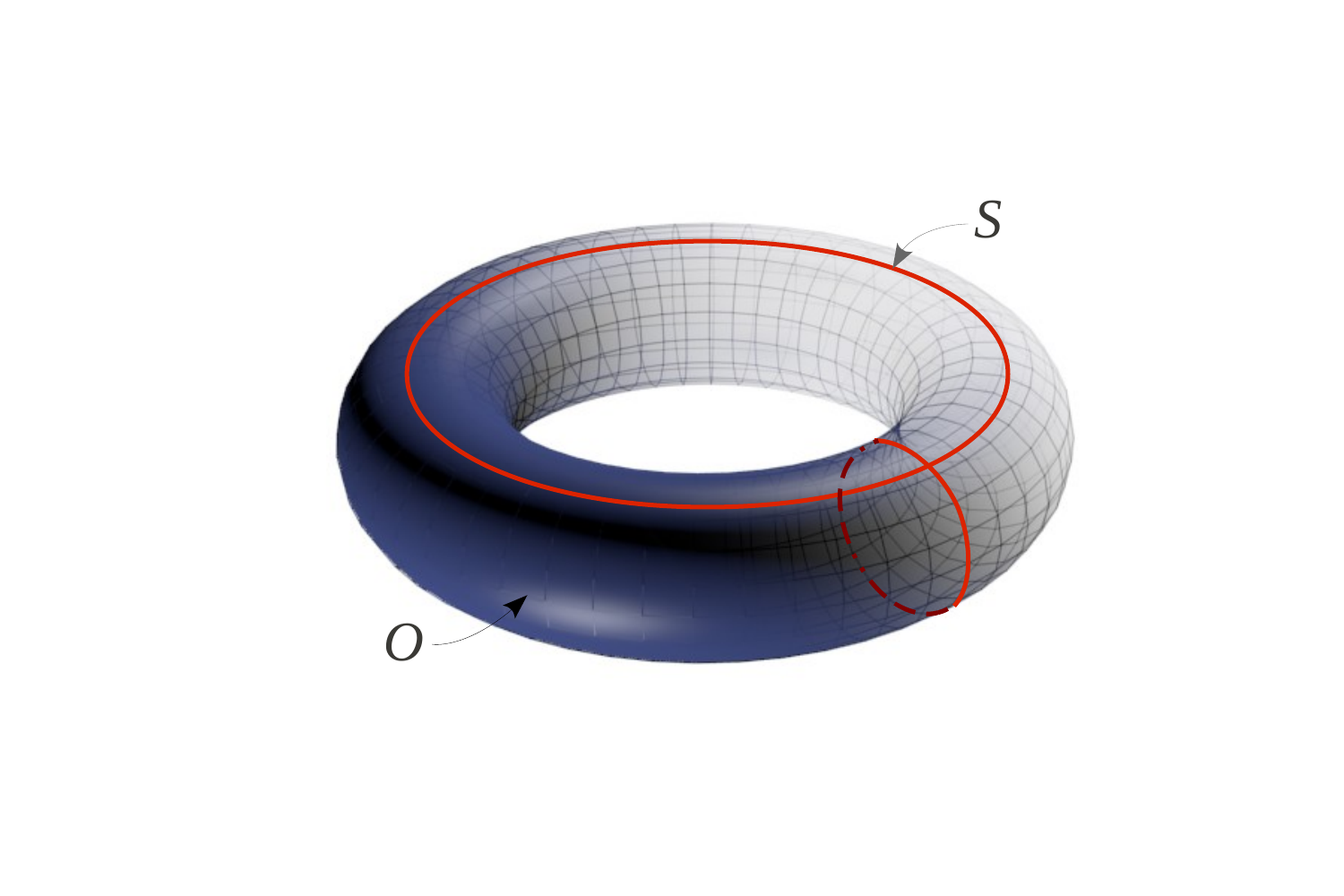}
    } \hspace{0.1in}
  \subfigure%[In this figure we choose a $S\subset O$ such that $H_{D-N}(S) \approxeq H_{D-N}(O)$. The replacement of $O$ by $S$, as in (a), does not alter the $(N-1)^{th}$ homology groups of the complement space, $(\mathbb{R}^D-O)$. However, the replacement does not preserve the class of a $(N-1)$-cycle in the complement space. This is because $H_{D-N}(O,S)\approxeq\!\!\!\!\!\!/ ~~0$ in this case. Thus this is \textbf{not} a valid replacement of $O$.]
  {
%      \label{fig:obstacle-eqvs-genrators-b}
      \includegraphics[width=0.4\textwidth, trim=95 50 80 50, clip=true]{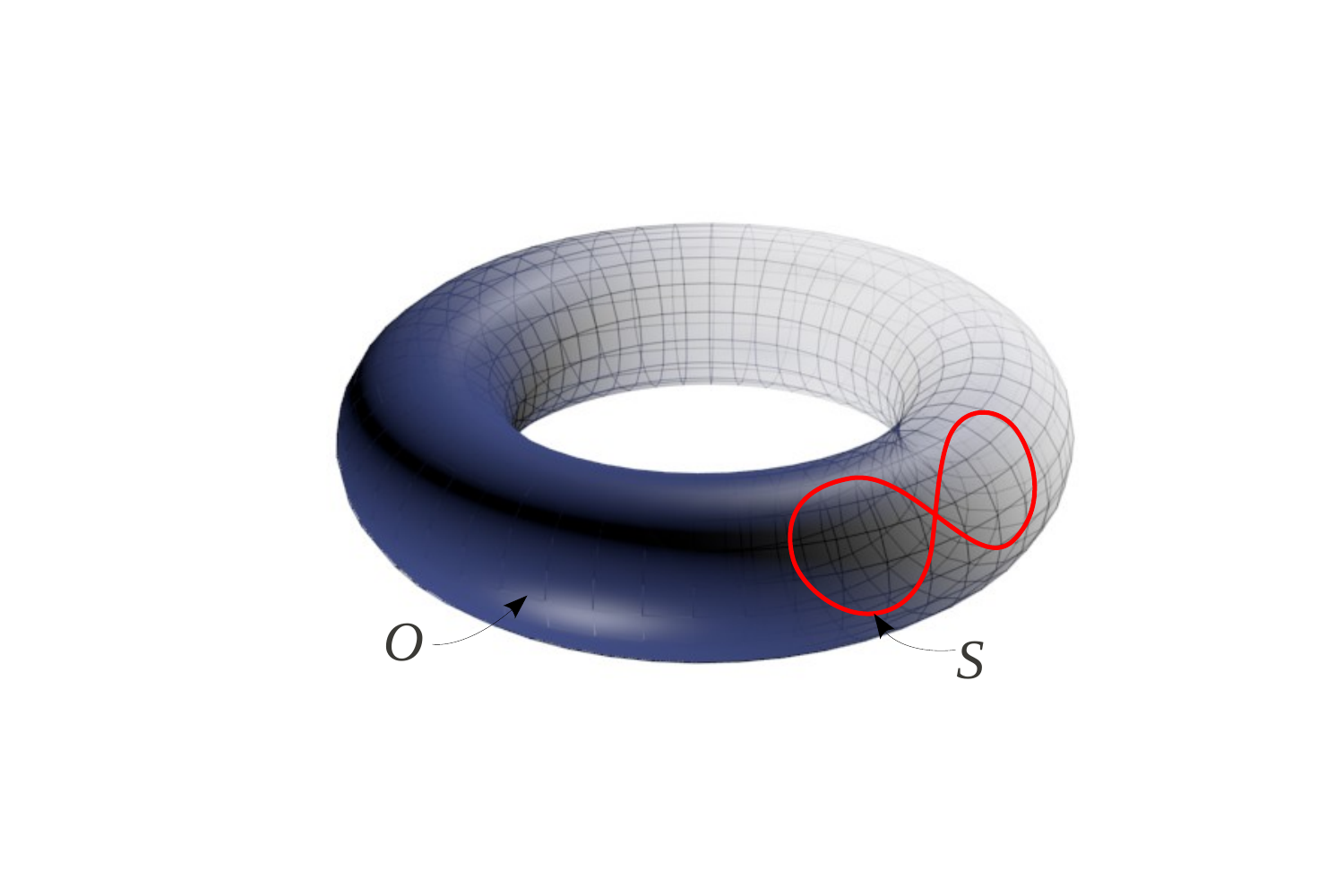}
      }
  \caption[A hollow (or thickened) torus as an obstacle in a $D=3$ dimensional space, with $N=2$ for the problem of robot path planning, with valid (left) and invalid (right) choices for equivalents.]{A hollow (or thickened) torus as an obstacle in a $D=3$ dimensional space, with $N=2$ for the problem of robot path planning (\emph{i.e.} we are interested in homology classes of $(N-1)=1$-dimensional manifolds, which are closed trajectories). It does not have a $(D-N)=1$-dimensional deformation retract or homotopy equivalent. However, we can replace it by its generating $1$-cycles (left). Other choices (right) are invalid, when $H_{D-N}(O,S)\approxeq\!\!\!\!\!\!/ ~~0$.
  %This is the consequence of Corollary~\ref{cor:generating-cycle-eqv}.
  } 
  \label{fig:obstacle-eqvs-genrators}
\end{figure*}

As preparation for the technical details involving linking numbers, we consider the replacement of our obstacles with their $(D-N)$-dimensional representatives. This is trivial for contractible obstacles in the plane (point representatives) and in $3$-dimensional space (cf. the \emph{skeletons} of \cite{planning:AURO:12}). The intuition is that replacing obstacles by their homotopy equivalents leaves the homology classes of trajectories in the complement unchanged (Figure~\ref{fig:obstacle-eqvs}); however, we have dimension constraints, and there exist simple obstacles that do not have a $(D-N)$-dimensional deformation retract (\emph{e.g.} for the $D=3, N=2$ case, a hollow torus does not have a $D-N=1$ dimensional homotopy equivalent). We therefore turn to $(D-N)$-dimensional equivalents faithful to homology in the desired dimension (Figure~\ref{fig:obstacle-eqvs-genrators}).

In the proposition and related corollaries that follow, we represent the ambient configuration space (without obstacles) by $\mathbb{R}^D$, an obstacle by $O$, and $S$ the \changedA{$(D-N)$-dimensional} equivalent of the obstacle with which we replace $O$ for computational simplicity.

% ---------------------------------

\begin{proposition} \label{prop:obstacle-equivalent}
Let $O$ be \changedA{a} compact, locally contractible \changedA{subspace} of $\mathbb{R}^D$. Let $S$ be \changedA{a}
%$(D-N)$-dimensional
compact, locally contractible \changedA{subspace} of $O$, such that the inclusion $i \colon S \hookrightarrow O$ induces
an isomorphism $i_* \colon H_{D-N}(S) \to H_{D-N}(O)$.  Then the inclusion map $\overline{i}:(\mathbb{R}^D-O) \hookrightarrow (\mathbb{R}^D-S)$ induces an isomorphism $\overline{i}_{*}:H_{N-1}(\mathbb{R}^D-O) \rightarrow H_{N-1}(\mathbb{R}^D-S)$.
\end{proposition}

\begin{quoteproof}
%\lo{One consequence of $O$ and $S$ being orientable, and our coefficients being in field $\mathbb{R}$, is that the homology groups of the spaces are freely and finitely generated (see pp. 198 of \cite{Hatcher:AlgTop}).

\vspace{0.1in}
% DAVE'S UPDATE
  Consider the following diagram.
  \[ \xymatrix{
    H^{D-N}(O) \ar^-{f}[r] \ar^{i^*}[d] &
    H_N(\mathbb{R}^D, \mathbb{R}^D - O) \ar^-{\partial}[r] \ar^{\overline{i}_*}[d] &
    H_{N-1}(\mathbb{R}^D - O) \ar^{\overline{i}_*}[d] \\
    H^{D-N}(S) \ar^-{f}[r] &
    H_N(\mathbb{R}^D, \mathbb{R}^D - S) \ar^-{\partial}[r] &
    H_{N-1}(\mathbb{R}^D - S)
  }\]
  The vertical arrows are induced by the inclusions $i$ and $\overline{i}$.  The arrows
  labeled $f$ are the isomorphisms given by proposition 3.46 of \cite{Hatcher:AlgTop} (it is here
  that we use the hypotheses that $O$ and $S$ be compact and locally contractible).  The
  arrows labeled $\partial$ are the boundary homomorphisms in the long exact sequence for the
  pairs $(\mathbb{R}^D, \mathbb{R}^D-O)$ and $(\mathbb{R}^D, \mathbb{R}^D-S)$.  These are
  also isomorphisms, by the contractibility of $\mathbb{R}^D$.

  The square on the right commutes by the naturality of the long exact
  sequence.  The square on the left commutes as well, and while this
  is not explicitly stated in \cite{Hatcher:AlgTop}, it follows easIly
  from the proof of Proposition 3.46, {\em ibid.}.

  The vertical arrow on the left is an isomorphism by hypothesis
  (using the Universal Coefficient Theorem over $\mathbb{R}$), and all
  the horizontal arrows are isomorphisms, so the vertical arrow on the
  right must also be an isomorphism.

\end{quoteproof}

\begin{figure*}
  \centering
  \subfigure[Both $S_1$ and $S_2$ are subsets of the solid torus, $O$. Moreover, each has the homotopy type of the solid torus. $\omega$ is a non-trivial cycle in $(\mathbb{R}^3-O)$.]{
      \label{fig:obstacle-eqv-noteqv-a}
      \includegraphics[width=0.3\textwidth, trim=0 0 0 0, clip=true]{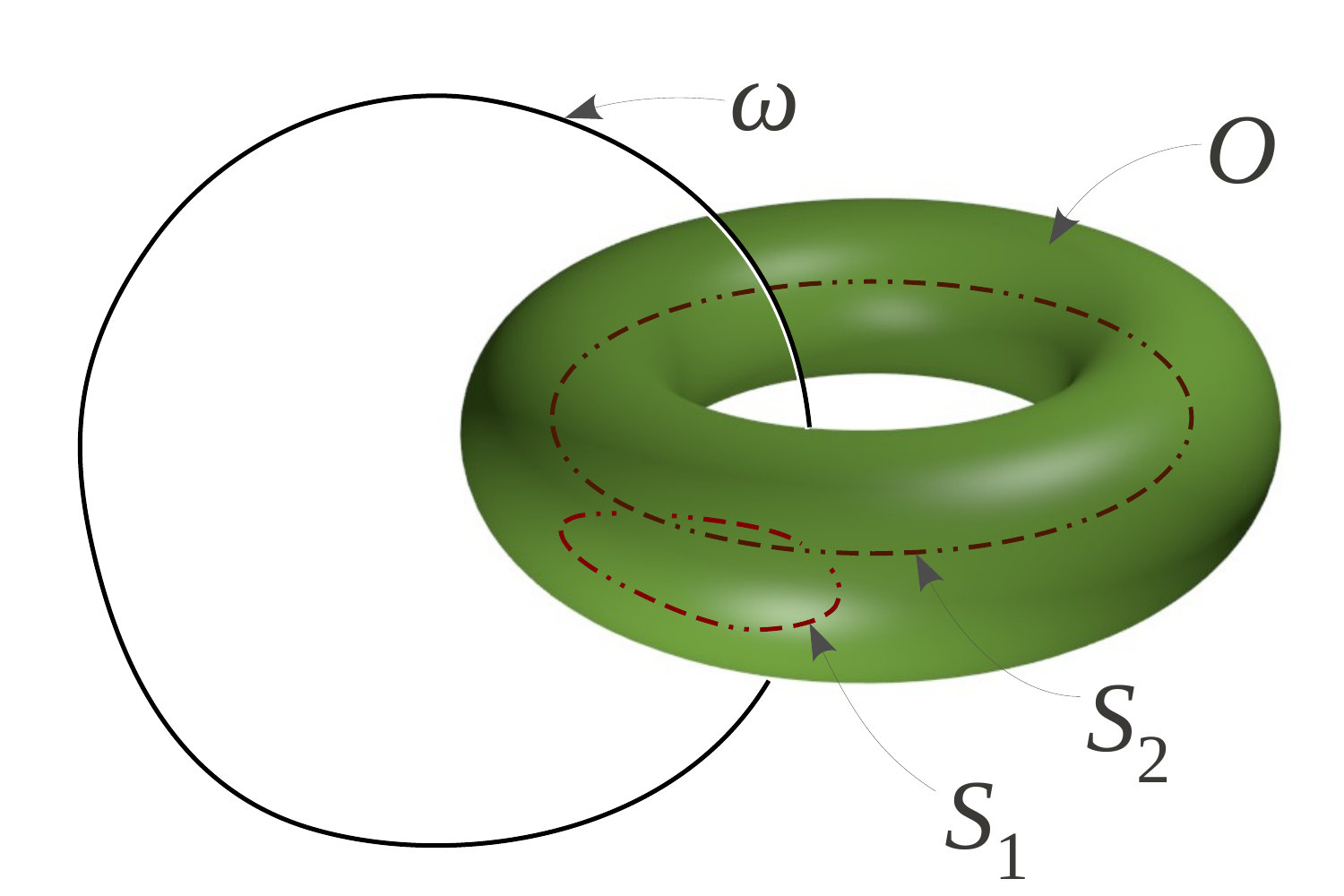}
      } \hspace{0.05in}
  \subfigure[$(\mathbb{R}^3-S_1)$ has homology groups isomorphic to those of $(\mathbb{R}^3-O)$. However, the cycle $\omega$ becomes trivial in $(\mathbb{R}^3-S_1)$. Thus $S_1$ is not a valid replacement of $O$.]{
      \label{fig:obstacle-eqv-noteqv-b}
      \includegraphics[width=0.3\textwidth, trim=0 0 0 0, clip=true]{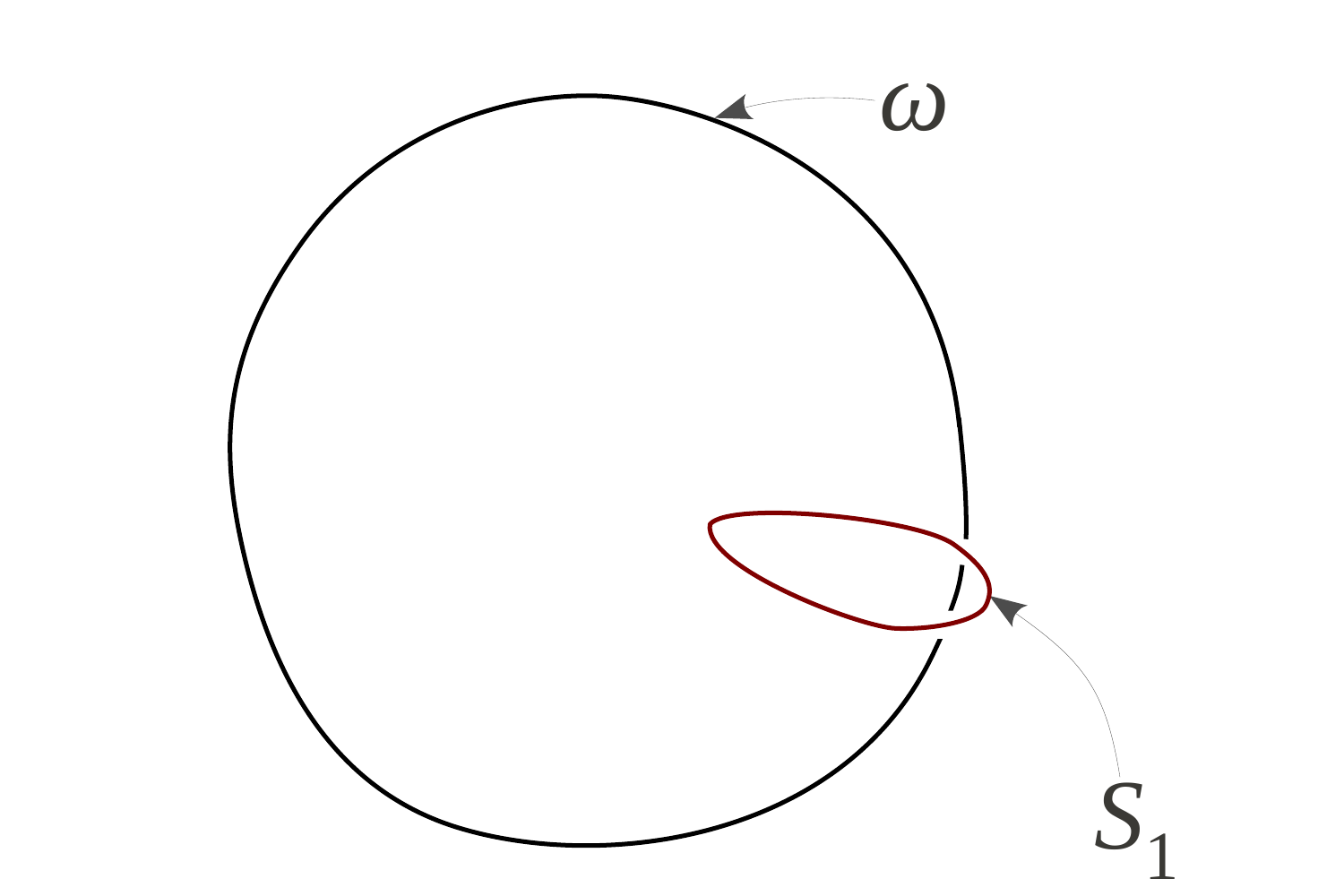}
      } \hspace{0.05in}
  \subfigure[$(\mathbb{R}^3-S_2)$ also has homology groups isomorphic to those of $(\mathbb{R}^3-O)$. Moreover, the cycle $\omega$ remain non-trivial in $(\mathbb{R}^3-S_2)$. $S_2$ is a valid replacement of $O$.]{
      \label{fig:obstacle-eqv-noteqv-c}
      \includegraphics[width=0.3\textwidth, trim=0 0 0 0, clip=true]{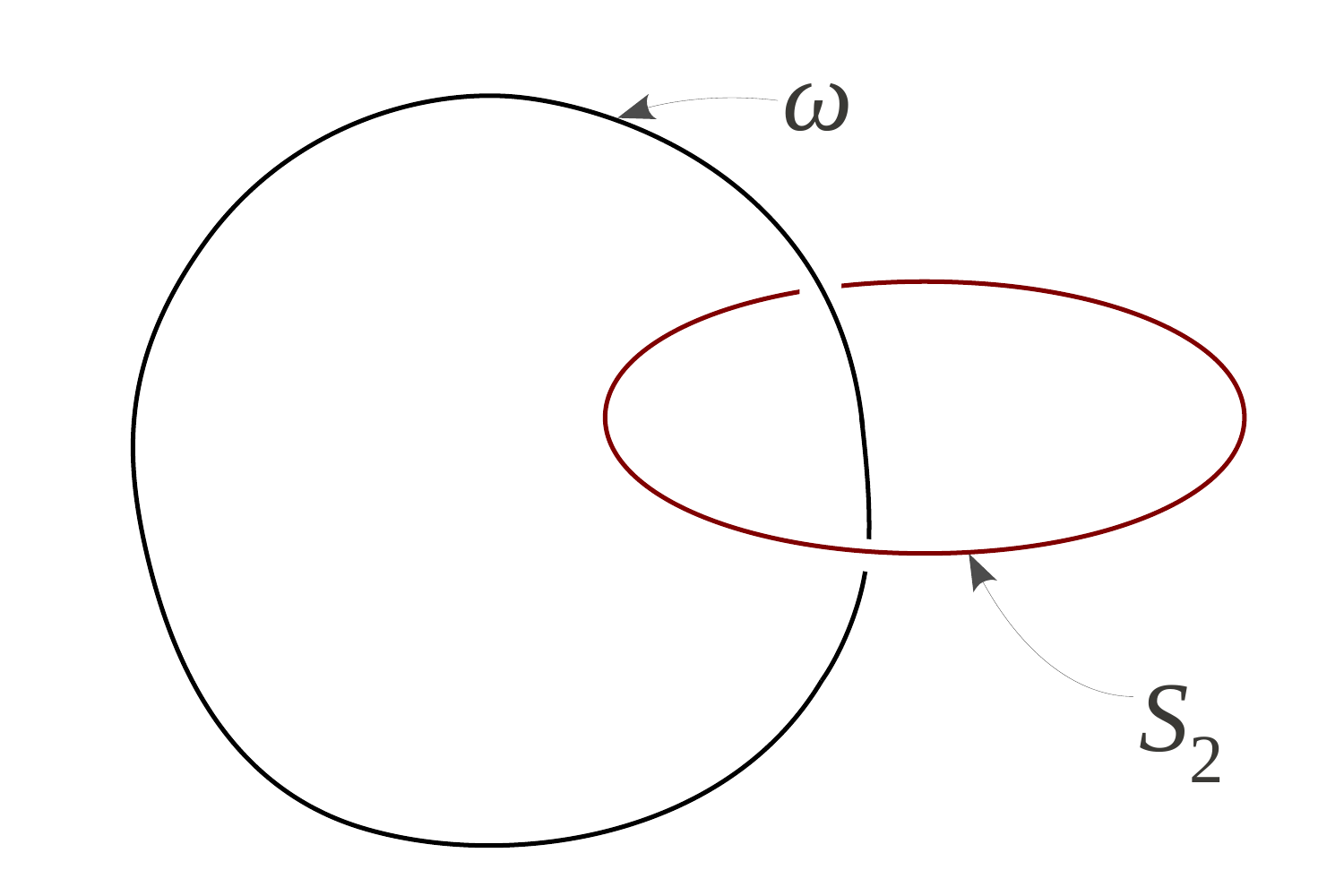}
      }
  \caption{A solid torus [left] with valid [right] and invalid [middle] equivalents. This is an example with $D=3,N=2$. The replacement needs to be such that the inclusion map $\overline{i}:(\mathbb{R}^D-O) \hookrightarrow (\mathbb{R}^D-S)$ induces the isomorphism.}
  \label{fig:obstacle-eqv-noteqv} 
\end{figure*}

In light of robot path planning, $O$ in the above proposition is a solid obstacle in the environment, and $S$ is its equivalent/replacement (in the terminology of \cite{planning:AURO:12} these are \emph{representative points} of obstacles on a $2$-dimensional plane, and \emph{skeletons} of obstacles in a $3$-dimensional Euclidean space).
The aim of the above proposition is to establish a relationship between the homology groups of the complement (or free) spaces, $(\mathbb{R}^D-O)$ and $(\mathbb{R}^D-S)$, from some known relationship between the spaces $O$ and $S$. %In fact, it is not just the homology groups that we are trying to establish relationship for, but homology classes of $(N-1)$-dimensional manifolds (the \emph{closed trajectories} in robot planning problem) in the complement space.
\lo{
For example, in Figure~\ref{fig:obstacle-eqv-noteqv}, the solid torus, which represents an obstacle in $\mathbb{R}^3$, can be replaced by $S_2$ for computations of homology class of closed loops like $\omega$ (formed by pair of trajectories). Under such a replacement, the homology class of $\omega$ in the ambient space remains unchanged. This replacement, according to Proposition~\ref{prop:obstacle-equivalent}, is justified by the facts that $H_1(S_2)\approxeq H_1(O)$ and $H_1(O,S_2)\approxeq 0$ (due to the fact that $S_2$ is a deformation retract of $O$).
Corollary~\ref{cor:deformation-retract-eqv} simply asserts that a deformation retract (which $S_2$ indeed is of $O$) in fact satisfies the required conditions of the proposition. On the other hand, although $S_1$ satisfies $H_1(S_1)\approxeq H_1(O)$ it does not satisfy $H_1(O,S_1)\approxeq 0$. Thus $S_1$ is not a valid replacement of $O$.}
\changedA{In the corollaries below, we suggest couple of approaches for identifying valid replacements, $S$, of a given obstacles, $O$.}

% ---------------------------------

The following is trivial, but stated formally for future reference.
\begin{corollary} \label{cor:deformation-retract-eqv}
If $S$ and $O$ are compact, locally contractible \changedA{subspaces} of $\mathbb{R}^D$ such that $S$ is a deformation retract of $O$, then the inclusion map $\overline{i}:(\mathbb{R}^D-O) \hookrightarrow (\mathbb{R}^D-S)$ induces isomorphisms $\overline{i}_{*}:H_{*}(\mathbb{R}^D-O) \rightarrow H_{*}(\mathbb{R}^D-S)$
\end{corollary}
\lo{
\begin{quoteproof}
% LOOK, IT'S HOMOTOPY EQUIVALENT IF YOU HAVE A DEFRET, RIGHT? NO NEED FOR L.E.S. ARGUMENTS I THINK...
If $S$ is a deformation retract of $O$ then the inclusion $S \hookrightarrow O$
induces isomorphisms $H_n(S) \approxeq H_n(O)$
for all $n$.  %, so that $H_n(O,S) \approxeq 0$ for all $n$ by a long exact sequence argument.
Thus the conditions of Proposition~\ref{prop:obstacle-equivalent} are automatically
satisfied, and so the result follows. 
%Deformation retraction implies homotopy equivalence, and that the inclusion induces isomorphisms $i_{*}:H_{D-N}(S)\xrightarrow{\approxeq} H_{D-N}(O)$ for all $N$. Using this in the long exact sequence of homology groups for pair $(O,S)$ we further have $H_{D-N}(O,S)\approxeq 0$ for all $N$ (This can be concluded by observing that $\overline{i}_{*}$ being an isomorphism requires that for exactness of the sequence, $\partial_{*}$ be a zero map, and $j_{*}$ be a surjection and a zero map at the same time. Thus $H_{D-N}(O,S)$ requires to be zero for all $N$.).
 %\changedA{Moreover, $(\mathbb{R}^D-S)$ is a manifold since $S$ is a manifold.}
% The result can then be concluded using Proposition~\ref{prop:obstacle-equivalent}.
\end{quoteproof}
}

\begin{corollary} \label{cor:generating-cycle-eqv}
% Given a compact, locally contractible and orientable \changedA{subspace} $O \subset \mathbb{R}^D$, \changedA{we choose a generating set of $(D-N)$-cycles in $Z_{D-N}(O)$, and denote them by $\overline{S}_k, k=1,2,\cdots,m$ (i.e., the homology classes of $\overline{S}_k, k=1,2,\cdots,m$, together forms a basis for $H_{D-N}(O)\approxeq \mathbb{R}^m$, which is possible since $H_{D-N}(O)$ is freely generated). }
Let $O \subset \mathbb{R}^D$ be compact and locally contractible.
Suppose there exists a set of pairwise-disjoint, connected, closed,
oriented $(D-N)$-dimensional manifolds $S_k \subseteq O$, $k =
1, \ldots, m$, such that the fundamental classes $[S_1], \ldots,
[S_m]$~form a basis for the homology group $H_{D-N}(O)$.
%\changedA{We assume that $O$ is such that it is possible to make the choice such that the images of each of the cycles, denoted as $S_k, k=1,2,\cdots,m$, are connected $(D-N)$-dimensional closed manifolds that are compact, locally contractible and orientable. This implies $H_{D-N}(S_k) \approxeq \mathbb{R}$.}
% We define \changedA{$\widetilde{\mathcal{S}} = \bigsqcup_{k=1,2,\cdots,m} S_k$.}
%let $\{S_k\}_{k=1,2,\cdots,m}$ be the images of generating cycles of $H_{D-N}(O)$
%(i.e. if $\{\overline{S}_k\}$ are $(D-N)$-cycles in $O$ such that the homology classes of $\{\overline{S}_k\}$ generate the group $H_{D-N}(O)$ freely, then $\{S_k\}$ are the images of $\{\overline{S}_k\}$), and let $\widetilde{\mathcal{S}} = \bigcup_{k=1,2,\cdots,m} S_k$.
Let $\widetilde{\mathcal{S}} = \bigcup_{k=1}^{m} S_k$.  Then the
inclusion map $\overline{i}:(\mathbb{R}^D-O) \hookrightarrow
(\mathbb{R}^D-\widetilde{\mathcal{S}})$ induces an isomorphism
$\overline{i}_{*:{\scriptscriptstyle
N-1}}:H_{N-1}(\mathbb{R}^D-O) \rightarrow
H_{N-1}(\mathbb{R}^D-\widetilde{\mathcal{S}})$.
\end{corollary}
\begin{quoteproof}
{
By construction, the inclusion induces an isomorphism
$H_{D-N}(\widetilde{\mathcal{S}}) \to H_{D-N}(O)$, 
and so the result follows from Proposition~\ref{prop:obstacle-equivalent}.
}
\end{quoteproof}

The consequence of the last two corollaries is that instead of computing homology classes of $(N-1)$ cycles in the original punctured space $(X-O)$, we can replace the obstacles $O$ with equivalents $S$
% and compute their homology in $(X-S)$, and yet, the results we obtain will be identical.
while preserving the relevant homology (cf. \cite{planning:AURO:12} for special cases).
%Corollary~\ref{cor:deformation-retract-eqv} says that we can replace obstacles with $(D-N)$-dimensional deformation retracts (Figure~\ref{fig:obstacle-eqvs}).
In cases where $(D-N)$-dimensional deformation retracts do not exist (e.g., Figure~\ref{fig:obstacle-eqvs-genrators}), Corollary~\ref{cor:generating-cycle-eqv} allows one to replace obstacles by $(D-N)$-dimensional equivalents (generating cycles of $(D-N)^{th}$ homology group).

\subsection{Reduced Problem Definition} \label{sec:reduced-problem-algtop}

Thus we have established that obstacles $\widetilde{\mathcal{O}}\subset\mathbb{R}^D$ (which represent illegal zones in robot planning problems) may be replaced by equivalents $\widetilde{\mathcal{S}}$ preserving the appropriate homology. We may (and do) choose the equivalents $\widetilde{\mathcal{S}}$ to be \changedA{a disjoint union} of \changedA{connected, closed, orientable $(D-N)$-dimensional  manifolds}. The reduced problem definition follows:

\begingroup \leftskip2em
\rightskip \leftskip \vspace{0.05in}
{\bf Given:} (1) the \emph{singularity manifolds} --- a disjoint collection $\widetilde{\mathcal{S}}=S_1\sqcup S_2\sqcup\cdots\sqcup S_m$ of $(D-N)$-dimensional ($N > 1$), connected, closed, orientable submanifolds, of $\mathbb{R}^D$; and (2) the \emph{candidate manifolds}) --- a collection of $(N-1)$-dimensional, closed, orientable manifolds in $(\mathbb{R}^D \setminus \widetilde{\mathcal{S}})$. \\
\indens {\bf Problem:} identify the homology classes of the candidate manifolds in the complement of the singularity manifolds.
Specifically, design a complete set of easily-computed invariants for these homology classes by finding a set of explicit generators for $H^{N-1}(\mathbb{R}^D \setminus \widetilde{\mathcal{S}})$ and integrating these generators over candidate manifolds.

\par
\endgroup

\vspace{0.1in}
In order to compute the action of the cocycles on the candidate manifolds, we represent them as $(N-1)$-cycles (\emph{i.e.} top-dimensional covering cycles). %Those, by inclusion, are $(N-1)$-cycles in $(\mathbb{R}^D \setminus \widetilde{\mathcal{S}})$.
Thus, given a candidate manifold $\omega$, we can use a \emph{cellular cover} of the manifold, $\overline{\omega}$, which is also an $(N-1)$-cycle in $(\mathbb{R}^D \setminus \widetilde{\mathcal{S}})$ under the inclusion map $\omega \hookrightarrow (\mathbb{R}^D \setminus \widetilde{\mathcal{S}})$ (a map that we will assume implicitly most often).
However, given two cycles $\overline{\omega}_1, \overline{\omega}_2 \in Z_{N-1}(\mathbb{R}^D \setminus \widetilde{\mathcal{S}})$, instead of checking if or not $\overline{\omega}_1 - \overline{\omega}_2$ is boundary in $H_{N-1}(\mathbb{R}^D \setminus \widetilde{\mathcal{S}})$, we will compute complete invariants $\phi_{\widetilde{\mathcal{S}}}(\overline{\omega}_1)$ and $\phi_{\widetilde{\mathcal{S}}}(\overline{\omega}_2)$, comparing them to make the desired assertion. In particular, we construct the function $\phi_{\widetilde{\mathcal{S}}}(\cdot)$ to be in form of an integration over $\overline{\omega}$ of some set of differential $(N-1)$-forms. Our strategy --- using integration and differential forms --- is a traditional method for understanding (co)homology of manifolds and submanifolds \cite{bott1982differential}.

% ---------------------------------------------

\section{Preliminaries on Linking Numbers} \label{sec:linking-num-pre}

Equipped with the notion of the $(D-N)$-dimensional replacements of the obstacles/punctures, $S_i$, we proceed towards computing the homology classes of $(N-1)$-cycles (in light of robot planning problem those are the closed trajectories) of $(\mathbb{R}^D \setminus \widetilde{\mathcal{S}})$. In this section we recall various notions of \emph{intersection} and \emph{linking number}, and from this:
\begin{itemize}
 \item[i.] Infer homology classes of the $(N-1)$-cycles in $(\mathbb{R}^D \setminus S_i)$ from linking data (Proposition~\ref{prop:linking-number-homology-injectivity}),
 \item[ii.] Computing the linking number using an integration over the $(N-1)$-cycle and a top-dimensional cycle of the $S_i$ (Proposition~\ref{prop:linking-number-integration}).
\end{itemize}
We illustrate the ideas using examples from robot planning problems.

\subsection{Definitions}

\begin{figure}
\centering
 \includegraphics[width=0.75\textwidth, trim=0 75 0 60, clip=true]{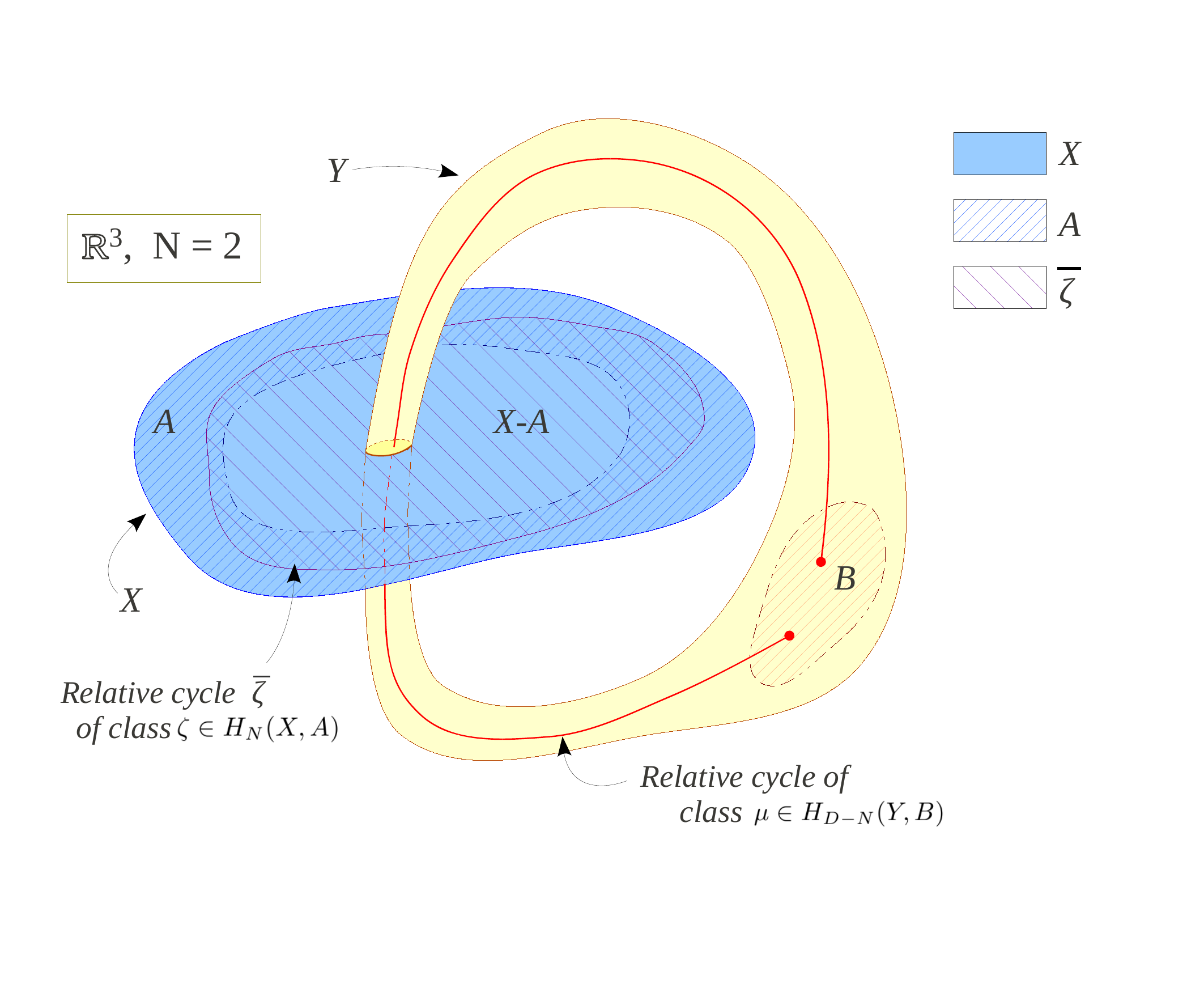}
  \caption{Illustration of intersection number in $\mathbb{R}^3$ with $N=2$ in light of Definition~\ref{def:intersection-number}.} 
\label{fig:intersection-number}
\end{figure}

Recall the definition of \emph{intersection number}:

\vspace{0.1in}
\begin{definition} [\textbf{\textit{Intersection Number}} -- Ch. VII, Def. 4.1 of \cite{dold1995lectures}] \label{def:intersection-number} ~
  Suppose $X$ and $Y$ are sub-manifolds of $\mathbb{R}^D$, and
  $A\subset X \subset \mathbb{R}^D$, $B\subset Y \subset \mathbb{R}^D$
  are such that $A \cap Y =\emptyset, X \cap B = \emptyset$
  (Figure~\ref{fig:intersection-number}).  Consider the map $p: (X
  \times Y, A \times Y \cup X \times B) \rightarrow (\mathbb{R}^D,
  \mathbb{R}^D - \{0\})$ given by $p(x,y) = x-y$.  The composition 
\[
 H_{N}(X,A) \times H_{D-N}(Y,B) \xrightarrow{\quad\times\quad} H_{D}(X\times Y, A\times Y \cup X \times B) \xrightarrow{(-1)^{D-N} p_{*}} H_{D}(\mathbb{R}^D,\mathbb{R}^D-\{0\}) \vspace{-0.05in}
\]
is called the \emph{intersection pairing} (where '$\times$' denotes the homology cross product -- see p. 268 of \cite{Hatcher:AlgTop}).
%(Note that the product `$\times$' for homology groups is the \emph{homology cross product}, which is more closely related to the tensor product of the homology groups rather than the cartesian product -- see p. 268 of \cite{Hatcher:AlgTop}).
We write 
\[
 \mathscr{I}(\zeta,\mu) = (-1)^{D-N} p_{*} (\zeta \times \mu), ~~~\textrm{for } \zeta \in H_{N}(X,A), ~\mu \in H_{D-N}(Y,B) \vspace{-0.05in}
\]
and call this element of $H_{D}(\mathbb{R}^D,\mathbb{R}^D-\{0\}) \approxeq \mathbb{R}$ the \emph{intersection number} of $\zeta$ and $\mu$.
\end{definition}

\vspace{0.1in}
\begin{definition} [\textbf{\textit{Linking Number}} -- Adapted from Ch. 10, Art. 77 of \cite{seifert1980seifert}] \label{def:linking-number} ~
We borrow definitions of $X,A,Y$ and $B$ from Definition~\ref{def:intersection-number}.
Recall from the \emph{long exact sequence} of the pair $(X,A)$ the \emph{connecting homomorphism} $\partial_{*}: H_N(X,A) \rightarrow H_{N-1}(A)$.
If $\varsigma \in H_{N-1}(A)$ is such that it can be written as $\varsigma = \partial_{*} \zeta$ for some  $\zeta \in H_{N}(X,A)$, and if $\mu \in H_{D-N}(Y,B)$, then the \emph{linking number} between $\varsigma$ and $\mu$ is defined as $\mathscr{L}(\varsigma,\mu) = \mathscr{I}(\zeta,\mu)$.
%
% \vspace{0.1in}
%\noindent\textbf{Note} that similar to the intersection number, %by the functoriality of homology,
%linking number can be defined between a cycle, $\overline{\varsigma}$, in $A$ and a relative cycle, $\overline{\mu}$, in $(Y,B)$.
\end{definition}

\begin{figure}[t]
\centering
 \includegraphics[width=0.85\textwidth, trim=0 30 0 0, clip=true]{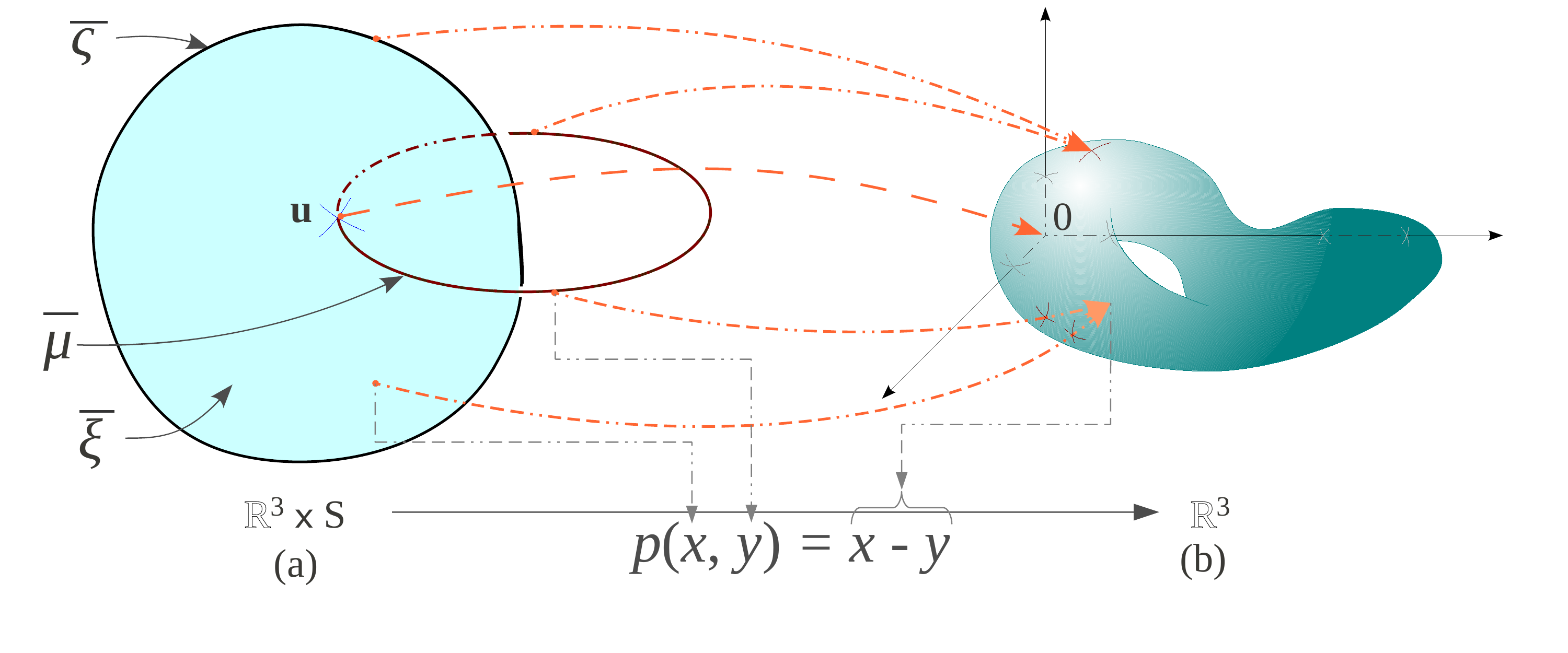} 
 \caption[A simplified illustration of intersection number and linking number in $\mathbb{R}^3$ with $N=2$.]{A simplified illustration (following from Figure~\ref{fig:obstacle-eqv-noteqv-c}) of intersection number and linking number in $\mathbb{R}^3$ with $N=2$. This is a special case of Definition~\ref{def:intersection-number} when $X = \mathbb{R}^3, ~A = \mathbb{R}^3 - S, ~~Y = S ~\textrm{ and } B=\emptyset$. \newline \emph{Figure (a) on the left}: The intersection number is computed between a $N$-chain, $\overline{\xi}$ (more precisely it is a relative cycle in $(X,A)$ that we consider -- the boundary of $\overline{\xi}$ trivialized), and the $(D-N)$-cycle, $\overline{\mu}$, that is a top-dimensional cycle on $S$. In this figure the said intersection number is $\pm 1$ due to the single intersection marked by the \lo{blue} `cross' at $\mathbf{u}$. Then, by definition, that is equal to the linking number between $\overline{\varsigma} = \partial \overline{\xi}$ and $\overline{\mu}$. \newline \emph{Figure (b) on the right}: The precise definition requires a mapping, $p$, from pair of points in the original space (one point from the $2$-chain, $\overline{\xi}$, embedded in the ambient space, $\mathbb{R}^3$, and another from $S$) to (a different copy of) $\mathbb{R}^3$.
 %$(\mathbb{R}^3-0)$ (\emph{i.e.} $\mathbb{R}^3$ without the origin).
 The intersection/linking number is then, informally, the number of times intersection points in the pre-image of $p$ (points like $\mathbf{u}$) maps to the origin, $0$ (with proper sign), in the image, or equivalently, the number of times
 the image of $\overline{\varsigma}\times\overline{\mu}$, under the action of $p$, \emph{wraps around} the origin. Thus, it is the homology class of the cycle $p(\overline{\varsigma}\times\overline{\mu})$ in the punctured Euclidean space $(\mathbb{R}^D-0)$.}
\label{fig:intersection-linking-number-simple} 
\end{figure}

\ls{
\subsubsection{Simplified Description of the Definitions}

Let us consider the simple case when $X = \mathbb{R}^3, ~A = \mathbb{R}^3 - S, ~~Y = S ~\textrm{ and } B=\emptyset$, and with $D=3,N=2$, which arises in robot path planning in $\mathbb{R}^3$ (Figure~\ref{fig:intersection-linking-number-simple}). Let $\overline{\mu}$ be a top-dimensional cycle on the $(D-N)$-dimensional manifold, $S$ (to be consistent with the notations used in the definition).
Intersection number, as the name suggests, informally speaking, counts the number of intersections between a $N$-chain $\overline{\xi}$ (in light of Definition~\ref{def:intersection-number}, it is represented by the relative cycle $\overline{\zeta}$), and a $(D-N)$-cycle $\overline{\mu}$.
%Consider the case where $D=3, N=2$ (one corresponding to robot trajectories in $\mathbb{R}^3$). If $S$ is the
Thus, in Figure~\ref{fig:intersection-linking-number-simple}, informally, the intersection number between $\overline{\xi}$ and $\overline{\mu}$ is $\pm 1$ (the sign depends on orientation).}
{These definitions, being based on homology classes, of course are applicable to cycle representatives. Figure~\ref{fig:intersection-linking-number-simple} illustrates the intuition behind these definitions using a simple example.}

\lo{
Now, if $\overline{\xi}$ has a boundary, say $\overline{\varsigma} = \partial \overline{\xi}$, the linking number between $\overline{\varsigma}$ and $\overline{\mu}$ is, by definition, the intersection number between $\overline{\xi}$ and $\overline{\mu}$.

In defining the intersection number, however, one does not talk about the $N$-chain $\overline{\xi}$. Instead, one talks about the corresponding relative cycle in $(\mathbb{R}^3,\mathbb{R}^3 - S)$ under the action of the quotient map $j:C_N(\mathbb{R}^3)\rightarrow C_N(\mathbb{R}^3,\mathbb{R}^3 - S)$ -- that is, the part of $\overline{\xi}$ that does not intersect with $S$, is trivialized (which is $\overline{\zeta} ~(=j(\overline{\xi}))$ in notation of Definition~\ref{def:intersection-number}).
This is because homology classes are not defined for chains, rather can be defined for cycles or relative cycles only. By trivializing the part of the chain that contains the boundary (\emph{i.e.} the one lying in $(\mathbb{R}^3 - S)$), we convert it into a relative cycle, thus enabling us to talk about its homology class ($\zeta$ in Definition~\ref{def:intersection-number}).

For constructing a precise algebraic definition of linking/intersection number, the relative cycles $\overline{\zeta}$ and $\overline{\mu}$ are then mapped to $\mathbb{R}^3$ via the map $p(x,y)=x-y$. The intersection number is then, informally, the number of times the intersection points (like $\mathbf{u}$ in Figure~\ref{fig:intersection-linking-number-simple}) map to (with proper sign) the origin in the co-domain of $p$. Linking number is then essentially the number of times the image of $\overline{\varsigma}\times\overline{\mu}$, under the action of $p$, \emph{wraps around} the origin in the co-domain of $p$. In other words, it is the homology class of the cycle $p(\overline{\varsigma}\times\overline{\mu})$ in the punctured Euclidean space $(\mathbb{R}^D-0)$.
}

% --------------------------------------------------------------------

%\vspace{0.1in}
\subsection{Propositions on Linking Number}

We state and two propositions\ls{, each followed by simplified explanation of the result of the propositions}{~related to linking numbers, and how they relate to homology class of cycles}. The first is well-known but stated for completeness.

\vspace{0.1in}
\begin{proposition} [Uniqueness of linking number] \label{prop:linking-number-independence}
 If $H_N(X)=H_{N-1}(X)=0$ holds, then
 $\mathscr{L}(\varsigma,\mu)$ is independent of the choice of $\zeta$ \changedA{in Definition~\ref{def:linking-number}} \emph{\cite{seifert1980seifert}}. %More precisely, under the said conditions,
% $\partial_{*}^{-1}$ exists,
%\changedA{$\partial_{*}: H_N(X,A) \rightarrow H_{N-1}(A)$ is bijective,}
%and thus $\mathscr{L}(\varsigma,\mu) = \mathscr{I}(\partial_{*}^{-1} \varsigma,\mu) = (-1)^{D-N}p_{*}(\partial_{*}^{-1} \varsigma \times \mu)$.
\end{proposition}
\lo{
\begin{quoteproof}
 From the long exact sequence for the pair $(X,A)$, using the condition $H_N(X)=H_{N-1}(X)=0$, it follows that $\partial_{*}: H_N(X,A) \rightarrow H_{N-1}(A)$ is an isomorphism (See p. 114 of \cite{Hatcher:AlgTop}).  Then the only choice of $\zeta$ is given by $\zeta = \partial_{*}^{-1}\varsigma$, so that $\mathscr{L}(\varsigma,\mu) = \mathscr{I}(\partial_{*}^{-1} \varsigma,\mu) = (-1)^{D-N}p_{*}(\partial_{*}^{-1} \varsigma \times \mu)$.  %Hence the result follows.
\noindent [Note that the contractibility of $X$ is sufficient for the said condition to hold.]
\end{quoteproof}
}
\lo{
\begin{figure*}
  \centering
  \subfigure[On $\mathbb{R}^2$ the linking number between a point $\overline{\mu}$ (a $0$-cycle) and a $1$-cycle $\overline{\varsigma}$ is uniquely determined, and is equal to the intersection number between $\overline{\mu}$ and a $2$-chain $\overline{\xi}$ such that $\overline{\varsigma}=\partial\overline{\xi}$. It can be shown that the choice of $\overline{\xi}$ does not matter.]{
      \label{fig:linking-number-uniqueness-a}
      \includegraphics[width=0.3\textwidth, trim=50 0 80 0, clip=true]{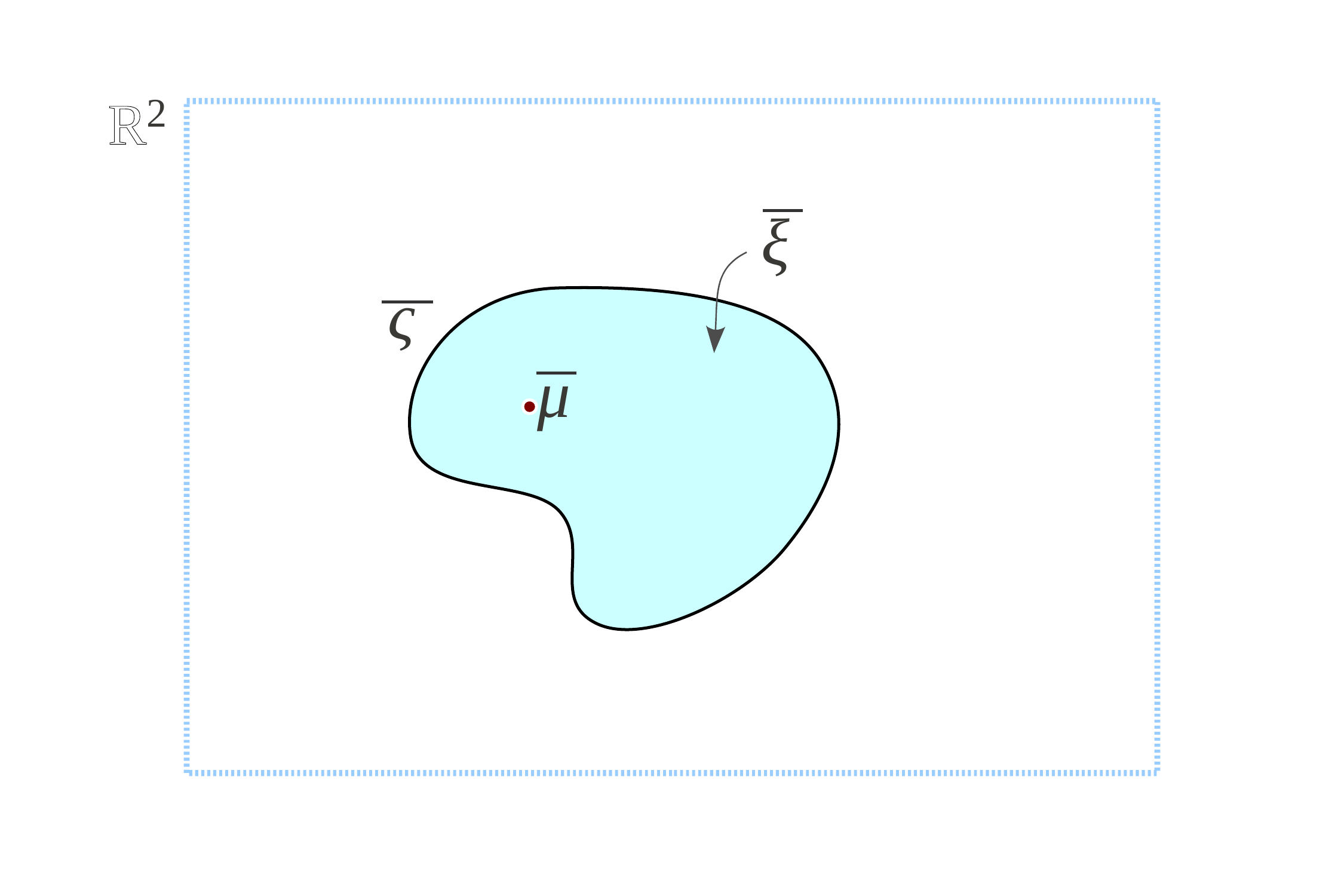}
      } \hspace{0.05in}
  \subfigure[Similarly, in $\mathbb{R}^3$ the linking number between a $1$-cycle, $\overline{\mu}$, and a closed cycle $\overline{\varsigma}$ is uniquely determined, and is equal to the intersection number between $\overline{\mu}$ and a $2$-chain $\overline{\xi}$ such that $\overline{\varsigma}=\partial\overline{\xi}$. It can be shown that the choice of $\overline{\xi}$ does not matter -- that is, we could choose $\overline{\xi}_1$ or $\overline{\xi}_2$ for computation of the intersection number, and the value that we would obtain will be the same.]{
      \label{fig:linking-number-uniqueness-b}
      \includegraphics[width=0.3\textwidth, trim=130 50 130 0, clip=true]{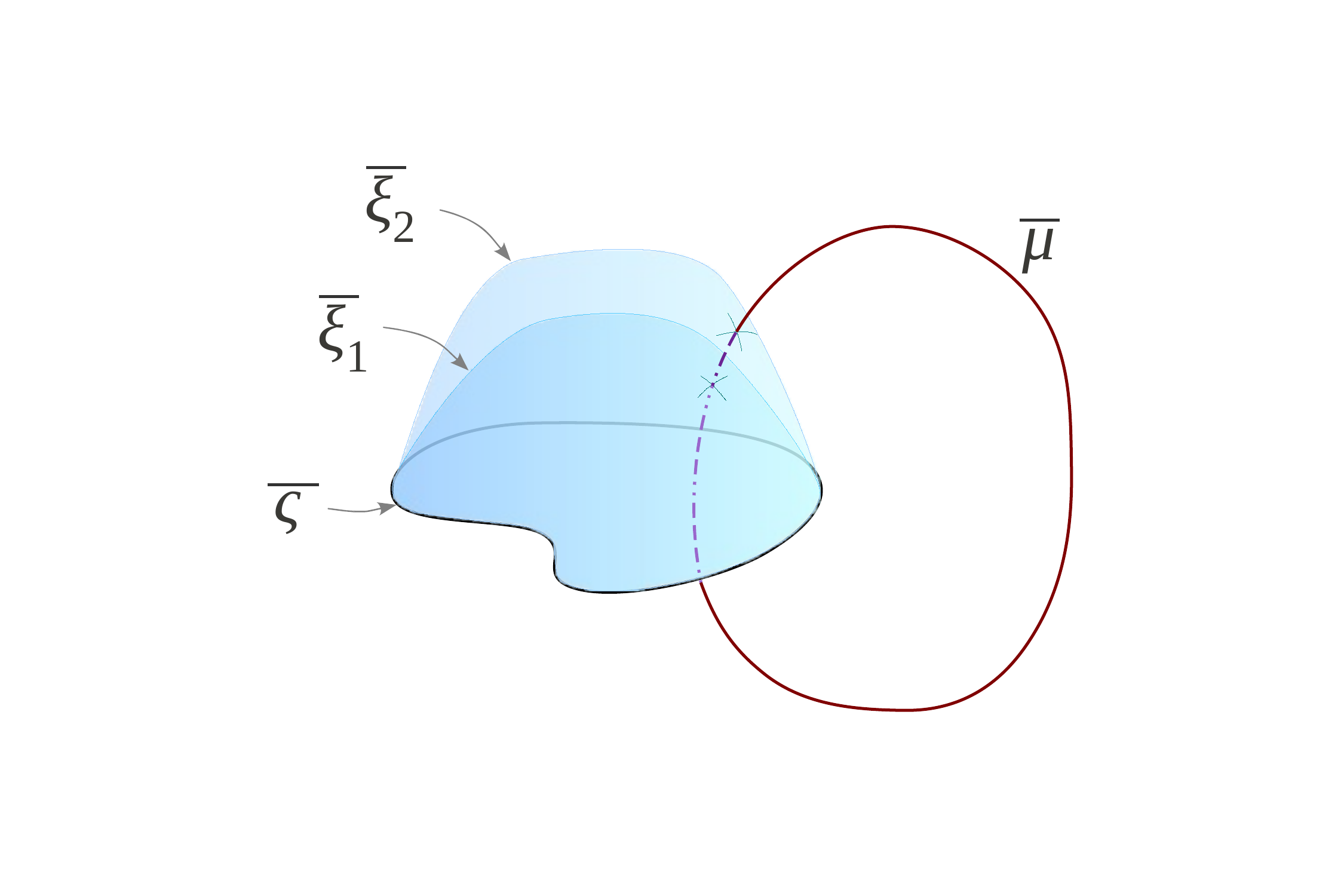}
      } \hspace{0.05in}
  \subfigure[However, if the ambient space, \changedA{$M$}, is not contractible (in this figure it is the $2$-sphere, $\mathbb{S}^2$), then the linking number between $\overline{\mu}$ and $\overline{\varsigma}$ is not unambiguously defined. This is because the choices of $\overline{\xi}$ such that $\overline{\varsigma}=\partial\overline{\xi}$ can be made in ways such that its intersection number with $\overline{\mu}$ is different for the different choices. In the figure, on $\mathbb{S}^2$, the boundary of both $\overline{\xi}_1$ and $\overline{\xi}_2$ are $\overline{\varsigma}$. However the intersection number between $\overline{\xi}_1$ and $\overline{\mu}$ is zero, while that between $\overline{\xi}_2$ and $\overline{\mu}$ is $\pm 1$.]{
      \label{fig:linking-number-uniqueness-c}
      \includegraphics[width=0.3\textwidth, trim=100 0 100 0, clip=true]{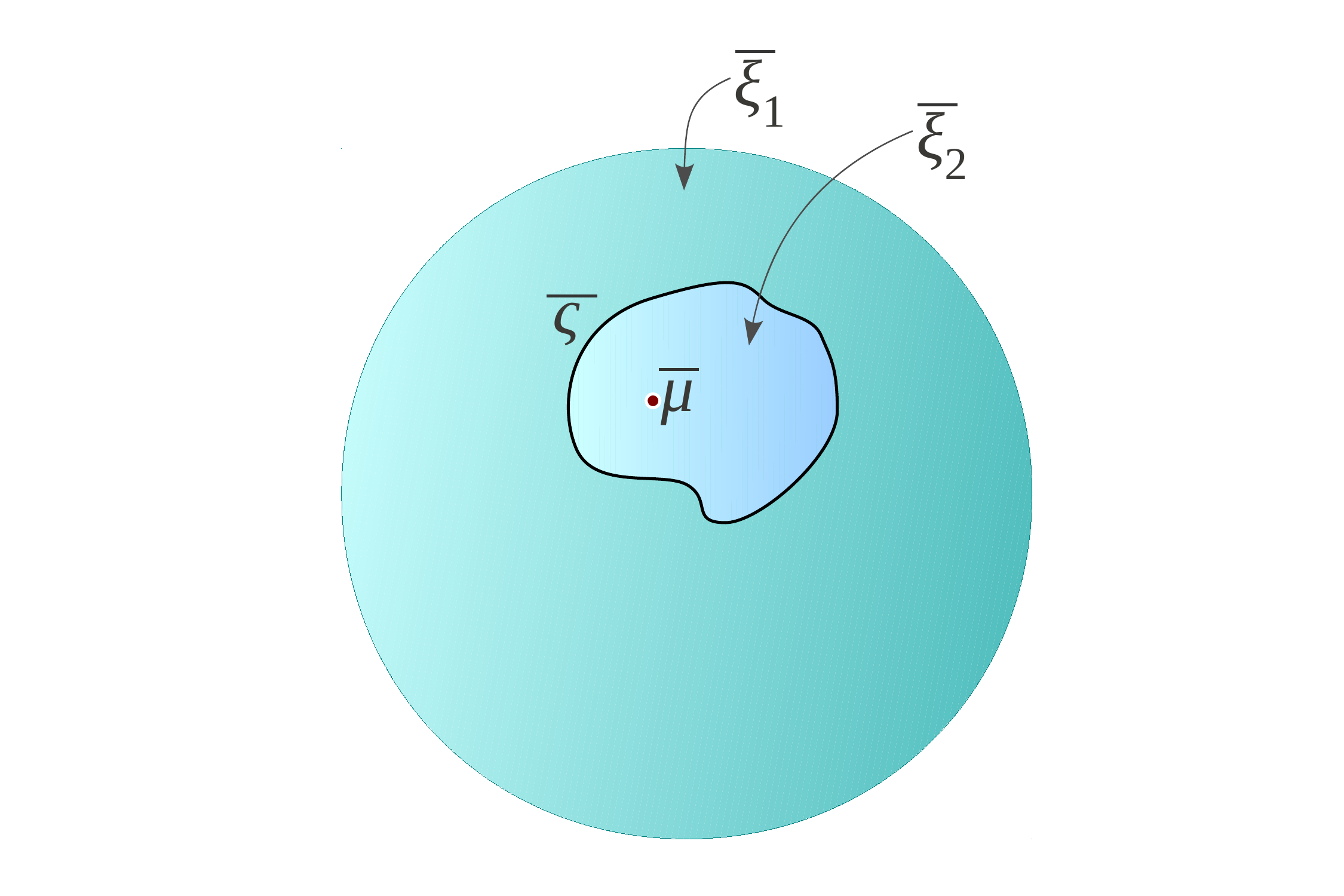}
      }
  \caption{Examples and counter-examples of uniqueness of linking number -- a consequence of Proposition~\ref{prop:linking-number-independence}.}
  \label{fig:linking-number-uniqueness}
\end{figure*}

The statement of the above proposition is about the uniqueness in the value of linking number.
Intersection number, according to definition, is between a $N$-chain (like $\overline{\xi}$ in the example of figure~\ref{fig:intersection-linking-number-simple} or \ref{fig:linking-number-uniqueness}) and a $(D-N)$-cycle (like $\overline{\mu}$ in the example of Figure~\ref{fig:intersection-linking-number-simple}). That is then, by definition, the linking number between the boundary of the $N$-chain, which is a $(N-1)$-cycle (like $\overline{\varsigma}$ in the example of Figure~\ref{fig:intersection-linking-number-simple}) and the $(D-N)$-cycle ($\overline{\mu}$).
However, it may be possible that there exists another $N$-chain, $\overline{\xi}'$, such that $\overline{\varsigma}$ is the boundary of that $N$-chain as well (\emph{i.e.} $\overline{\varsigma}$ is a common boundary between $\overline{\xi}$ and $\overline{\xi}'$ -- as illustrated in Figure~\ref{fig:linking-number-uniqueness-b}). Then, if the intersection number between $\overline{\xi}'$ and $\overline{\mu}$ is not same as that between $\overline{\xi}$ and $\overline{\mu}$, the definition of the linking number between $\overline{\varsigma}$ and $\overline{\mu}$ becomes ambiguous.
This is exemplified in Figure~\ref{fig:linking-number-uniqueness-c}.
Proposition~\ref{prop:linking-number-independence} precisely gives the condition under which such ambiguity is not present.
}

\vspace{0.1in}
\begin{proposition} [Connection to homology of $A$] \label{prop:linking-number-homology-injectivity}
 Consider a fixed non-zero $\mu \in H_{D-N}(Y,B)$.
 If, in addition to the condition of Proposition~\ref{prop:linking-number-independence}, we have $H_N(X,A)\approxeq H_{N-1}(A) \approxeq \mathbb{R}$, and if there exists at least one $(N-1)$-cycle in $A$ such that its linking number with $\mu$ is non-zero,
 then the value of $\mathscr{L}(\varsigma,\mu)$ tells us which element of $H_{N-1}(A)$ is the chosen $\varsigma$. In other words, the map $\mathcal{H} \equiv \mathscr{L}(\cdot,\mu): H_{N-1}(A) \rightarrow H_{D-1}(\mathbb{R}^D,\mathbb{R}^D-\{0\}) \approxeq \mathbb{R}$ is an injective homomorphism.
\end{proposition}
\begin{quoteproof}
 The map $\mathcal{H}$ is given by $\mathcal{H}(\varsigma) = (-1)^{D-N}p_{*}(\partial_{*}^{-1} \varsigma \times \mu)$. This clearly is a group homomorphism between $H_{N-1}(A)$ and $H_{D-1}(\mathbb{R}^D,\mathbb{R}^D-\{0\})$. Since by hypothesis, both the domain and the co-domain of $\mathcal{H}$ are isomorphic to $\mathbb{R}$, $\mathcal{H}$ can either be a trivial homomorphism (\emph{i.e.} maps everything in its domain to $0$ in its co-domain), or it can be an injection. The former possibility is ruled out by the hypothesis of existence of at least one $(N-1)$-cycle in $A$ with non-zero linking number with $\mu$. Thus the result follows.
\end{quoteproof}

The result implies that the linking number with $\mu$ is a \emph{complete invariant} for the homology class $\varsigma$.

\lo{So far we have been talking about intersection number and linking number. However what we are really interested in is the homology class of $\overline{\varsigma}$ in $A$ (in light of robot planning problems, that is the homology class of the closed trajectories in $(\mathbb{R}^D-S)$, as illustrated in figures~\ref{fig:intersection-linking-number-simple} ans \ref{fig:linking-number-uniqueness}).
The result of the above proposition establishes a relationship between the linking number between $\overline{\mu}$ and $\overline{\varsigma}$ (see figure~\ref{fig:intersection-linking-number-simple} or \ref{fig:linking-number-uniqueness}), and the homology class of $\overline{\varsigma}$.
It says that under certain conditions, the linking number will precisely tell us about the homology class of $\overline{\varsigma}$ (\emph{i.e.} a \emph{complete invariant}).}

% --------------------------------------------------------------------

\subsection{Computation of Intersection/Linking Number for Given Cycles}

We describe how to compute the linking number between the cycles $\overline{\varsigma}$ and $\overline{\mu}$.  As discussed in the beginning of this paper, we would like to be able to compute the homology class of $(N-1)$-cycles (top-dimensional cycles on $(N-1)$-dimensional manifolds) as an explicit number (or a set of numbers). Equipped with Proposition~\ref{prop:linking-number-homology-injectivity}, that problem can be converted to the problem of computation of the linking numbers.

\lo{
A $n$-form can always be integrated on an oriented $n$-cycle. However, the value of the integration may not tell us anything about the homology class of the cycle. For example, in $(\mathbb{R}^2-0)$, \emph{i.e.} the plane with the origin removed, $\d x$ is a differential $1$-form. However it is \emph{exact} and evaluates to $0$ on every closed curve. On the other hand $\d \theta = \frac{x \d y + y \d x}{x^2 + y^2} ~(=I\!m(\frac{\d z}{z})))$, which is \emph{closed} but not \emph{exact} \cite{bott1982differential} in $(\mathbb{R}^2-0)$, in fact tell us about the homology class of closed loops.

The purpose of the proposition below is to design a differential from, integration of which, along with the conditions of Proposition~\ref{prop:linking-number-homology-injectivity}, captures the homology class of $(N-1)$-cycles in $A$ (which, in light of robot planning problem, are the punctured spaces $(\mathbb{R}^D-S)$, $S$ being path-connected).
In order to achieve this for arbitrary $A$ (in robot planning, for example, we can arbitrary representatives, $S$, for the obstacles), we exploit the transformation, $p$, in the definition of linking numbers. Thus, the closed, non-exact differential form that we have to choose is one from the co-domain of $p$, namely $(\mathbb{R}^D-0)$ (a space which is much simpler and well-known than, say, $(\mathbb{R}^D-S)$), and then \emph{pull it back} to the original space by $p$. Thus we have the following proposition.
 %Consider the example in Figure~\ref{fig:diff-forms-illustration}
\vspace{0.1in}}

Let $\eta_0 \in ~\Omega^{D-1}_{dR}(\mathbb{R}-\{0\})$ be a closed differential form % in $(\mathbb{R}-\{0\})$ such that $[\eta_0] \in H^{D-1}_{dR}(\mathbb{R}-\{0\}) \approxeq \mathbb{R}$ is a generator of $H^{D-1}_{dR}(\mathbb{R}-\{0\})$, agreeing with the usual orientation.
that represents the standard generator of $H^{D-1}(\mathbb{R}^D - \{0\})$.  Let $j_* \colon H_{D_N}(Y) \to H_{D-N}(Y,B)$ denote the quotient map.

\begin{proposition} \label{prop:linking-number-integration}
 Assume the same hypotheses as in Proposition~\ref{prop:linking-number-independence}.  Fix $\mu \in H_{D-N}(Y,B)$, and suppose there
 exists a class $u \in H_{D-N}(Y)$ such that $j_*(u) = \mu$.  Then for any $\varsigma \in H_{N-1}(A)$, the linking number $\mathscr{L}(\varsigma, \mu)$ is uniquely determined by the value of the integral
\begin{equation} \label{eq:linking-number-integration}
(-1)^{D-N} \int_{\varsigma \times u} p^{*}(\eta_0).
\end{equation}
%where $u \in H_{D-N}(Y)$ is any class such that ${j_*}(u) = \mu$, where $j_*$ is the map $C_{D-N}(Y) \rightarrow C_{D-N}(Y,B)$. (See Thm. 2.13 of \emph{\cite{Hatcher:AlgTop}} -- note that for a given $\overline{\mu}$, in general, there can be many possible choices for $\overline{u}$).
\end{proposition}
\begin{quoteproof}
First, note that the map
\[
H_D(\mathbb{R}^D, \mathbb{R}^D - \{0\}) \xrightarrow{\partial_*} H_{D-1}(\mathbb{R}^D - \{0\})
\xrightarrow{\int_{\cdot}\eta_0} \mathbb{R}
\]
is an isomorphism, so that every element $m \in H_D(\mathbb{R}^D, \mathbb{R}^D -
\{0\})$ is uniquely determined by the value of the integral
$\int_{\partial_* m} \eta_0$.

Choose a class $\zeta \in H_{N}(X,A)$ such that $\partial_*(\zeta) =
\varsigma$.  Then, by definition,
\[
\mathscr{L}(\varsigma, \mu) = \mathscr{I}(\zeta,\mu) = (-1)^{D-N} p_*(\zeta \times \mu) \in H_D(\mathbb{R}^D, \mathbb{R}^D - \{0\}).
\]
Now, consider the diagram below.\vspace{-0.05in}
\[
\xymatrix@C=0pt@R=12pt{
  &
  H_N(X,A) \otimes H_{D-N}(Y) 
  \ar[dl]^(.4){1 \otimes j_*} \ar[dr]_(.4){\partial_* \otimes 1} \ar[dd]^{\times}
  &
  \\
  H_N(X,A) \otimes H_{D-N}(Y,B)
  \ar[dd]^{\times}
  &
  &
  H_{N-1}(A) \otimes H_{D-N}(Y) 
  \ar[dd]^{\times}
  \\
  &
  H_D(X \times Y, A \times Y) 
  \ar[dl]^(.4){j_*} \ar[dr]_(.4){\partial_*} \ar[dd]^{p_*}
  &
  \\
  H_D(X\times Y, A\times Y \cup X \times B) 
  \ar[dr]_{p_*}
  &
  &
  H_{D-1}(A \times Y) 
  \ar[dd]^{p_*}
  \\
  &
  H_D(\mathbb{R}^D, \mathbb{R}^D - \{0\})
  \ar[dr]_(.4){\partial_*}
  &
  \\
  &
  &
  H_{D-1}(\mathbb{R}^D - \{0\})
}
\]
It is a standard fact that every part of this diagram commutes, and as
a consequence we have that
\[
\partial_* p_*(\zeta \times \mu)
= \partial_* p_*(\zeta \times j_*u) = p_*(\partial_* \zeta \times u) = p_*(\varsigma \times u)
\]
%By the commutativity of the diagram above, 
%$\partial_* \mathscr{L}(\varsigma, \mu) = p_*(\varsigma \times u)$.
Finally, by the naturality of integration, we have
\[
\int_{\partial_* \mathscr{L}(\varsigma, \mu)} \eta_0 = 
(-1)^{D-N} \int_{p_*(\varsigma \times u)} \eta_0 =
(-1)^{D-N} \int_{\varsigma \times u} p^*(\eta_0).
\]
Thus the integral on the right uniquely determines the
value of the linking number $\mathscr{L}(\varsigma, \mu)$.
\end{quoteproof}

\changedA{Note that linking number, by definition, is defined between a cycle in $A$ and a relative cycle in $(Y,B)$. However, for computing the integration of Equation~\eqref{eq:linking-number-integration}, the cycles we choose are from $A$ and $Y$. Thus it is possible to use the standard notion of integration over chains~\cite{bott1982differential}. However, if $B=\emptyset$, a relative cycle in $(Y,B)$ becomes a cycle in $Y$.}

% =======================================================

\section{Construction and Explicit Computation} \label{sec:homology-invariant-formula}
\subsection{Construction of the Complete Invariant}

We specialize the results of the previous section to match the description of the \emph{reduced problem definition} in Section~\ref{sec:reduced-problem-algtop}. At present, we consider the case where there is a single path-connected component of $\widetilde{\mathcal{S}}$, namely $S$. In connection to the definitions stated in Section~\ref{sec:linking-num-pre} (cf. Figure~\ref{fig:intersection-linking-number-simple}), we set
\[
    X = \mathbb{R}^D, ~A = \mathbb{R}^D - S, ~~Y = S ~\textrm{ and } B=\emptyset
\]
Moreover, since $Y\equiv S$ is a $(D-N)$-dimensional closed, connected and oriented manifold, we have $H_{D-N}(S) \approxeq \mathbb{R}$.  We thus choose $\overline{\mu}=\overline{S} \in Z_{D-N}(S)$ to be a cycle representing the fundamental class of $S$, i.e. the generator $\mathbf{1} \in H_{D-N}(S)$. % such that $[\overline{S}] = \mathbf{1} \in H_{D-N}(S)$ is a generator of $H_{D-N}(S)$ (which is isomorphic to $\mathbb{R}$ since $S$ is a closed, path-connected and orientable $(D-N)$-dimensional manifold).
Also, note that since $B=\emptyset$, the map $j': Z_{D-N}(Y) \rightarrow Z_{D-N}(Y,B)$ is the identity map. So in this case $[\overline{S}] \in H_{D-N}(S,B) \equiv H_{D-N}(S)$.

\begin{figure}
\centering
 \includegraphics[width=0.35\textwidth, trim=80 50 80 30, clip=true]{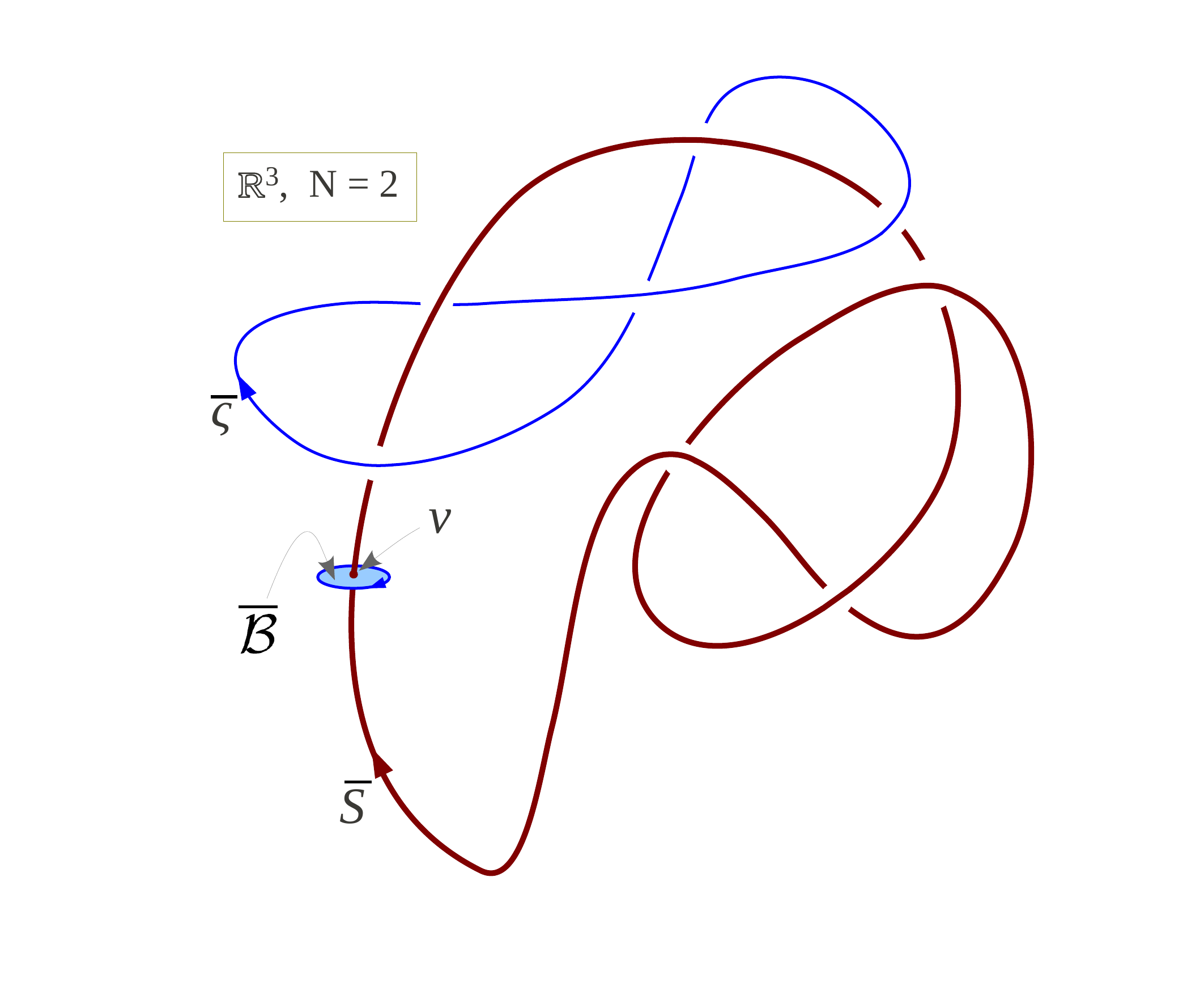}
 \caption{The specific problem under consideration, illustrated for $D=3, N=2$.}
\label{fig:problem}
\end{figure}

For this choice it is easy to verify that the conditions of Propositions \ref{prop:linking-number-independence}, \ref{prop:linking-number-homology-injectivity} and \ref{prop:linking-number-integration} hold.
\begin{itemize}
 \item[i.] \emph{Proposition~\ref{prop:linking-number-independence}:} $H_N(\mathbb{R}^D) = H_{N-1}(\mathbb{R}^D) = 0$ follows from contractibility of $\mathbb{R}^D$.
 \item[ii.] \emph{Proposition~\ref{prop:linking-number-homology-injectivity}:}
 \begin{itemize}
   \item[a.] By Alexander duality \cite{Hatcher:AlgTop}, $ H_N (\mathbb{R}^D, \mathbb{R}^D - S) \approxeq H^{D-N} (S)$. Using Poincar\'e Duality for $S$, %(which is a $(D-N)$-dimensional closed, orientable manifold),
$H^{D-N} (S) \approxeq H_0(S) \approxeq \mathbb{R}$. % (since $S$ has a single connected component).
   Finally, from the long exact sequence for the pair $(\mathbb{R}^D, \mathbb{R}^D - S)$, using the contractibility of $\mathbb{R}^D$, we have, $H_N (\mathbb{R}^D, \mathbb{R}^D - S) \approxeq H_{N-1}(\mathbb{R}^D - S)$. Combining these three isomorphisms we have, \begin{equation}\label{eq:single-puncture-order-1} H_N (\mathbb{R}^D, \mathbb{R}^D - S) \approxeq H_{N-1}(\mathbb{R}^D - S) \approxeq \mathbb{R} \end{equation}.
   \item[b.] Consider a point $v\in S$. Since $\overline{S}$ covers $S$, this point is also in (the image of) $\overline{S}$. Since $S$ is $(D-N)$-dimensional, we can choose a small $N$-ball, $\mathcal{B}$, centered at $v$ such that it intersects $S$ transversely only at $v$. Let $\overline{\mathcal{B}} \in C_{N}(\mathbb{R}^D)$ be a top-dimensional non-zero chain that covers $\mathcal{B}$. Clearly the intersection number between $\overline{S}$ and $j(\overline{\mathcal{B}})$ (where $j: \mathbb{R}^D \rightarrow \mathbb{R}^D/(\mathbb{R}^D-S)$ is the quotient map) is non-zero. Thus the linking number between $\partial \overline{\mathcal{B}} \big|_{(\mathbb{R}^D - S)}$ (which, by our construction, is a $(N-1)$-cycle in $(\mathbb{R}^D - S)$) and $\overline{S}$ is non-zero. Thus there exists at least one $(N-1)$-cycle in $(\mathbb{R}^D-S)$ that has non-zero linking number with $\overline{S}$ (see Figure~\ref{fig:problem}).
 \end{itemize}
 \item[iii.]  \emph{Proposition~\ref{prop:linking-number-integration}:} Follows from the fact that $B = \emptyset$.
\end{itemize}

 {\bf Construction:} A complete invariant for homology classes of $(N-1)$-cycles, $\overline{\omega} \in Z_{N-1}(\mathbb{R}^D-S)$, is, by Proposition~\ref{prop:linking-number-homology-injectivity}, the linking number between $\overline{\omega}$ and $\overline{S}$. Using Proposition~\ref{prop:linking-number-integration}, the complete invariant, $\phi_S$, for the homology classes of such chains is given by the integral
\begin{eqnarray} \label{eq:pullback-integral}
 \mathcal{\phi}_{S}(\overline{\omega}) & = & (-1)^{D-N} \int_{\overline{\omega} \times \overline{S}} p^{*}(\eta_0) \nonumber \\
 & = & (-1)^{D-N} \int_{\overline{\omega}} \int_{\overline{S}} p^{*}(\eta_0) ~~~\textrm{ [Fubini theorem]}
\end{eqnarray}

% ---------------------------------------------

\subsection{Computation of $\phi_S$}% using a particular choice for $\eta_0$}

Let $\mathbf{x} \in (\mathbb{R}^D \setminus S) \subset \mathbb{R}^D$ be the coordinate variable describing points in $(\mathbb{R}^D \setminus S)$, and let $\mathbf{x'} \in S \subset \mathbb{R}^D$ be the one describing points in $S$. Thus we have $p(\mathbf{x},\mathbf{x'}) = \mathbf{x} - \mathbf{x'}$.
%Let $\mathbf{s} \in (\mathbb{R}^D \setminus \{0\}) \subset \mathbb{R}^D$ be the natural coordinate variable describing points in the space $(\mathbb{R}^D \setminus \{0\})$.
 A well-known \cite{William:Clifford, Differential:flanders:1989} explicit generator for the deRham cohomology $H_{dR}^{D-1}(\mathbb{R}^D\setminus \{0\})$) is,
\begin{equation}
 \eta_0%(\mathbf{s})
 =
 \sum_{k=1}^D \mathcal{G}_k %(\mathbf{s}) ~
 ~(-1)^{k+1} ~~\d s_1 \wedge \cdots \wedge \d s_{k-1} \wedge \d s_{k+1} \wedge \cdots \wedge \d s_D
 =
 \sum_{k=1}^D \mathcal{G}_k %(\mathbf{s}) ~
 ~(-1)^{k+1} ~\bigwedge_{i=1 \atop i\neq k}^{D} \d s_i
 \label{eq:omega-in-s}
\end{equation}
where,
\begin{equation}
 \mathcal{G}_k(\mathbf{s}) = \frac{1}{A_{D-1}}~ \frac{s_k}{\left( s_1^2 + s_2^2 + \cdots + s_D^2 \right)^{D/2}} \label{eq:G-k}
\end{equation}
for $\mathbf{s} = (s_i) \in (\mathbb{R}^D\setminus \{0\})$, and $A_{D-1} = \frac{ D \pi^{\frac{D}{2}}}{\Gamma (\frac{D}{2} + 1)}$, the $(D-1)$-volume of the $(D-1)$-dimensional unit sphere.

The pullback of $\eta_0$ under $p$ is given by the following formula,
\begin{eqnarray}
 \eta(\mathbf{x},\mathbf{x'}) & = & p^{*}(\eta_0) ~=~ \eta_0 \big|_{\mathbf{s} = \mathbf{x} - \mathbf{x'}} \changedC{~=~ \sum_{k=1}^D \mathcal{G}_k ~(-1)^{k+1} } 
 \bigwedge_{i=1 \atop i\neq k}^D d(x_i-x_i')
 %\nonumber
 %\\
 %& = & \sum_{k=1}^D \mathcal{G}_k (\mathbf{x} - \mathbf{x'}) ~~(-1)^{k+1} ~~\d (x_1-x'_1) \wedge \d (x_2-x'_2) \wedge \cdots \nonumber \\
 %    & &  \qquad\qquad\qquad \wedge \d (x_{k-1}-x'_{k-1}) \wedge \d (x_{k+1}-x'_{k+1}) \wedge \cdots \wedge \d (x_D-x'_D)
\label{eq:eta}\end{eqnarray}

Now consider the quantity of interest, $\phi(\overline{\omega}) = \int_{\mathbf{x} \in \overline{\omega}} \int_{\mathbf{x'} \in \overline{S}}  \eta(\mathbf{x}, \mathbf{x'})$. On $\overline{\omega}\times\overline{S}$, at most $(N-1)$ unprimed differentials can be independent, and at most $(D-N)$ primed differentials can be independent (since $\mathbf{x}$ represents a point on the image of the $(N-1)$ chain $\overline{\omega}$ and $\mathbf{x'}$ represents a point on the image of the $(D-N)$ chain $\overline{S}$). Thus we can conveniently drop all the terms in the expansion of $\eta$ (which is a $(D-1)$-differential form on $(\mathbb{R}^D-S)\times S$) that do not satisfy these conditions on maximum number of primed/unprimed differentials. Thus we obtain a simpler differential form $\tilde{\eta}$,

\begin{equation}
 \tilde{\eta}(\mathbf{x},\mathbf{x'})  =  \sum_{k=1}^D
    \left(
            \mathcal{G}_k (\mathbf{x}-\mathbf{x'}) ~~(-1)^{k+1+D-N}
            \sum_{\tau_i \in \{0,1\} \atop \tau_1+\cdots+\tau_D=D-N}
            \bigwedge_{i=1 \atop i\neq k}^D dx_i^{(\tau_i)}
     \right)
 %\cdot \nonumber \\
 %       & &  \qquad~~~ \!\!\!\!\!\!\!\!\! \sum_{\mysubstack{\tau_i \in \{0,1\}}{\tau_1+\cdots+\tau_D=D-N}} \!\!\!\!\!\!\!\! \d x^{(\tau_1)}_1 \wedge \d x^{(\tau_2)}_2 \wedge \cdots \wedge \d x^{(\tau_{k-1})}_{k-1} \wedge \d x^{(\tau_{k+1})}_{k+1} \wedge \cdots \wedge \d x^{(\tau_D)}_D  \Bigg)  \nonumber \\
\label{eq:eta-tilde}
\end{equation}
[where, $x_i^{(\tau)}$ represents $x'_i$ if $\tau=1$, otherwise represents $x_i$ if $\tau=0$.]

\noindent This differential form, though simpler, has the property that
\begin{equation}
 \phi_S(\overline{\omega}) \quad=\quad (-1)^{D-N} \int_{\mathbf{x} \in \overline{\omega}} \int_{\mathbf{x'} \in \overline{S}} \eta(\mathbf{x}, \mathbf{x'}) \quad=\quad (-1)^{D-N} \int_{\mathbf{x} \in \overline{\omega}} \int_{\mathbf{x'} \in \overline{S}} \tilde{\eta}(\mathbf{x}, \mathbf{x'})
\end{equation}
Finally, we re-write the formula for $\tilde{\eta}$ using a new notation as follows,
\begin{eqnarray}
 \tilde{\eta}(\mathbf{x},\mathbf{x'}) & = & (-1)^{D-N} \sum_{k=1}^D \Bigg( \mathcal{G}_k (\mathbf{x}-\mathbf{x'}) ~~(-1)^{k+1} \cdot \nonumber \\
    & & \qquad \sum_{~~\rho \in {part}^{D-N} (\mathcal{N}^D_{-k})} \!\!\!\!\!\! \textrm{sgn}(\rho) ~~\d x'_{\rho_l(1)} \wedge \cdots \wedge \d x'_{\rho_l(D-N)} \wedge \d x_{\rho_r(1)} \wedge \cdots \wedge \d x_{\rho_r(N-1)} \Bigg) \nonumber \\ \label{eq:def-eta-bar-partition}
\label{eq:eta-tilde-part} \end{eqnarray}
where,
\begin{enumerate}
 \item $\mathcal{N}^D_{-k} = [1,2,\cdots,k-1,k+1,\cdots,D ]$ is an ordered set,
 \item ${part}^w(\mathcal{A})$ is the set of all $2$ partitions of the ordered set $\mathcal{A}$, such that the first partition contains $w$ elements, and each of the partitions contain elements in order. The sign of an element from the set is the permutation sign of the ordered set formed by concatenating the two partitions of the element. \lo{\footnote{Let us consider an ordered set $\mathcal{A} = [a_1,a_2,\cdots,a_q]$ with $a_1\leq a_2\leq \cdots \leq a_q$ (where the inequality sign signifies order of arrangement and not necessarily the order of magnitude). We represent the set of all ordered $2$-partitions of the set $\mathcal{A}$ into $w$ and $q-w$ elements as $~{part}^w(\mathcal{A})$, such that for a $\rho = [\rho_l,\rho_r] \in {part}^w(\mathcal{A})$, $\rho_l$ and $\rho_r$ are ordered sets of $w$ and $q-w$ elements respectively, with the properties that $\rho_l \cap \rho_r = \emptyset$, $\rho_l(1)\leq\rho_l(2)\leq\cdots\leq\rho_l(w)$ and $\rho_r(1)\leq\rho_r(2)\leq\cdots\leq\rho_r(q-w)$. Then the sign of the partition, $~\textrm{sgn}(\rho)$, is defined as the permutation sign of the ordered set $\rho_l \sqcup \rho_r$. For example,

 \begin{minipage}[r]{0.9\textwidth} {\small ${part}^3([1,3,6,9,5]) = \big\{ ~
 \left[ [ 1,3,6], [ 9,5] \right],$
  $\left[ [ 1,3,9], [ 6,5] \right],$ $
 \left[ [ 1,3,5], [ 6,9] \right],
 $ $\left[ [ 1,6,9], [ 3,5] \right],$
  $ \left[ [ 1,6,5], [ 3,9] \right],$ $
 \left[ [ 1,9,5], [ 3,6] \right],
 $ $\left[ [ 3,6,9], [ 1,5] \right],$
  $ \left[ [ 3,6,5], [ 1,9] \right],$ $
 \left[ [ 3,9,5], [ 1,6] \right],
 $ $\left[ [ 6,9,5], [ 1,3] \right]~
\big\}$.}\end{minipage}

\noindent Then if $\rho = \left[ [ 1,6,5], [ 3,9] \right] \in {part}^3([1,3,6,9,5])$, we write $\rho_l = [ 1,6,5]$ and $\rho_r = [ 3,9]$. Also, the $j^{th}$ element of $\rho_b,~b\in\{l,r\}$ is written as $\rho_b(j)$. Thus, in the example, $\rho_l(2) = 6$.}}
\end{enumerate}

\noindent Thus, \sh{after some simplification, }the final formula for the complete invariant for homology class of $\overline{\omega}\in Z_{N-1}(\mathbb{R}^D-S)$ is, \vspace{-0.05in}
% \ls{\begin{eqnarray}
%  \phi_S(\overline{\omega}) & = & (-1)^{D-N} \int_{\mathbf{x} \in \overline{\omega}} \int_{\mathbf{x'} \in \overline{S}} \tilde{\eta}(\mathbf{x}, \mathbf{x'}) \nonumber \\
%  & = & \int_{\mathbf{x} \in \overline{\omega}}~~ \sum_{k=1}^D \!\!\!\!\!\!\sum_{~~~~~~\rho \in {part}^{D-N} (\mathcal{N}^D_{-k})} \nonumber \\
%        & & \qquad \left( (-1)^{k+1}  \int_{\mathbf{x'} \in \overline{S}} \mathcal{G}_k (\mathbf{x}-\mathbf{x'}) ~~\textrm{sgn}(\rho) ~~\d x'_{\rho_l(1)} \wedge \cdots \wedge \d x'_{\rho_l(D-N)} \right) \nonumber \\
%        & & \qquad\qquad\qquad\qquad\qquad\qquad\qquad\qquad\qquad\qquad \wedge \d x_{\rho_r(1)} \wedge \cdots \wedge \d x_{\rho_r(N-1)}  \nonumber \\
%  & = & \int_{\mathbf{x} \in \overline{\omega}}~~ \sum_{k=1}^D \!\!\!\!\!\!\sum_{~~~~~~\rho \in {part}^{D-N} (\mathcal{N}^D_{-k})} \!\!\!\!\!\! U^k_{\rho} (\mathbf{x};S) \wedge \d x_{\rho_r(1)} \wedge \cdots \wedge \d x_{\rho_r(N-1)}  \label{eq:int-eta-final}
% \end{eqnarray}}
{\begin{eqnarray}
 \phi_S(\overline{\omega}) & = & (-1)^{D-N} \int_{\mathbf{x} \in \overline{\omega}} \int_{\mathbf{x'} \in \overline{S}} \tilde{\eta}(\mathbf{x}, \mathbf{x'}) \nonumber \\
% & = & \int_{\mathbf{x} \in \overline{\omega}}~~ \sum_{k=1}^D \!\!\!\!\!\!\sum_{~~~~~~\rho \in {part}^{D-N} (\mathcal{N}^D_{-k})} \nonumber \\
%       & & \qquad \left( (-1)^{k+1}  \int_{\mathbf{x'} \in \overline{S}} \mathcal{G}_k (\mathbf{x}-\mathbf{x'}) ~~\textrm{sgn}(\rho) ~~\d x'_{\rho_l(1)} \wedge \cdots \wedge \d x'_{\rho_l(D-N)} \right) \nonumber \\
%       & & \qquad\qquad\qquad\qquad\qquad\qquad\qquad\qquad\qquad\qquad \wedge \d x_{\rho_r(1)} \wedge \cdots \wedge \d x_{\rho_r(N-1)}  \nonumber \\
 & = & \int_{\mathbf{x} \in \overline{\omega}}~~ \sum_{k=1}^D \!\!\!\!\!\!\sum_{~~~~~~\rho \in {part}^{D-N} (\mathcal{N}^D_{-k})} \!\!\!\!\!\! U^k_{\rho} (\mathbf{x};S) \wedge \d x_{\rho_r(1)} \wedge \cdots \wedge \d x_{\rho_r(N-1)}  \label{eq:int-eta-final}
\end{eqnarray}}
where, 
\begin{equation}
  U^k_{\rho} (\mathbf{x};S) = (-1)^{k+1} ~~\textrm{sgn}(\rho) \int_{\mathbf{x'} \in \overline{S}} \mathcal{G}_k (\mathbf{x}-\mathbf{x'}) ~~\d x'_{\rho_l(1)} \wedge \cdots \wedge \d x'_{\rho_l(D-N)} \label{eq:U-final}
\end{equation}
and by convention, $\overline{S}$ is a top-dimensional cycle covering $S$ such that $[\overline{S}] = \mathbf{1} \in H_{D-N}(S)$.

Also, note that the quantity inside the integral in the formula for $\phi_S$ is a differential $(N-1)$-form in $(\mathbb{R}^D-S)$. Thus we can integrate it over $\overline{\omega}$. We represent the differential $(N-1)$-form by $\psi_S$
 \begin{equation} \label{eq:omg-final} \vspace{-0.05in}
 \psi_S = \sum_{~~~~~~\rho \in {part}^{D-N} (\mathcal{N}^D_{-k})} \!\!\!\!\!\! U^k_{\rho} (\mathbf{x};S) \wedge \d x_{\rho_r(1)} \wedge \cdots \wedge \d x_{\rho_r(N-1)}
\end{equation}

\vspace{0.05in}
\changedA{It should be noted that the $\eta_0$ we used in \eqref{eq:omega-in-s} is just a particular choice, but this choice is the only symmetric one (up to a scalar multiple) under rotations about the origin.  This symmetry enables us to write a clean formula in coordinates, but in general any closed and non-exact form $\eta_0$ would work. The resulting invariant would differ from ours by a constant multiple.  %One could very well choose some other $\eta_0 \in Z_{dR}^{D-1}(\mathbb{R}^D - \widetilde{\mathcal{S}})$ (for example, the bump $D$-form \cite{bott1982differential}), and obtain a different formula for $\psi_S$ which would be an equally valid complete invariant.
}

% --------------------------------------

\subsection{Incorporating Multiple Connected Components of $\widetilde{\mathcal{S}}$}

So far we have worked in the case of a single connected obstacle $S$. %with a single connected component of the puncture, namely, $S$.
However, recall that the original space under consideration was $(\mathbb{R}^D - \widetilde{\mathcal{S}})$, with \changedA{$\widetilde{\mathcal{S}} = \bigsqcup_{i=1}^m S_i$}, such that each $S_i$ is a path connected, closed, locally contractible and orientable \changedA{$(D-N)$-manifold}.
%Moreover, by hypothesis, \changedA{$S_i \cap S_j = \emptyset, ~\forall i\neq j$}.
A straightforward induction argument computes the homology of the smaller space, $(\mathbb{R}^D - \widetilde{\mathcal{S}})$, in terms of the larger spaces, $(\mathbb{R}^D \setminus S_k)$.

\begin{proposition} \label{prop:direct-sum-homology}
 $H_{N-1}(\mathbb{R}^D \setminus \widetilde{\mathcal{S}}) ~\approxeq~ \bigoplus_{k=1}^m H_{N-1}(\mathbb{R}^D \setminus S_k) ~\approxeq~ \mathbb{R}^m$, where the first isomorphism is induced by the direct sum of the inclusion maps $\tilde{i}_k: (\mathbb{R}^D - \widetilde{\mathcal{S}}) \hookrightarrow (\mathbb{R}^D - S_k)$.
\end{proposition}
\begin{quoteproof}
Recall that the spaces $S_i$ are pairwise disjoint, so that for any $p$
\begin{eqnarray} \label{eq:punctured-relative-homology}
 (\mathbb{R}^D-S_p) ~\cup~ (\mathbb{R}^D\setminus \changedA{\sqcup_{i=p+1}^m S_i}) &=& \mathbb{R}^D %\changedA{\setminus ~(S_p ~\cap~ (\sqcup_{i=p+1}^m S_i))}
\nonumber \\
 (\mathbb{R}^D-S_p) ~\cap~ (\mathbb{R}^D\setminus \changedA{\sqcup_{i=p+1}^m S_i}) &=& \mathbb{R}^D\setminus \changedA{\sqcup_{i=p}^m S_i} \nonumber
\end{eqnarray}
%\changedA{Now, each of $S_i$ being a $(D-N)$-dimensional compact, locally contractible and orientable manifold embedded in $\mathbb{R}^D$, by transversality theorem (and upon performing small perturbations of $S_i$ if required), $(S_p ~\cap~ (\sqcup_{i=p+1}^m S_i))$ is at most a $\big( 2(D-N)-D = \big) D-2N$ dimensional compact, locally contractible and orientable spaces.
%Therefore by Proposition 3.46 of \cite{Hatcher:AlgTop},
%\begin{eqnarray}
% H_n \big(\mathbb{R}^D, ~\mathbb{R}^D \!- (S_p ~\cap~ (\sqcup_{i=p+1}^m S_i)) \big) & \approxeq & H^{D-n} \big(S_p ~\cap~ (%\sqcup_{i=p+1}^m S_i) \big) \nonumber \\
% & \approxeq & 0, ~\textrm{for all $n \leq 2N-1$}
%\end{eqnarray}
%By hypothesis, $N < N+1 \leq 2N-1$ (since we assumed $N>1$).
%Thus from the long exact sequence for the pair $(\mathbb{R}^D, ~\mathbb{R}^D \!- (S_p ~\cap~ (\sqcup_{i=p+1}^m S_i)) \big)$, using \eqref{eq:punctured-relative-homology} and the contractibility of $\mathbb{R}^D$, we obtain
%\begin{equation} \label{eq:null-hom-punctured}
% H_{N} \big(\mathbb{R}^D \!- (S_p ~\cap~ (\sqcup_{i=p+1}^m S_i)) \big) ~\approxeq~ H_{N-1} \big(\mathbb{R}^D \!- (S_p ~\cap%~ (\sqcup_{i=p+1}^m S_i)) \big) ~\approxeq~ 0
%\end{equation}
%}
%\changedA{Substituting \eqref{eq:null-hom-punctured} into
From the Mayer-Vietoris sequence~\cite{Hatcher:AlgTop} for the triad $\big( \mathbb{R}^D; ~\mathbb{R}^D\!-S_p, ~\mathbb{R}^D\!\setminus \sqcup_{i=p+1}^m S_i \big)$, one obtains an isomorphism %the following exact sequence
%}
%for $p=1,2,\cdots, m-1$.
%For an arbitrary $p$, using the contractibility of $\mathbb{R}^D$, the part of the sequence that is of interest to us is,
\begin{equation} \label{eq:mayer-puncture-seq}
\qquad
%0 ~~\xrightarrow{~~~~}~~
    H_{N-1}(\mathbb{R}^D\setminus \sqcup_{i=p}^m S_i)
~~\xrightarrow{(\tilde{u}_{p*}, \tilde{v}_{p*})}~~
    H_{N-1}(\mathbb{R}^D-S_p) \oplus H_{N-1}(\mathbb{R}^D\setminus \sqcup_{i=p+1}^m S_i)
%~~\xrightarrow{~~~~}~~ 0
\qquad
\end{equation}
%Due to exactness, the middle map (which is induced by the respective inclusion maps) is an isomorphism.
Note that $\tilde{u}_{1*} = \tilde{i}_{1*}$ and, $\tilde{v}_{1*} \circ \tilde{v}_{2*} \circ \cdots \tilde{v}_{(p-1)*} \circ \tilde{u}_{p*} = \tilde{i}_{p*}$.

By induction on $p$, we obtain a sequence of isomorphisms
%Starting from $p=1$, if we successively apply the isomorphism to the second homology group in sequence~\eqref{eq:mayer-puncture-seq}, we obtain the following sequence of isomorphisms,
\begin{eqnarray}
 H_{N-1}(\mathbb{R}^D\setminus \widetilde{\mathcal{S}})  & \xrightarrow[\approxeq]{~~(\tilde{i}_{1*}, \tilde{v}_{1*})~~} & H_{N-1}(\mathbb{R}^D-S_1) \oplus H_{N-1}(\mathbb{R}^D\setminus \sqcup_{i=2}^m S_i) \nonumber \\
  & \xrightarrow[\approxeq]{(\tilde{i}_{1*}, \tilde{i}_{2*}, \tilde{v}_{2*})} & H_{N-1}(\mathbb{R}^D-S_1) \oplus H_{N-1}(\mathbb{R}^D-S_2) \oplus H_{N-1}(\mathbb{R}^D\setminus \sqcup_{i=3}^m S_i) \nonumber \\
 & \cdots  & \cdots  \nonumber \\
 &  \xrightarrow[\approxeq]{~~\oplus_{k=1}^{m} \tilde{i}_{k*}~~} & \bigoplus_{k=1}^{m} H_{N-1}(\mathbb{R}^D-S_k)
\end{eqnarray}
The fact that this is isomorphic to $\mathbb{R}^m$ follows from Equation~\eqref{eq:single-puncture-order-1}, where we showed $H_{N-1}(\mathbb{R}^D-S_k) \approxeq \mathbb{R}$.
\end{quoteproof}

\vspace{0.1in}
The following theorem hence follows directly from Propositions~\ref{prop:direct-sum-homology} and Equation~\eqref{eq:pullback-integral}.

\begin{theorem}\label{thm:complete-invariant}
 For any $\overline{\omega} \in Z_{N-1}(\mathbb{R}^D \setminus \widetilde{\mathcal{S}})$, a complete invariant for the homology class of $\overline{\omega}$ is given by, 
 \begin{equation} \phi_{\widetilde{\mathcal{S}}}(\overline{\omega}) ~~\stackrel{\scriptscriptstyle def.}{=}~~ \left[ \begin{array}{c} \phi_{S_1}(\overline{\omega}) \\ \phi_{S_2}(\overline{\omega}) \\ \vdots \\ \phi_{S_m}(\overline{\omega}) \end{array} \right] \label{eq:final-psi-eqn} \end{equation}
 where, $\phi_{S_i}$ is given by the formula in Equation~\eqref{eq:int-eta-final}.
\end{theorem}
% Thus, for any $\overline{\omega} \in Z_{N-1}(\mathbb{R}^D \setminus \widetilde{\mathcal{S}})$, a complete invariant for the homology class of $\overline{\omega}$ is given by,
%  \begin{equation} \phi_{\widetilde{\mathcal{S}}}(\overline{\omega}) ~~\stackrel{\scriptscriptstyle def.}{=}~~ \left[ \begin{array}{c} \phi_{S_1}(\overline{\omega}) \\ \phi_{S_2}(\overline{\omega}) \\ \vdots \\ \phi_{S_m}(\overline{\omega}) \end{array} \right] \label{eq:final-psi-eqn} \end{equation}
%  where, $\phi_{S_i}$ is given by the formula in Equation~\eqref{eq:int-eta-final}.

Note that we have implicitly assumed a inclusion map $\tilde{i}_k: (\mathbb{R}^D - \widetilde{\mathcal{S}}) \hookrightarrow (\mathbb{R}^D - S_k)$ being applied on $\overline{\omega}$ for computation of the $k^{th}$ component. For simplicity we do not write it explicitly, since the map is identity as far as computation is concerned.

Thus, $[\overline{\omega}_1] = [\overline{\omega}_2]$ if and only if ~$\phi_{\widetilde{\mathcal{S}}}(\overline{\omega}_1) = \phi_{\widetilde{\mathcal{S}}}(\overline{\omega}_2)$, for any $\overline{\omega}_1, \overline{\omega}_2 \in Z_{N-1}(\mathbb{R}^D - \widetilde{\mathcal{S}})$.

% Following the same method as in Proposition \ref{prop:isomorphism-homology-int}, we obtain $H_{N-1}(\mathbb{R}^D \setminus \widetilde{\mathcal{S}}) ~~\approx~~ H_0 (\widetilde{\mathcal{S}})$. Since $\widetilde{\mathcal{S}}$ has $m$ path-connected components, $S_1,S_2,\cdots,S_m$, by Proposition 2.6 of \cite{Hatcher:AlgTop}, $ H_0(\widetilde{\mathcal{S}}) ~\approx~ \bigoplus_{k=1}^m H_0(S_k)$. But by Proposition \ref{prop:isomorphism-homology-int}, $H_0(S_k) ~\approx~ H_{N-1}(\mathbb{R}^D \setminus S_k) ~\approx~ \mathbb{Z}$. Hence the proposed result follows.

\begin{figure*}
  \centering
  \subfigure[$D=2, N=2$]{
      \label{fig:abstract:a}
      \includegraphics[width=0.3\textwidth, trim=100 100 150 130, clip=true]{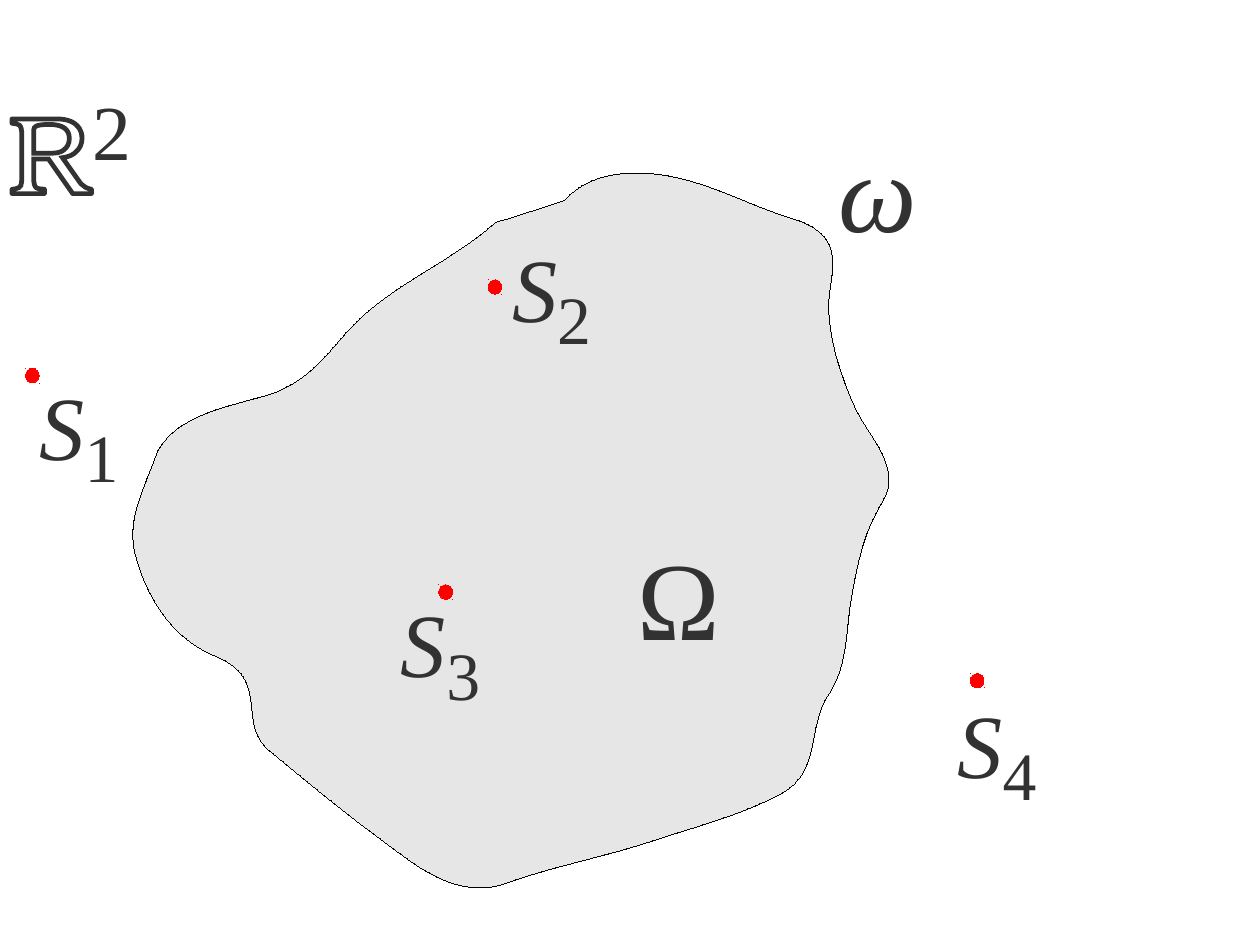}
      }
  \subfigure[$D=3, N=2$]{
      \label{fig:abstract:b}
      \includegraphics[width=0.3\textwidth, trim=50 150 150 60, clip=true]{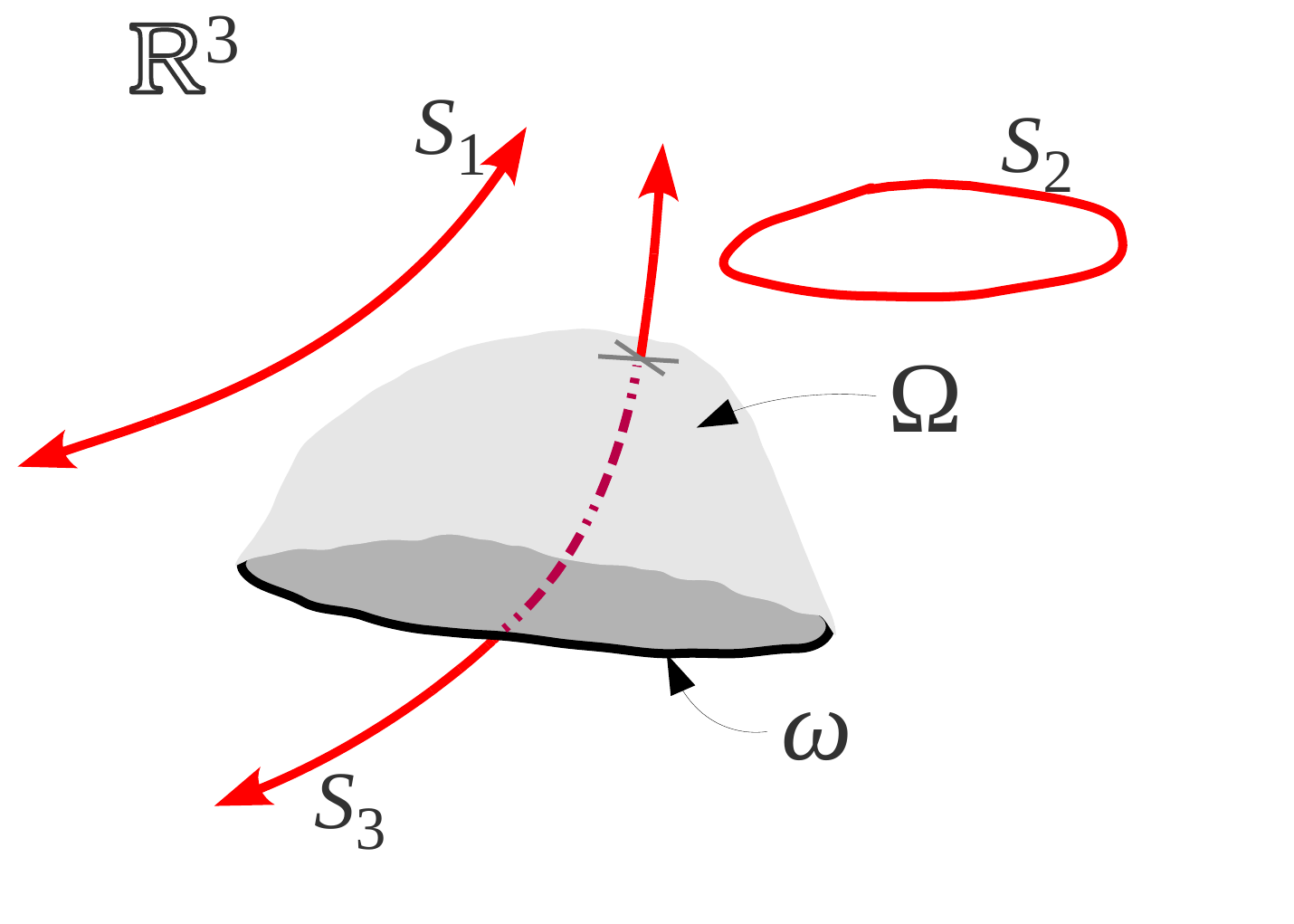}
      }
  \subfigure[$D=3, N=3$]{
      \label{fig:abstract:c}
      \includegraphics[width=0.3\textwidth, trim=50 100 150 130, clip=true]{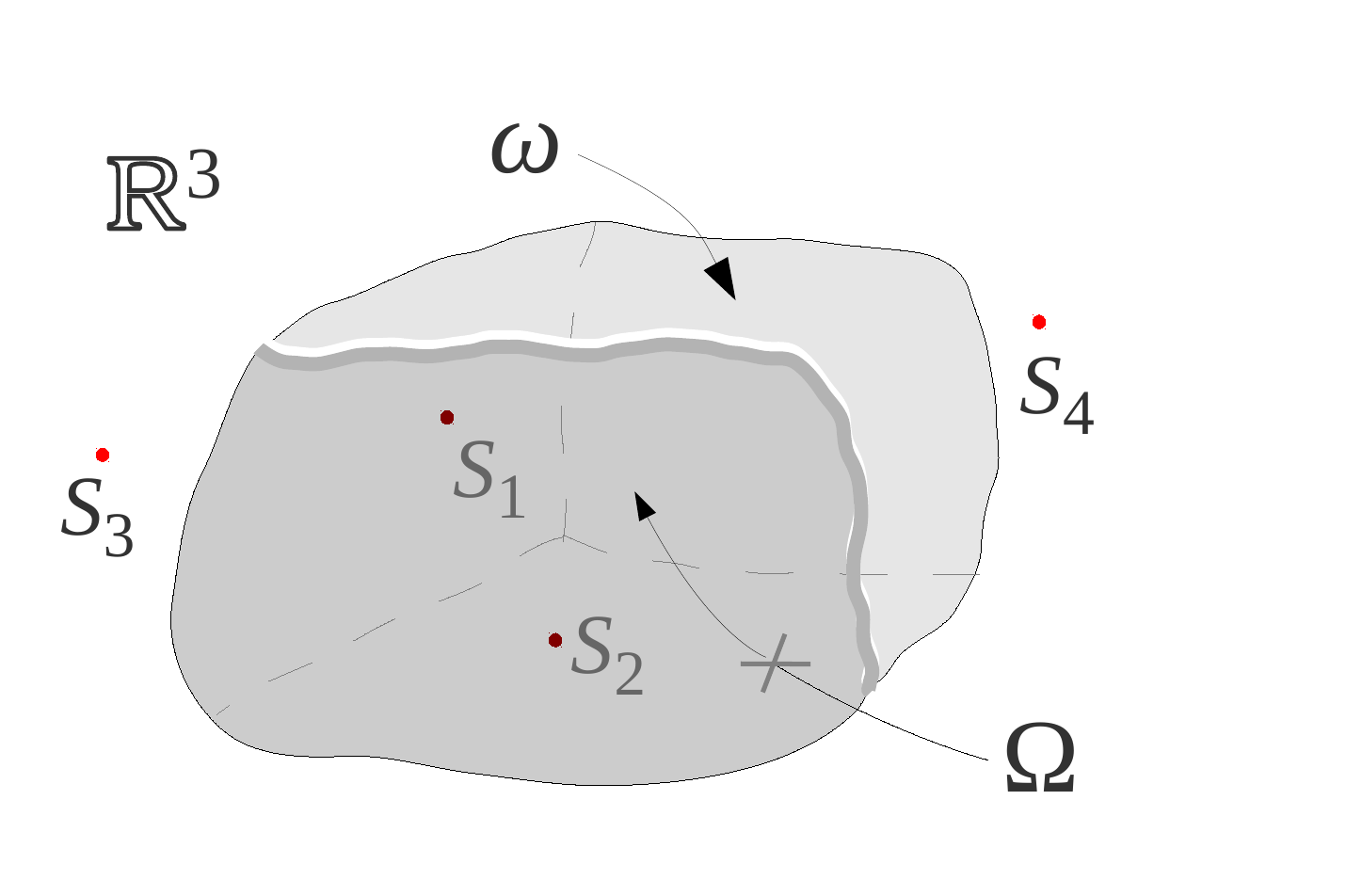}
      }
  \caption[Schematic illustration of some lower dimensional cases of the problem.]{Schematic illustration of some lower dimensional cases of the problem: (a) the Residue theorem, (b) Ampere's law, and (c) Gauss' theorem.} 
  \label{fig:abstract}
\end{figure*}

% -------------------
\Ls{
\section{Validations in Low Dimensions} \label{sec:theo-validation}

In this section we illustrate the forms that equations \eqref{eq:U-final} and \eqref{eq:omg-final} take under certain special cases. We compare those with the well-known formulae from complex analysis, electromagnetism and electrostatics that are known to give homology class invariants.
Once again, we demonstrate all the computations using a single connected component of $\widetilde{\mathcal{S}}$.
}
{We remark that for low values of $D$ and $N$, one recovers well-known integral formulae (Figure~\ref{fig:abstract}). In particular, with $D=2,N=2$ we obtain $\psi_{\mathbf{S}} = \frac{1}{2\pi} ~~Im \left( \frac{1}{z-\mathbf{S}_c} \d z \right)$ --- the differential $1$-form in the Residue Theorem from complex analysis; with $D=3,N=2$ we obtain $\psi_S = \frac{1}{4 \pi} \left( \int_S \frac{\d \mathbf{l'} \times (\mathbf{x}-\mathbf{x'})}{|\mathbf{x}-\mathbf{x'}|^3} \right) \cdot [\d x_1 ~ \d x_2 ~ \d x_3 ]^T$ (where $\d \mathbf{l'} = [\d x'_1 ~ \d x'_2 ~ \d x'_3 ]^T$) --- the differential $1$-form in Ampere's Law; and with $D=3,N=3$ we obtain $\psi_S = \left( \frac{1}{4\pi} \frac{\mathbf{x}-\mathbf{S}}{|\mathbf{x}-\mathbf{S}|^3} \right) ~~\cdot~~ [ ~\d x_2 \wedge \d x_3 ~,~~ \d x_3 \wedge \d x_1 ~,~~ \d x_1 \wedge \d x_2 ]^T$ --- the differential $2$-form in Gauss' divergence theorem.
}

\Lo{
% -------------------------------------------------

\subsection{$D=2,N=2$ :}

This particular case has parallels with the \emph{Cauchy integral theorem} and the \emph{Residue theorem} from \emph{Complex analysis}.
This formula was used in \cite{planning:AURO:12} for designing a $H$-signature in the $2$-dimensional case.
Here a singularity manifold, $S$, is a $D-N=0$-dimensional manifold, \emph{i.e.} a point, the coordinate of which we represent by $\mathbf{S}=[s_1,s_2]^T$ (Figure~\ref{fig:abstract:a}).

\noindent Thus, the partitions in (\ref{eq:omg-final}) for the different values of $k$ are as follows,\\
For $k=1$, ~~${part}^0(\{2\}) = \Big\{ \{\{\},\{2\}\} \Big\} $, \\
For $k=2$, ~~${part}^0(\{1\}) = \Big\{ \{\{\},\{1\}\} \Big\} $

\noindent Thus,
\[
 U^1_1(\mathbf{x}) = \frac{1}{2\pi}~ (-1)^{2-2+1+1} (1) \frac{x_1-S_1}{|\mathbf{x}-\mathbf{S}|^2} = \frac{1}{2\pi}~ \frac{x_1-s_1}{|\mathbf{x}-\mathbf{S}|^2}
\]
\[
 U^2_1(\mathbf{x}) = \frac{1}{2\pi}~ (-1)^{2-2+2+1} (1) \frac{x_2-S_2}{|\mathbf{x}-\mathbf{S}|^2} = - \frac{1}{2\pi}~ \frac{x_2-s_2}{|\mathbf{x}-\mathbf{S}|^2}
\]
where the subscripts of $U$ indicate the index of the partition used (in the lists above). Also, note that integration of a $0$-form on a $0$-dimensional manifold is equivalent to evaluation of the $0$-form at the point.

\noindent Thus,
\begin{eqnarray}
 \psi_{\mathbf{S}} & = & U^1_1(\mathbf{x}) \d x_2 + U^2_1(\mathbf{x}) \d x_1 \nonumber \\
  & = & \frac{1}{2\pi} ~\frac{(x_1-s_1) \d x_2 - (x_2-s_2) \d x_1}{|\mathbf{x}-\mathbf{S}|^2} \nonumber \\
  & = & \frac{1}{2\pi} ~~Im \left( \frac{1}{z-\mathbf{S}_c} \d z \right) \nonumber
\end{eqnarray}
where in the last expression we used the complex variables, $z = x_1 + i x_2$ and $\mathbf{S}_c = s_1 + i s_2$. In fact, from complex analysis (Residue theorem and Cauchy integral theorem) we know that $\int_{\gamma} \frac{1}{z-\mathbf{S}_c} \d z $ (where $\gamma$ is a closed curve in $\mathbb{C}$) is $2\pi i$ if $\gamma$ encloses $\mathbf{S}_c$, but zero otherwise. This is just the fact that
\[
 \int_{\gamma} \psi_{\mathbf{S}} = \int_{\textrm{Ins}(\gamma)} \d \psi_{\mathbf{S}} = \left\{ \begin{array}{l}
                                                                 \pm 1, \textrm{ if $\textrm{Ins}(\gamma)$ contains $S$}\\
                                                                 0, \textrm{ otherwise}
                                                                \end{array} \right.
\]
where $\textrm{Ins}(\gamma)$ represents the inside region of the curve $\gamma$, \emph{i.e.} the area enclosed by it.

% -------------------------------------------------

\subsection{$D=3,N=2$ :}

This particular case has parallels with the \emph{Ampere's Law} and the \emph{Biot-Savart Law} from \emph{Electromagnetism}.
%This formula was used in \cite{planning:AURO:12} for designing a $H$-signature in the $3$-dimensional case.
Here a singularity manifold, $S$, is a $D-N=1$-dimensional manifold, which, in electromagnetics, represents a current-carrying line/wire.% (the \emph{skeletons} according to terminology of \cite{planning:AURO:12} - Figure~\ref{fig:abstract:b}).

\noindent The partitions in (\ref{eq:omg-final}) for the different values of $k$ are as follows,\\
For $k=1$, ~~${part}^1(\{2,3\}) = \Big\{ \{\{2\},\{3\}\} ~,~ \{\{3\},\{2\}\} \Big\} $, \\
For $k=2$, ~~${part}^1(\{1,3\}) = \Big\{ \{\{1\},\{3\}\} ~,~ \{\{3\},\{1\}\} \Big\} $, \\
For $k=3$, ~~${part}^1(\{1,2\}) = \Big\{ \{\{1\},\{2\}\} ~,~ \{\{2\},\{1\}\} \Big\} $, \\

\noindent Thus,
{\small \[
 U^1_1(\mathbf{x}) = \frac{1}{4\pi}~ (-1)^{3-2+1+1} (1) \int_S \frac{x_1-x'_1}{|\mathbf{x}-\mathbf{x'}|^3} \d x'_2 = - \frac{1}{4\pi}~ \int_S \frac{x_1-x'_1}{|\mathbf{x}-\mathbf{x'}|^3} \d x'_2
\]
\[
 U^1_2(\mathbf{x}) = \frac{1}{4\pi}~ (-1)^{3-2+1+1} (-1) \int_S \frac{x_1-x'_1}{|\mathbf{x}-\mathbf{x'}|^3} \d x'_3 = \frac{1}{4\pi}~ \int_S \frac{x_1-x'_1}{|\mathbf{x}-\mathbf{x'}|^3} \d x'_3
\]

\[
 U^2_1(\mathbf{x}) = \frac{1}{4\pi}~ (-1)^{3-2+2+1} (1) \int_S \frac{x_2-x'_2}{|\mathbf{x}-\mathbf{x'}|^3} \d x'_1 = \frac{1}{4\pi}~ \int_S \frac{x_2-x'_2}{|\mathbf{x}-\mathbf{x'}|^3} \d x'_1
\]
\[
 U^2_2(\mathbf{x}) = \frac{1}{4\pi}~ (-1)^{3-2+2+1} (-1) \int_S \frac{x_2-x'_2}{|\mathbf{x}-\mathbf{x'}|^3} \d x'_3 = - \frac{1}{4\pi}~ \int_S \frac{x_2-x'_2}{|\mathbf{x}-\mathbf{x'}|^3} \d x'_3
\]

\[
 U^3_1(\mathbf{x}) = \frac{1}{4\pi}~ (-1)^{3-2+3+1} (1) \int_S \frac{x_3-x'_3}{|\mathbf{x}-\mathbf{x'}|^3} \d x'_1 = - \frac{1}{4\pi}~ \int_S \frac{x_3-x'_3}{|\mathbf{x}-\mathbf{x'}|^3} \d x'_1
\]
\[
 U^3_2(\mathbf{x}) = \frac{1}{4\pi}~ (-1)^{3-2+3+1} (-1) \int_S \frac{x_3-x'_3}{|\mathbf{x}-\mathbf{x'}|^3} \d x'_2 = \frac{1}{4\pi}~ \int_S \frac{x_3-x'_3}{|\mathbf{x}-\mathbf{x'}|^3} \d x'_2
\] }
where, as before, the subscripts of $U$ indicate the index of the partition used (in the lists above).\\
Thus,
\begin{eqnarray}
 \psi_S & = & U^1_1(\mathbf{x}) \d x_3 + U^1_2(\mathbf{x}) \d x_2
              + U^2_1(\mathbf{x}) \d x_3 + U^2_2(\mathbf{x}) \d x_1
              + U^3_1(\mathbf{x}) \d x_2 + U^3_2(\mathbf{x}) \d x_1  \nonumber \\
        & = & ( U^2_2(\mathbf{x}) + U^3_2(\mathbf{x}) ) \d x_1
              + ( U^1_2(\mathbf{x}) + U^3_1(\mathbf{x}) ) \d x_2
              + ( U^1_1(\mathbf{x}) + U^2_1(\mathbf{x}) ) \d x_3  \nonumber \\
        & = & \left[ \begin{array}{c}
                            U^2_2(\mathbf{x}) + U^3_2(\mathbf{x}) \\
                            U^1_2(\mathbf{x}) + U^3_1(\mathbf{x}) \\
                            U^1_1(\mathbf{x}) + U^2_1(\mathbf{x})
                      \end{array} \right] \cdot\wedge
                         \left[ \begin{array}{c}\d x_1 \\ \d x_2 \\ \d x_3 \end{array} \right] \nonumber \\
        & = &  \frac{1}{4 \pi} \int_S \left[ \begin{array}{c}
                               -\frac{x_2-x'_2}{|\mathbf{x}-\mathbf{x'}|^3} \d x'_3 + \frac{x_3-x'_3}{|\mathbf{x}-\mathbf{x'}|^3} \d x'_2 \\
                               \frac{x_1-x'_1}{|\mathbf{x}-\mathbf{x'}|^3} \d x'_3 - \frac{x_3-x'_3}{|\mathbf{x}-\mathbf{x'}|^3} \d x'_1 \\
                               -\frac{x_1-x'_1}{|\mathbf{x}-\mathbf{x'}|^3} \d x'_2 + \frac{x_2-x'_2}{|\mathbf{x}-\mathbf{x'}|^3} \d x'_1
                                                   \end{array} \right] \cdot\wedge
                         \left[ \begin{array}{c}\d x_1 \\ \d x_2 \\ \d x_3 \end{array} \right] \nonumber \\
        & = &  \frac{1}{4 \pi} \int_S \frac{\d \mathbf{l'} \times (\mathbf{x}-\mathbf{x'})}{|\mathbf{x}-\mathbf{x'}|^3} \cdot\wedge
                         \left[ \begin{array}{c}\d x_1 \\ \d x_2 \\ \d x_3 \end{array} \right] \nonumber
\end{eqnarray}
where, bold face indicates column $3$-vectors and the cross product ``$\times$''$:\mathbb{R}^3\times\mathbb{R}^3\to\mathbb{R}^3$ is the elementary cross product operation of column $3$-vectors. The operation ``$\cdot\wedge$'' between column vectors implies element-wise wedge product followed by summation. Also, $\d \mathbf{l'} = [\d x'_1 ~ \d x'_2 ~ \d x'_3 ]^T$.
%The last expression uses the notion of \emph{cross product of column vectors} (which, in this context, is completely independent of the notion of exterior product), and
It is not difficult to identify the integral in the last expression, $ \mathbf{B} = \frac{1}{4 \pi} \int_S \frac{\d \mathbf{l'} \times (\mathbf{x}-\mathbf{x'})}{|\mathbf{x}-\mathbf{x'}|^3} $ with the \emph{Magnetic Field vector} created by unit current flowing through $S$, computed using the \emph{Biot–Savart law}.
Thus, if $\gamma$ is a closed loop, the statement of the \emph{Ampère's circuital law} gives, $\int_{\gamma} \mathbf{B}\cdot\d \mathbf{l} = \int_{\gamma} \psi_S = I_{encl}~$, the \emph{current enclodes by the loop}.

% -------------------------------------------------

\subsection{$D=3,N=3$ :}

This particular case has parallels with the \emph{Gauss's law} in \emph{Electrostatics}, and in general the \emph{Gauss Divergence theorem}.
Here a singularity manifold, $S$, is a $D-N=0$-dimensional manifold, \emph{i.e.} a point, the coordinate of which is represented by $\mathbf{S} = [S_1,S_2,S_3]^T$, which in the light of \emph{Electrostatics}, is a point charge. The candidate manifolds are $2$-dimensional surfaces (Figure~\ref{fig:abstract:c}).

\noindent The partitions in (\ref{eq:omg-final}) for the different values of $k$ are as follows,\\
For $k=1$, ~~${part}^0(\{2,3\}) = \Big\{ \{\{\},\{2,3\}\} \Big\} $, \\
For $k=2$, ~~${part}^0(\{1,3\}) = \Big\{ \{\{\},\{1,3\}\} \Big\} $, \\
For $k=3$, ~~${part}^0(\{1,2\}) = \Big\{ \{\{\},\{1,2\}\} \Big\} $, \\

\noindent Here, $D-N=0$ implies the integration of (\ref{eq:U-final}) once again becomes evaluation of $0$-forms at $\mathbf{S}$. Thus,
{\small \[
 U^1_1(\mathbf{x}) = \frac{1}{4\pi}~ (-1)^{3-3+1+1} (1) \frac{x_1-S_1}{|\mathbf{x}-\mathbf{S}|^3} = \frac{1}{4\pi}~ \frac{x_1-S_1}{|\mathbf{x}-\mathbf{S}|^3}
\]
\[
 U^2_1(\mathbf{x}) = \frac{1}{4\pi}~ (-1)^{3-3+2+1} (1) \frac{x_2-S_2}{|\mathbf{x}-\mathbf{S}|^3} = - \frac{1}{4\pi}~ \frac{x_2-S_2}{|\mathbf{x}-\mathbf{S}|^3}
\]
\[
 U^3_1(\mathbf{x}) = \frac{1}{4\pi}~ (-1)^{3-3+3+1} (1) \frac{x_3-S_3}{|\mathbf{x}-\mathbf{S}|^3} = \frac{1}{4\pi}~ \frac{x_3-S_3}{|\mathbf{x}-\mathbf{S}|^3}
\] }

\noindent Thus,
\begin{eqnarray}
 \psi_S & = & U^1_1(\mathbf{x}) ~\d x_2 \wedge \d x_3  ~~+~~  U^2_1(\mathbf{x}) ~\d x_1 \wedge \d x_3  ~~+~~  U^3_1(\mathbf{x}) ~\d x_1 \wedge \d x_2 \nonumber \\
    & = & \frac{1}{4\pi} \left( \frac{x_1-S_1}{|\mathbf{x}-\mathbf{S}|^3} ~\d x_2 \wedge \d x_3  ~~+~~
                                \frac{x_2-S_2}{|\mathbf{x}-\mathbf{S}|^3} ~\d x_3 \wedge \d x_1  ~~+~~
                                \frac{x_3-S_3}{|\mathbf{x}-\mathbf{S}|^3} ~\d x_1 \wedge \d x_2  ~~+~~  \right) \nonumber \\
    & = & \left( \frac{1}{4\pi} \frac{\mathbf{x}-\mathbf{S}}{|\mathbf{x}-\mathbf{S}|^3} \right) ~~\cdot\wedge~~ [ ~\d x_2 \wedge \d x_3 ~,~~ \d x_3 \wedge \d x_1 ~,~~ \d x_1 \wedge \d x_2 ]^T
\end{eqnarray}
The quantity $\mathbf{E} = \frac{1}{4\pi} \frac{\mathbf{x}-\mathbf{S}}{|\mathbf{x}-\mathbf{S}|^3}$ can be readily identified with the electric field created by an unit point charge at $\mathbf{S}$. If $\mathcal{A}$ is a closed surface, then $\int_{\mathcal{A}} \mathbf{E}\cdot \d \mathbf{A} = \int_{\mathcal{A}} \psi_S = Q_{encl}~$, the charge enclosed by $\mathcal{A}$.
}

% ====================================================

\section{Examples and Applications} \label{sec:example-applications}

We implemented the general formula for computing $\psi_{\mathcal{S}}(\omega)$ in {\tt C++} for arbitrary $D$ and $N$.
The singularity manifolds, $S$, and the candidate manifold, $\omega$, are discretized to create simplicial complexes $\overline{S}$ and $\overline{\omega}$ respectively, thus enabling us to compute the integral in equations \eqref{eq:int-eta-final} and \eqref{eq:U-final} as a sum of integrals over simplices. In the following section, for simplicity, we use the same notation for the manifolds and their simplicial equivalents.
We used the {\em Armadillo linear programming library} \cite{Armadillo} for all vector and matrix operations, and the {\em GNU Scientific Library} \cite{GSL:Manual} for all the numerical integrations.

\subsection{An Example for $D=5,N=3$} \label{sec:5d-validation}

\Lo{In Section~\ref{sec:theo-validation} we have shown that the general formulation we proposed in Section~\ref{sec:homology-invariant-formula} indeed reduces to known formulae that gives us the homology class invariants for certain low dimensional cases.} We present numerical validation for a simple case of dimension greater than three: $D=5$ and $N=3$. The candidate manifold is of dimension $N-1=2$. We consider a $2$-sphere centered at the origin in $\mathbb{R}^5$ as the candidate manifold: let $\omega(R_C) = \{ \mathbf{x}~~ | ~~x_1^2 + x_2^2 + x_3^2 = R_C^2, ~~x_4=0, ~~x_5=0 \}$ be the boundary of the ball
$\Omega(R_C) = \{ \mathbf{x}~~ | ~~x_1^2 + x_2^2 + x_3^2 \leq R_C^2, ~~x_4=0, ~~x_5=0 \}$.
The candidate manifold  $\omega(R_C)$ is easily parametrized  \lo{(see Figure~\ref{fig:sphere-triangulation})} via spherical coordinates $\theta$ and $\phi$. %$\theta\in[-\frac{\pi}{2},\frac{\pi}{2}]$ and $\phi\in[0,2\pi]$, as follows:
\Ls{\begin{equation}
\begin{array}{l}
 x_1 = R_C \cos(\theta)\cos(\phi) \\
 x_2 = R_C \cos(\theta)\sin(\phi) \\
 x_3 = R_C \sin(\theta) \\
 x_4 = 0 \\
 x_5 = 0
\end{array} \label{eq:num-an-omega}
\end{equation}}{$x_1 \!=\! R_C \cos(\theta)\cos(\phi), x_2 \!=\! R_C \cos(\theta)\sin(\phi), x_3 \!=\! R_C \sin(\theta), x_4 \!=\! 0, x_5 \!=\! 0$. }
Let the singularity manifold $S$ be the $2$-torus \lo{(Figure~\ref{fig:torus-triangulation}) }as follows:
\Ls{\begin{equation}
 \begin{array}{l}
  x_1 = 0 \\
  x_2 = 0 \\
  x_3 = \left( R_T + r \cos(\phi') \right) \cos(\theta') - (R_T + r) \\
  x_4 = \left( R_T + r \cos(\phi') \right) \sin(\theta') \\
  x_5 = r \sin(\phi)
 \end{array} \label{eq:num-an-S}
\end{equation}}{$x_1 \!=\! 0, x_2 \!=\! 0, x_3 \!=\! \left( R_T + r \cos(\phi') \right) \cos(\theta') - (R_T + r), x_4 \!=\! \left( R_T + r \cos(\phi') \right) \sin(\theta'), x_5 \!=\! r \sin(\phi)$, }
with $R_T > r$ and the parameters $\theta'\in[0,2\pi]$ and $\phi'\in[0,2\pi]$. For all examples that follow, we choose $r=0.8, R_T=1.6$.

\lo{
\begin{figure*}
  \centering
  \subfigure[Triangulation of the singularity manifold projected on the space of $x_3, x_4, x_5$.]{
      \label{fig:torus-triangulation}
      \includegraphics[width=0.4\textwidth, trim=100 200 100 200, clip=true]{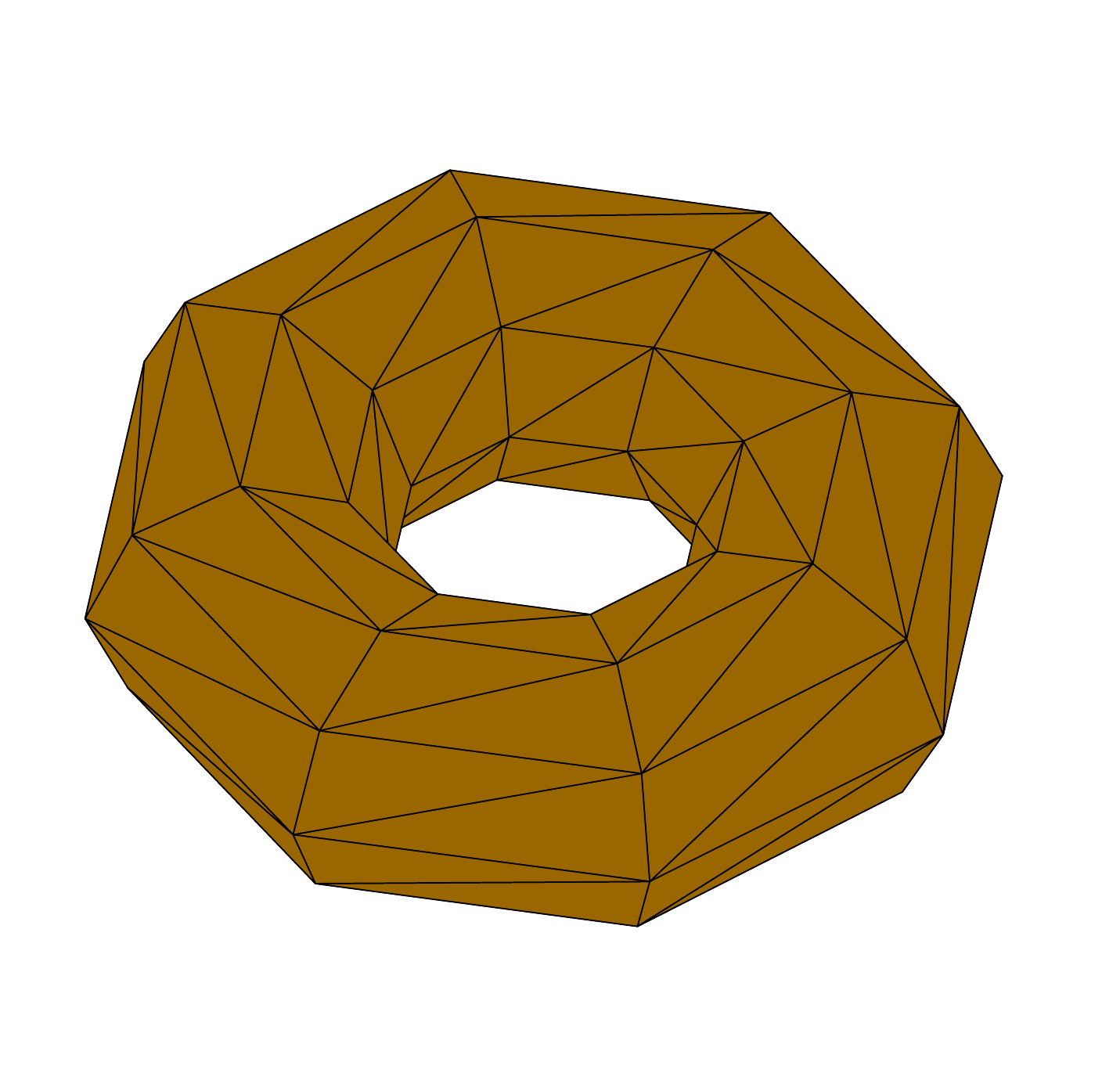}
      } \hspace{0.02in}
  \subfigure[Triangulation of a candidate manifold projected on the space of $x_1, x_2, x_3$.]{
      \label{fig:sphere-triangulation}
      \includegraphics[width=0.4\textwidth, trim=100 200 100 200, clip=true]{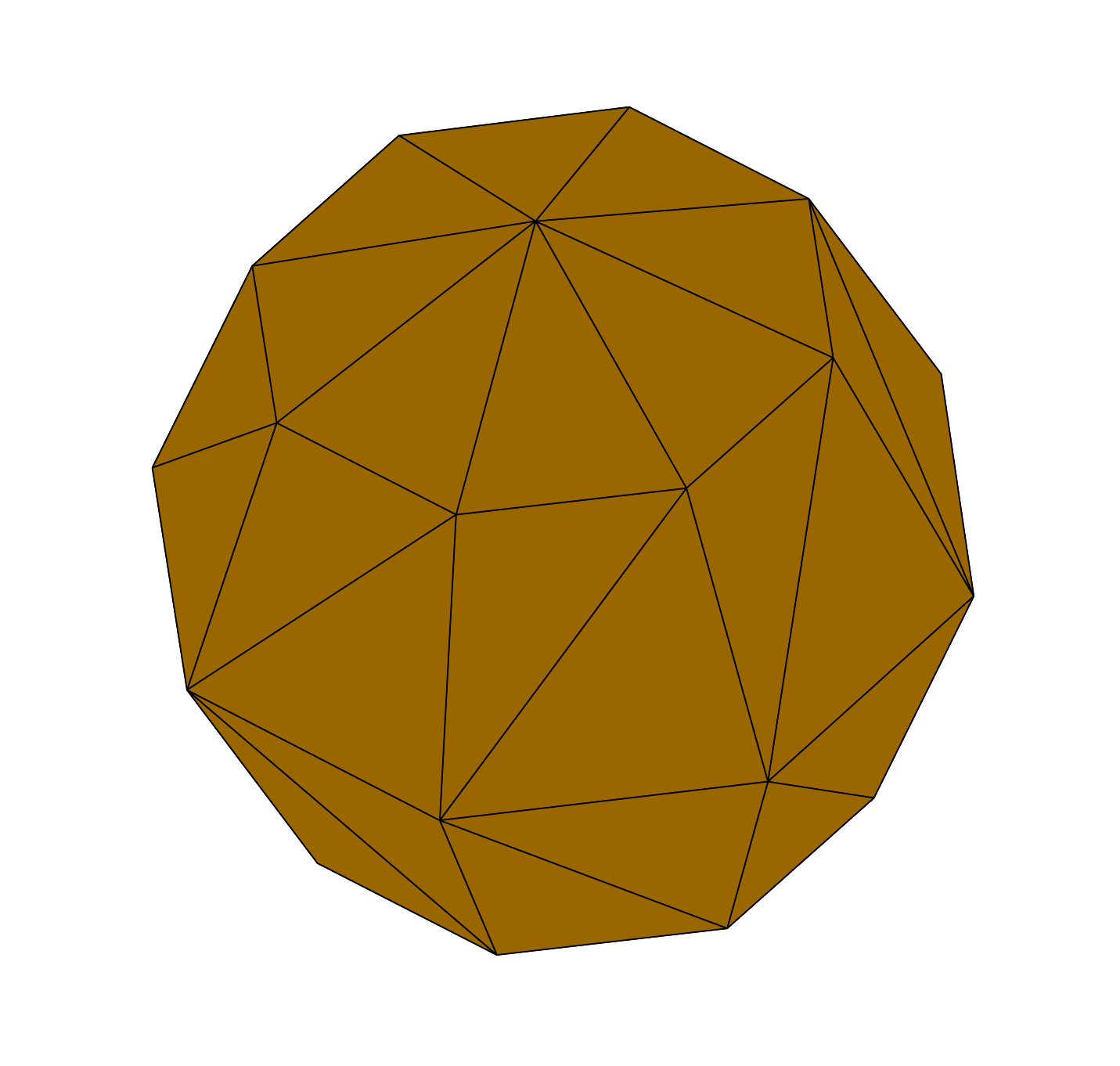}
      }
  \caption{A coarse triangulation using parameters $\theta'$, $\phi'$, $\theta$ and $\phi$ for creating a simplicial complex for the example in Section~\ref{sec:5d-validation}.}
  \label{fig:D5N3-example}
\end{figure*}}

Consider the particular candidate manifold $\omega(1)$ (\emph{i.e.} $R_C=1$). By numerical computation of integrals in \eqref{eq:int-eta-final} and \eqref{eq:U-final}, the value of $\phi_S(\omega(1))$ that we obtain for the above example is $-1$.
In order to interpret this result we first observe that $\omega(1)$ does not intersect $S$\Lo{~(\emph{i.e.} there is no common solution for (\ref{eq:num-an-omega}) and (\ref{eq:num-an-S}) with $R_C=1.0,r=0.8, R_T=1.6$)}. However on $S$\lo{~(Equations \eqref{eq:num-an-S})}, when $x_1=x_2=x_4=x_5=0$, $x_3$ can assume the values $0$, $-2r$, $-2R_T$ and $-2(R_T + r)$. Thus, if $2r > R_C$, $S$ intersects $\Omega(R_C)$ (the ball whose boundary is $\omega(R_C)$) only at one point, the origin. A simple computation of the tangents reveals that the intersection is transverse. Since that is a single transverse intersection with $\Omega(R_C)$, the linking number between $\omega(R_C)$ and $S$ (\emph{i.e.} intersection number between $\Omega(R_C)$ and $S$ according to Definition~\ref{def:linking-number}) is $\pm 1$ for all $R_C<2r$, just as indicated by the value of $\phi_S(\omega(1))$. The sign is not of importance since that is determined by our choice of orientation. In fact, with different values of $R_C, r$ and $R_T$, as long as $R_T > r > \frac{R_C}{2}$, we obtain the same value of $-1$ for $\phi_S(\omega(R_C))$.
%So do we obtain by perturbation of the pose and deformation of the sphere or torus.

However with $R_C=2$ for the candidate manifold, and the singularity manifold remaining the same (\emph{i.e.} $r=0.8, R_T=1.6$), the value of $\phi_S(\omega(2))$ we obtain numerically is $0$. In this case, the points at which $S$ intersect $\Omega(2)$ are the origin and the point $(x_1=x_2=x_4=x_5=0,x_3=-0.8)$. Of course, in the family of candidate manifolds $\omega(R_C),~R_C\in[1,2]$, we can easily observe that $\omega(1.6)$ indeed intersects $S$, thus indicating that $\omega(1)$ and $\omega(2)$ are possibly in different homology classes.

Next, consider the following family of candidate manifolds:
\Ls{\begin{equation}
 \omega'(T_C) = \{ \mathbf{x}~~ | ~~x_1^2 + x_2^2 + x_3^2 = 2.0, ~~x_4=0, ~~x_5=T_C \}
\end{equation}}{$\omega'(T_C) = \{ \mathbf{x}~~ | ~~x_1^2 + x_2^2 + x_3^2 = 2, ~~x_4=0, ~~x_5=T_C \}$, }
and a corresponding $\Omega'(T_C)$ such that $\omega'(T_C) = \partial \Omega'(T_C)$:
\Ls{\begin{equation}
 \Omega'(T_C) = \{ \mathbf{x}~~ | ~~x_1^2 + x_2^2 + x_3^2 \leq 2, ~~x_4=0, ~~x_5=T_C \}
\end{equation}}{$\Omega'(T_C) = \{ \mathbf{x}~~ | ~~x_1^2 + x_2^2 + x_3^2 \leq 2.0, ~~x_4=0, ~~x_5=T_C \}$.}
With the same $S$ as before, if $T_C>r$, clearly there is no intersection between $\Omega'(T_C)$ and $S$. Thus it is not surprising that indeed by numerical computation, we found that $\phi_S(\omega'(1)) = 0$.

Now, since we computed $\phi_S(\omega(2)) = 0$ (although $\Omega(2)$ intersects $S$ at $2$ points) and $\phi_S(\omega'(1)) = 0$ (and $\Omega'(1)$ does not intersect $S$), it suggests that $\omega(2)$ and $\omega'(1)$ are in the same homology class. We verify this by observation.
%find a path in $\mathbf{\varpi}^2_5$ that connects $\chi_S(\omega(2))$ and $\chi(\omega'(1))$.
None from the family of candidate manifolds $\omega'(T_C),~\forall T_C\in[0,1]$ intersect $S$, and each is a $2$-sphere.
Thus $\omega'$ defines an embedding of $\mathbb{S}^2\times I$ in $\mathbb{R}^5 - S$ such that $\omega'(0)\sqcup -\omega'(1)$ is its boundary. It follows that $\omega'(0)$ and $\omega'(1)$ are homologous.
%Thus, there is a path in $\mathbf{\varpi}^2_5$ connecting $\omega'(0)$ and $\omega'(1)$, \emph{i.e.} they are $\chi$-homotopic by definition.
However, $\omega(2)=\omega'(0)$. Thus it follows that $\omega(2)$ and $\omega'(1)$ are homologous.

\subsection{Application to Graph Search-based Robot Path Planning with Topological Constraints}

One consequence of $\phi_{\widetilde{\mathcal{S}}}$ being a cocycle is that it is a linear function.
%This is evident from the fact that $\phi_{\widetilde{\mathcal{S}}}$ was described by direct sum of integrations over the cycles in Equations~\eqref{eq:int-eta-final} and \eqref{eq:final-psi-eqn}.
As a result, if we have a cycle  $\overline{\omega}$ that can be expressed as a sum of chains, \emph{i.e.} $\overline{\omega} = \sum_i \overline{\tau}_i$, with $\overline{\tau}_i \in C_{N-1}(\mathbb{R}^D \setminus \widetilde{\mathcal{S}})$, then we can write
\begin{equation} 
 \phi_{\widetilde{\mathcal{S}}}(\overline{\omega}) = \sum_i \phi_{\widetilde{\mathcal{S}}}(\overline{\tau}_i)
\end{equation}
where by $\phi_{\widetilde{\mathcal{S}}}(\overline{\tau}_i)$ we simply mean the vector formed by evaluation of the integrals in Equations~\eqref{eq:final-psi-eqn}.

%The following is the obvious definition, included for completeness, motivated by path-planning relative to endpoints.
% RG WE SHOULD NOT BE DEFINING THIS...AAARGH
\begin{remark} %[Homology classes of chains with the same boundary] 
\label{def:homology-chains}
%\textbf{Remark:} 
Given $(N-1)$-chains, $\overline{\tau}_1$ and $\overline{\tau}_2$ in \changedA{$X$}, such that $\partial \overline{\tau}_1 = \partial \overline{\tau}_2$, \changedC{by an abuse of terminology} \changedC{in the following discussions, we will say} that they are in the same homology class if $\overline{\tau}_1 -\overline{\tau}_2$ is null-homologous in \changedA{$X$}. \changedC{It should however be remembered that homology classes are not formally defined for chains, and are defined only for cycles or relative cycles.}
\end{remark}

That is, \changedA{in the context of our problem where $X=(\mathbb{R}^D \setminus \widetilde{\mathcal{O}})$,} $\overline{\tau}_1 \approx \overline{\tau}_2$ iff $\phi_{\widetilde{\mathcal{S}}}(\overline{\tau}_1 -\overline{\tau}_2) = 0$ (where $\widetilde{\mathcal{S}}$ is the equivalent of $\widetilde{\mathcal{O}}$ satisfying the property of Proposition~\ref{prop:obstacle-equivalent}).
% This is an equivalence relation since by the property of the integration in the definition of $\phi_{\widetilde{\mathcal{S}}}$,
%\begin{equation}
% \phi_{\widetilde{\mathcal{S}}}(\overline{\tau}_1 -\overline{\tau}_2) = 0 ~~~~\Leftrightarrow~~~~ \phi_{\widetilde{\mathcal{S}}}(\overline{\tau}_1) = \phi_{\widetilde{\mathcal{S}}}(\overline{\tau}_2)
%\end{equation}
In context of robot path planning problem, the candidate manifolds are all $1$-dimensional. Thus we have $N=2$.
While trajectories connecting two points in a configuration space $(\mathbb{R}^D \setminus \widetilde{\mathcal{O}})$ themselves are not closed manifolds, two trajectories connecting the same points together form a closed manifold. %Thus for homology classes of trajectories, we will use Definition~\ref{def:homology-chains}.

\begin{figure*}
\centering
 \subfigure[A graph created by uniform square discretization of an environment. The \changedA{\ls{brown}{dark}} cells represent obstacles. Each vertex is connected to its $8$ neighbors (except inaccessible vertices).]{ \label{fig:graph-8-connected}
  \includegraphics[width=0.35\textwidth, trim=120 280 220 250, clip=true]{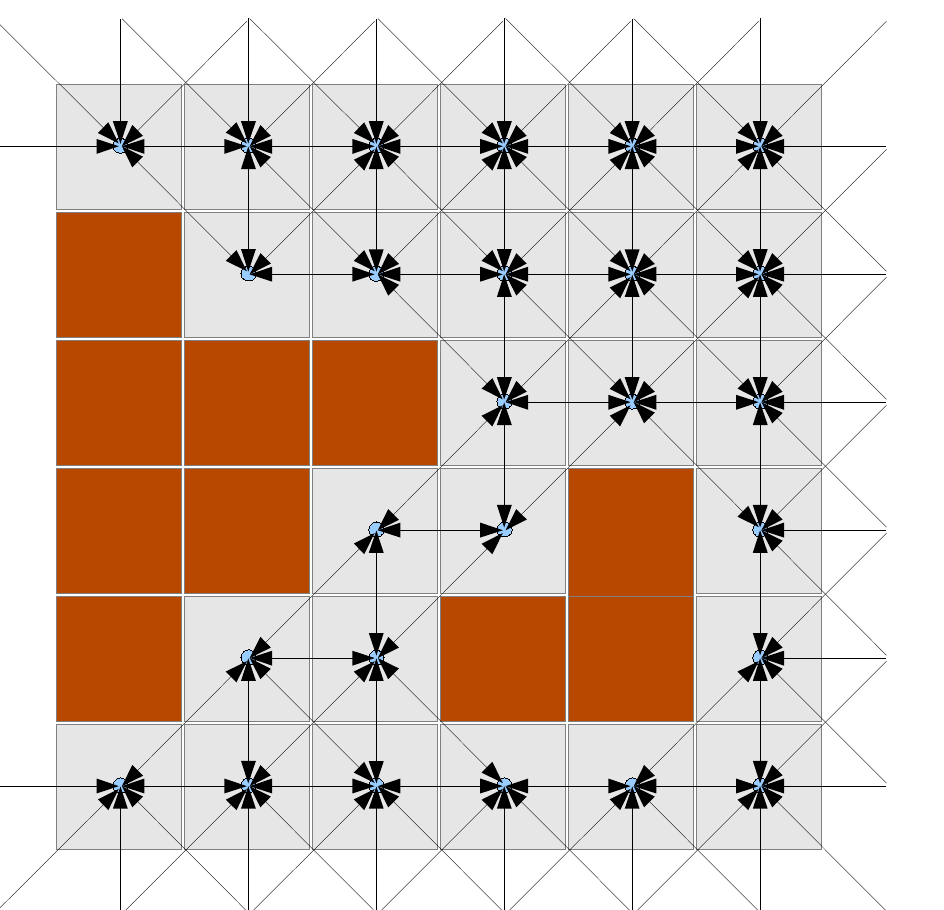}
 } \hspace{0.2in}
 \subfigure[A trajectory in the continuous configuration space can be approximated by a path in the graph.]{ \label{fig:graph-path}
  \includegraphics[width=0.35\textwidth, trim=220 120 220 160, clip=true]{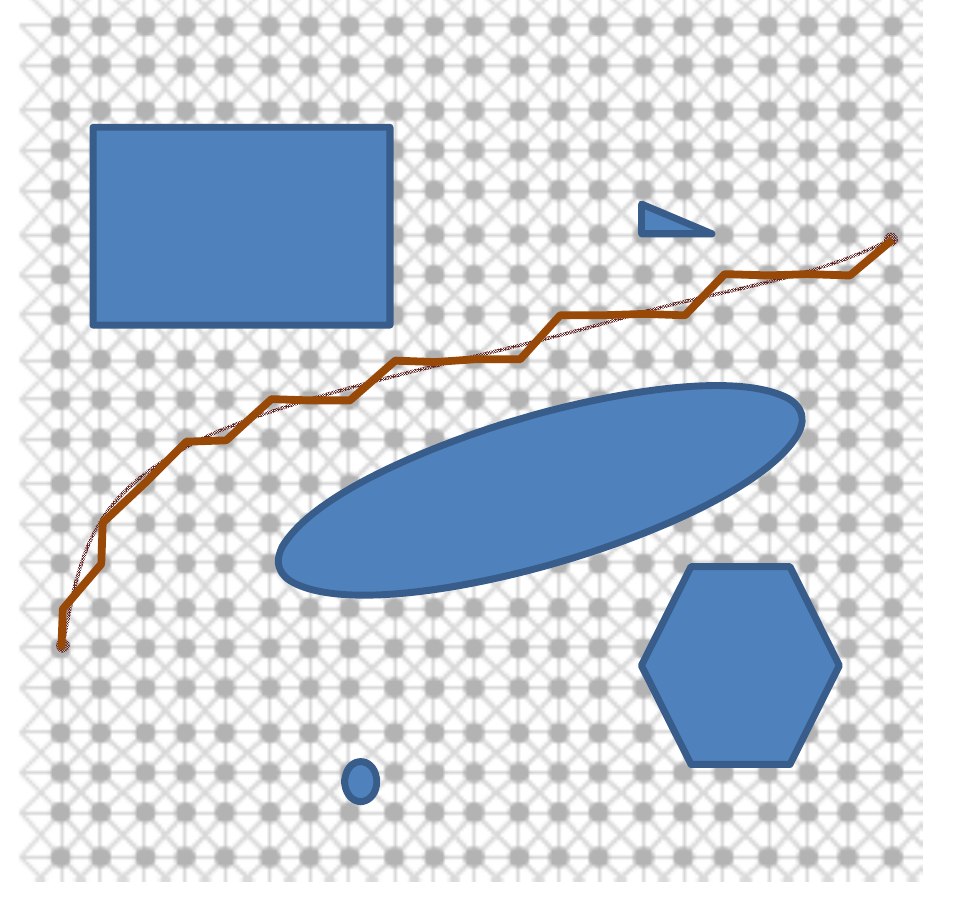}
 }
 \caption{Graph, $\mathcal{G}$, created by uniform discretization of an environment. This specific type of graph shown in the figures is referred to as the $8$-connected grid.} 
\label{fig:discrete}
\end{figure*}

\vspace{0.1in}
Next we outline the basic graph construction for search-based planning with topological constraints ({\em cf.} the $H$-augmented graph of \cite{planning:AURO:12}).
Discrete graph search techniques for robot path planning problems are widely used and have been shown to be complete and efficient \cite{Ste95b,max:planning:08}.
Given a $D$-dimensional configuration space, the standard starting point is to discretize the configuration space, place vertices inside each discrete cell, and establish edges between the neighboring vertices to create a directed graph, $\mathcal{G}=(\mathcal{V},\mathcal{E})$ (Figure~\ref{fig:graph-8-connected}). The discretization itself can be quite arbitrary and non-uniform in general.
% RG DON"T NEED THIS 
%The vertices $\mathbf{v}=[x_1,x_2,\cdots,x_D]\in\mathcal{V}$ represents the coordinate of the centroid of a discretized cell, and 
A directed edge $[\mathbf{v}_1, \mathbf{v}_2]\in\mathcal{E}$ connects vertices $\mathbf{v}_1$ to $\mathbf{v}_2$ iff
there is a single action of the robot that can take it from state $\mathbf{v}_1$ to state $\mathbf{v}_2$.
Since an edge $[\mathbf{v}_1, \mathbf{v}_2]\in\mathcal{E}$ is a $1$-dimensional manifold embedded in $(\mathbb{R}^D \setminus \widetilde{\mathcal{S}})$, we can evaluate the function $\phi_{\widetilde{\mathcal{S}}}$ on (a top-dimensional covering chain on) it we and write it as $\phi_{\widetilde{\mathcal{S}}}([\mathbf{v}_1, \mathbf{v}_2])$.
Likewise, a path, $\lambda$, in the graph (Figure~\ref{fig:graph-path}) can be represented by a covering chain $\overline{\lambda} \in H_{N-1}(\mathbb{R}^D \setminus \widetilde{\mathcal{S}})$, and $\phi_{\widetilde{\mathcal{S}}}$ can be evaluated on it. For simplicity, we often write $\phi_{\widetilde{\mathcal{S}}}(\lambda)$ to indicate this quantity, which is made possible due to the assumption that such covering chains are essentially constructed out of simplices with unit coefficients.
The weight/cost of each edge is the cost of traversing that edge by the robot (typically the metric length of the edge). We write $w([\mathbf{v}_1, \mathbf{v}_2])$ to represent the weight of an edge.
Inaccessible coordinates (lying inside obstacles or outside a specified workspace) do not constitute nodes of the graph. A path in this graph represents a trajectory of the robot in the configuration space. The triangulation of any path in the graph essentially consists of the directed edges of the graph that make up the path.

Suppose we are given a fixed start and a fixed goal coordinate, represented by $\mathbf{v}_s,\mathbf{v}_g\in(\mathbb{R}^D \setminus \widetilde{\mathcal{O}})$ respectively, for the robot (by the boldface $\mathbf{v}$'s, with a slight abuse of notation, we will indicate both the vertex in the graph as well as the coordinate of the vertex in the original configuration space).
%These two points together form the boundary of any $N-1=1$-dimensional trajectory of a robot (see Figures~\ref{fig:robot-trajs}). In accordance to our previous discussion, those points form the $N-2=0$-dimensional reference manifold, $\lambda_{sg}$. That is,
%\[ \lambda_{sg} = \mathbf{v}_s \sqcup \mathbf{v}_g \]
We next construct an augmented graph, $\widehat{\mathcal{G}}= \{ \widehat{\mathcal{V}}, \widehat{\mathcal{E}}$, from the graph $\mathcal{G}$ in order to incorporate the information regarding the homology class of trajectories leading from the given start coordinate to the goal coordinate, as follows.
\begin{itemize} 
 %\vspace{0.9in}
 \item[1.]  \[ 
     \widehat{\mathcal{V}} = \left\{ \{\mathbf{v},\mathbf{c}\} \left|
                       \begin{array}{l}
                        \mathbf{v} \in \mathcal{V}, \textrm{ and,} \\
                        \mathbf{c} \textrm{ is a $m$-vector of reals such that } \mathbf{c} = \phi_{\widetilde{\mathcal{S}}}(\overline{\lambda}) \\ %\textrm{ for some $1$ dimensional }\\
                                   \quad \textrm{ for some $1$-chain, $\overline{\lambda}$, with boundary } \mathbf{v}_s \sqcup -\mathbf{v} \\
									\quad\quad \textrm{(\emph{i.e.} $\overline{\lambda}$ is a covering chain of some path in $\mathcal{G}$ connecting $\mathbf{v}_s$ to $\mathbf{v}$)}. %, \textrm{ and,} \\
%                        \mathbf{h} \in \mathcal{A} ~(\textrm{equivalently, } \mathbf{h} \notin \mathcal{B}) \\
%                                   \qquad\qquad\qquad \textrm{ when } \mathbf{v}=\mathbf{v}_g
                       \end{array}
                     \right. \!\!\!\right\}
    \]

 \item[2.] An edge $[ \{\mathbf{v},\mathbf{c}\}, \{\mathbf{v'},\mathbf{c}'\} ]$ exists in $\widehat{\mathcal{E}}$ for $[\mathbf{v},\mathbf{c}] \in \widehat{\mathcal{V}}$ and $[\mathbf{v'},\mathbf{c}'] \in \widehat{\mathcal{V}}$, iff
   \begin{itemize}
    \item [i.] The edge $[\mathbf{v}, \mathbf{v'}] \in \mathcal{E}$, and,
    \item[ii.] $\mathbf{c}' = \mathbf{c} + \phi_{\widetilde{\mathcal{S}}}([\mathbf{v}, \mathbf{v'}])$.
               % where, $\mathcal{H}(\mathbf{v} \rightarrow \mathbf{v'})$ is the \emph{homotopy signature} of the edge $\{\mathbf{v} \rightarrow \mathbf{v'}\} \in \mathcal{E}$.
   \end{itemize}

 \item[3.] The cost/weight associated with an edge $[ \{\mathbf{v},\mathbf{c}\}, \{\mathbf{v'},\mathbf{c}'\} ]$ is same as the cost/weight associated with edge $[\mathbf{v}, \mathbf{v'}] \in \mathcal{E}$. That is, the weight function we use is $\widehat{w}([ \{\mathbf{v},\mathbf{c}\}, \{\mathbf{v'},\mathbf{c}'\} ]) = w([\mathbf{v}, \mathbf{v}'])$.

\end{itemize}
It can be noted that $\{\mathbf{v}_s,\mathbf{0}\}$ is in $\widehat{\mathcal{V}}$ (where $\mathbf{0}$ is an $m$-vector of zeros).

For finding a least cost path in $\widehat{\mathcal{G}}$ that belongs to a particular homotopy class, we can use a heuristic graph search algorithm (\emph{e.g.} weighted A*) \cite{Hart-Astar,cormen2001,KoeLik-DLite}. In particular, we used the YAGSBPL library \cite{yagsbpl} for constructing the graph and performing A* searches in it. Starting from the start vertex $\{\mathbf{v}_s,\mathbf{0}\}$ we expand the vertices in $\widehat{\mathcal{G}}$. The process of vertex expansion eventually leads to vertices of the form $\{\mathbf{v}_g,\mathbf{c}_i\}$, where $\mathbf{c}_i=\phi_{\widetilde{\mathcal{S}}}(\lambda_{sg})$ for some path $\lambda_{sg}$ in $\mathcal{G}$ connecting $\mathbf{v}_s$ to $\mathbf{v}_g$.
Each of these vertices in $\widehat{\mathcal{G}}$ correspond to an unique homology class of the path taken to reach $\mathbf{v}_g$ from $\mathbf{v}_s$.
Let those vertices in the order in which we expand them be $\{\mathbf{v}_g,\mathbf{c}_1\}$, $\{\mathbf{v}_g,\mathbf{c}_2\}$, etc. Say during the search process, we expand the vertex $\{\mathbf{v}_g,\mathbf{c}_j\}\in\widehat{\mathcal{V}}$.
Depending on whether we are trying to search for a particular homology class of trajectories or exploring multiple homology classes, %the particular problem
we can choose to take one of the following actions:
%If $\mathbf{c}_1$ is a desired or an admitted value for the $\chi$-value of the trajectory we are searching for, we have two options:
\begin{itemize}
 \item[i.] If $\mathbf{c}_j$ is the desired value (or an admitted value) for the $\phi_{\widetilde{\mathcal{S}}}$-value of the trajectory we are searching for, we store the path up to $\{\mathbf{v}_g,\mathbf{c}_j\}$ in $\widehat{\mathcal{G}}$, and stop the search algorithm.
 \item[ii.] If $\mathbf{c}_j$ is an admitted value for the $\phi_{\widetilde{\mathcal{S}}}$-value of the trajectory we are searching for, we store the path up to $\{\mathbf{v}_g,\mathbf{c}_j\}$ in $\widehat{\mathcal{G}}$, and continue expanding vertices in $\widehat{\mathcal{G}}$.
 \item[iii.] If $\mathbf{c}_j$ is not an admitted value for the $\phi_{\widetilde{\mathcal{S}}}$-value of the trajectory we are searching for, we continue expanding vertices in $\widehat{\mathcal{G}}$.
% \item[ii] We can store the path up to $\{\mathbf{v}_g,\mathbf{c}_1\}$, and continue expanding states in $\mathcal{G}_\chi$, which will eventually lead us to the expansion of another node, $\{\mathbf{v}_g,\mathbf{c}_2\}$. Then, once again we can store the path to $\{\mathbf{v}_g,\mathbf{c}_2\}$ and stop the search or we can continue expansion of states.
\end{itemize}
Clearly, the projection of any of the stored trajectories onto $\mathcal{G}$ are paths in $\mathcal{G}$ connecting $\mathbf{v}_s$ to $\mathbf{v}_g$.
Since both $\widehat{\mathcal{G}}$ and $\mathcal{G}$ use the same cost function, if $\left\{ \{\mathbf{v}_s,\mathbf{0}\}, \{\mathbf{v}^{1*},\mathbf{c}^{1*}\}, \{\mathbf{v}^{2*},\mathbf{c}^{2*}\}, \cdots, \{\mathbf{v}_g,\mathbf{c}_j\} \right\}$ is the $j^{th}$ stored path using an optimal search algorithm (\emph{e.g} A*), then $\left\{ \mathbf{v}_s,\mathbf{v}^{1*}, \mathbf{v}^{2*}, \cdots, \mathbf{v}_g \right\}$ is the optimal path in $\mathcal{G}$ with $\phi_{\widetilde{\mathcal{S}}}$-value of $\mathbf{c}_j$ (\emph{i.e.} least cost path belonging to the particular homology class).
Thus we can explore the different homology classes of the trajectories connecting $\mathbf{v}_s$ to $\mathbf{v}_g$.

%If we know the exact goal node in $\mathcal{G}_\chi$ that we want to reach, say $\{\mathbf{v}_g,\mathbf{c}_g\}$,
If $\mathbf{c}_g$ is the desired value of $\phi_{\widetilde{\mathcal{S}}}$ evaluated on the trajectory we are searching for,
we follow the above process of expanding the vertices using the graph search algorithm until we expand $\{\mathbf{v}_g,\mathbf{c}_g\}$.
Given two paths $\lambda_1,\lambda_2$ in $\mathcal{G}$, and if $\overline{\lambda}_1,\overline{\lambda}_2$ are their respective covering chains, since $\overline{\lambda}_1 \sqcup -\overline{\lambda}_2 \in C_{N-1}(\mathbb{R}^D \setminus \widetilde{\mathcal{S}})$, we notice that $( \phi_{\widetilde{\mathcal{S}}}(\overline{\lambda}_1) - \phi_{\widetilde{\mathcal{S}}}(\overline{\lambda}_1) ) \in \mathbb{Z}^m$ (with unit coefficients on the simplices that constitute the chains, and with the choice of $\phi_{\widetilde{\mathcal{S}}}$ as described in Equations~\eqref{eq:int-eta-final} and \eqref{eq:final-psi-eqn}). Thus, if we know the value of a $\mathbf{c}_j = \chi_{\widetilde{\mathcal{S}}}(\overline{\lambda}_j)$, we can construct another $m$-vector that is a valid value for $\phi_{\widetilde{\mathcal{S}}}$ evaluated on some other trajectory connecting $\mathbf{v}_s$ to $\mathbf{v}_g$ as $\mathbf{c}_{j'} = \mathbf{c}_j + \zeta$ for some $\zeta \in \mathbb{Z}^m$. This we can hence set as $\mathbf{c}_g$ for finding the least cost path in the new homology class.

%\subsubsection{Cost and Heuristics function}

A consequence of the point $3$ in the definition of $\mathcal{G}_\chi$ is that any \emph{admissible heuristics} (which is a lower bound on the cost to the goal vertex)
%for any heuristic graph search algorithm (\emph{e.g.} weighted A*)
in $\mathcal{G}$ will remain admissible in $\widehat{\mathcal{G}}$. That is, if $h(\mathbf{v},\mathbf{v}')$ was the heuristic function in $\mathcal{G}$, we can define $\widehat{h}(\{\mathbf{v},\mathbf{c}\},\{\mathbf{v}',\mathbf{c}'\}) = h(\mathbf{v},\mathbf{v}')$ as the heuristic function in $\widehat{\mathcal{G}}$. % for any $\mathbf{c}'\in\chi_{\widetilde{\mathcal{S}}}(\Lambda;\lambda_{sg})$.
As a consequence, if we keep expanding vertices in $\widehat{\mathcal{G}}$ as described in the previous section, the order in which we will encounter states of the form $\{\mathbf{v}_g,\mathbf{c}_i\}$ is the order of the costs of the least cost paths in the different homology classes.
%connecting  $\{\mathbf{v}_s,\mathbf{0}\}$ to  $\{\mathbf{v}_g,\mathbf{c}_i\}$.
%This is because typically in search algorithms like A*, the states are expanded in order of
% However, if we know the exact goal node in $\mathcal{G}_\chi$ that we want to reach, say $\{\mathbf{v}_g,\mathbf{c}_g\}$, then it is possible to design a cost function and corresponding heuristic that can possibly expedite the search to the desired goal node in $\mathcal{G}_\chi$. In particular, given the goal node $\{\mathbf{v}_g,\mathbf{c}_g\}$ in $\mathcal{G}_\chi$, we propose the following cost and heuristic function
% \begin{eqnarray}
%  w_\chi(\{ \{\mathbf{v},\mathbf{c}\} \rightarrow \{\mathbf{v'},\mathbf{c}'\} \}) = w(\{\mathbf{v} \rightarrow \mathbf{v}'\}) + H \|
% \end{eqnarray}

%However,
%Paths in

\subsubsection{Planning in Low Dimensional Configuration Spaces}

\begin{figure*}
  \begin{center}
  \setlength{\fboxsep}{0.01in}
    \subfigure[Class $1$]{
      \fbox{\includegraphics[width=0.17\textwidth, trim=0 0 0 0, clip=true]{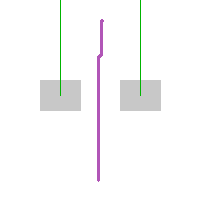}}
    }
    \subfigure[Class $2$]{
      \fbox{\includegraphics[width=0.17\textwidth, trim=0 0 0 0, clip=true]{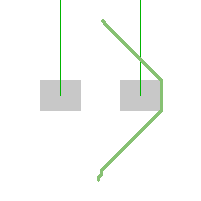}}
    }
    \subfigure[Class $3$]{
      \fbox{\includegraphics[width=0.17\textwidth, trim=0 0 0 0, clip=true]{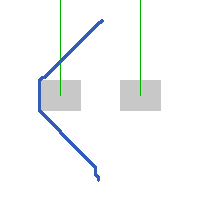}}
    }
    \subfigure[Class $4$]{
      \fbox{\includegraphics[width=0.17\textwidth, trim=0 0 0 0, clip=true]{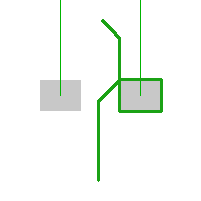}}
    }
    \subfigure[Class $5$]{
      \fbox{\includegraphics[width=0.17\textwidth, trim=0 0 0 0, clip=true]{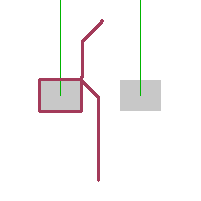}}
    }
    \subfigure[Class $6$]{
      \fbox{\includegraphics[width=0.17\textwidth, trim=0 0 0 0, clip=true]{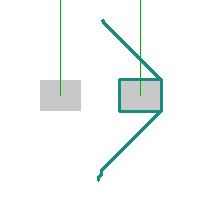}}
    }
    \subfigure[Class $7$]{
      \fbox{\includegraphics[width=0.17\textwidth, trim=0 0 0 0, clip=true]{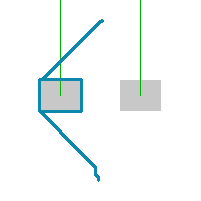}}
    }
    \subfigure[Class $8$]{
      \fbox{\includegraphics[width=0.17\textwidth, trim=0 0 0 0, clip=true]{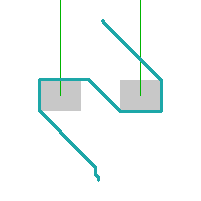}}
    }
    \subfigure[Class $9$]{
      \fbox{\includegraphics[width=0.17\textwidth, trim=0 0 0 0, clip=true]{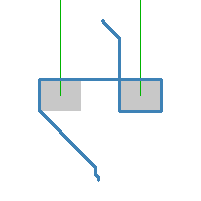}}
    }
    \subfigure[Class $10$]{
      \fbox{\includegraphics[width=0.17\textwidth, trim=0 0 0 0, clip=true]{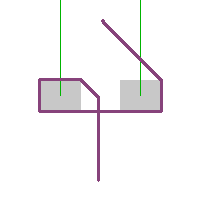}}
    }
  \end{center}
   \caption{The first $10$ homology classes of trajectories in order of length/cost. The gray regions are the obstacles. The trajectories are in different homotopy classes as well.} \label{fig:homology-bump} 
\end{figure*}

Figure~\ref{fig:homology-bump} shows a $2$-dimensional region punctured by two obstacles.
The graph $\mathcal{G}$ is constructed by uniform square discretization ($200\times 200$), placing a vertex in each cell, and by connecting the free/accessible neighboring vertices (Figure~\ref{fig:graph-8-connected}).
During the search of graph $\widehat{\mathcal{G}}$, we adopt the action `ii.' whenever we encounter a vertex of the form $\{\mathbf{v}_g,\mathbf{c}_j\}\in\widehat{\mathcal{V}}$, until we have stored $10$ paths.
\changedA{One can choose the \emph{bump $1$-form} \cite{bott1982differential} for constructing $\psi_{\widetilde{\mathcal{S}}}$ as discussed earlier. The \emph{supports} of that form are illustrated in the figure as the thin \lo{green }rays.}

\begin{figure}
  \begin{center}
    \includegraphics[width=0.55\textwidth, trim=50 130 10 130, clip=true]{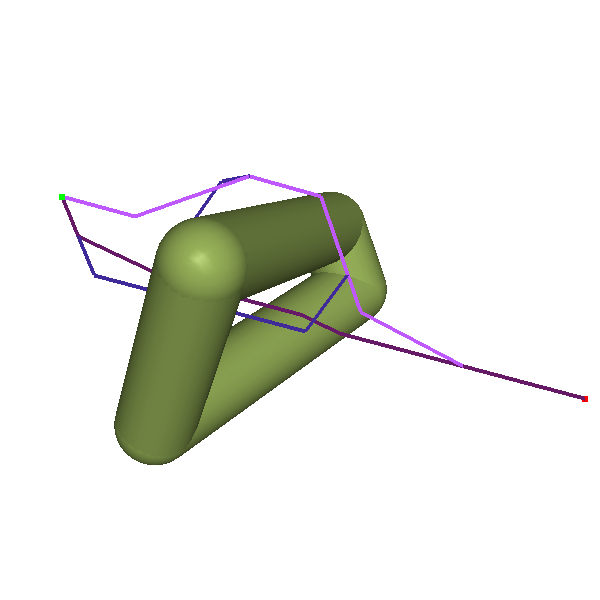}
  \end{center}
  \caption{Exploration of $3$ homology classes of robot trajectories for a $D=3$-dimensional configuration space.\label{fig:search-D3N2}}
\end{figure}

Figure~\ref{fig:search-D3N2} demonstrates an example of search for $3$ homology classes in a configuration space with $D=3$. The graph $\mathcal{G}$ is created by uniform discretization of the region of interest into $16\times 16\times 16$ cubic cells, and connecting the vertices corresponding to each cell to their immediate $26$ neighbors.

\lo{It is to be noted that the trajectories obtained in the search process as above are also in distinct homotopy classes. Although we missed some homotopy classes (\emph{e.g.} one between Figures~\ref{fig:homology-bump}(h) and \ref{fig:homology-bump}(i)), the trajectories we obtained are in different homotopy class. This is due to the Hurewicz map \cite{Hatcher:AlgTop}.
}

\begin{figure}
  \begin{center}
  \setlength{\fboxsep}{0.01in}
    \subfigure[$t=1s$]{
      \fbox{\includegraphics[width=0.42\textwidth, trim=40 90 40 90, clip=true]{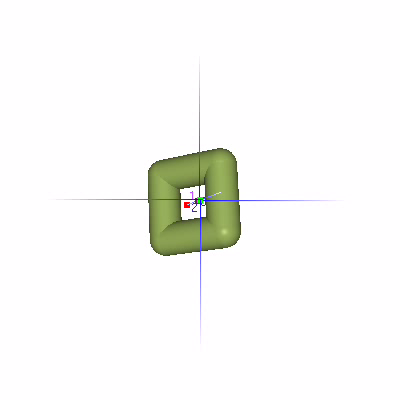}}
    } \hspace{-0.05in}
    \subfigure[$t=4s$]{
      \fbox{\includegraphics[width=0.42\textwidth, trim=10 90 70 90, clip=true]{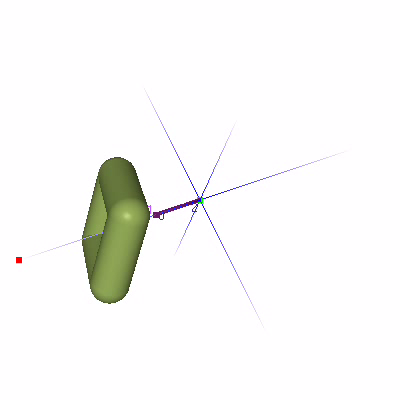}}
    } 
    \subfigure[$t=7s$]{
      \fbox{\includegraphics[width=0.42\textwidth, trim=30 90 50 90, clip=true]{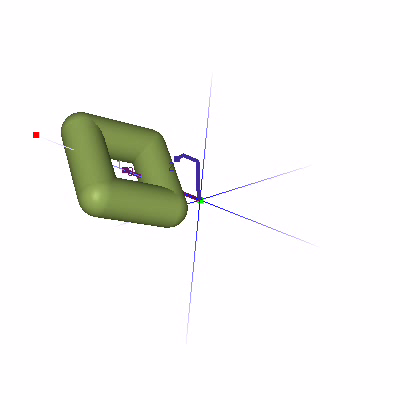}}
    } \hspace{-0.05in}
    \subfigure[$t=10s$]{
      \fbox{\includegraphics[width=0.42\textwidth, trim=40 90 40 90, clip=true]{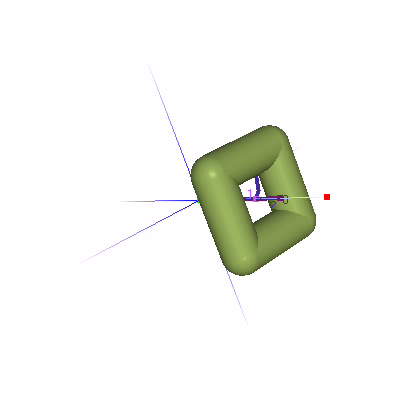}}
    } 
    \subfigure[$t=13s$]{
      \fbox{\includegraphics[width=0.42\textwidth, trim=70 110 10 70, clip=true]{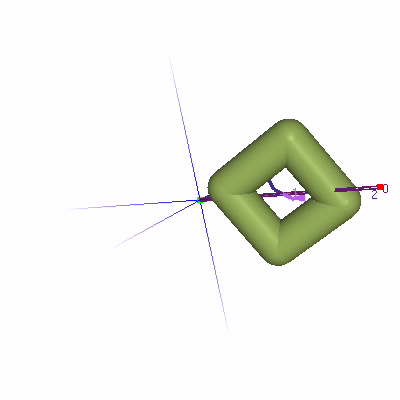}}
    } \hspace{-0.05in}
    \subfigure[$t=16s$]{
      \fbox{\includegraphics[width=0.42\textwidth, trim=60 120 20 60, clip=true]{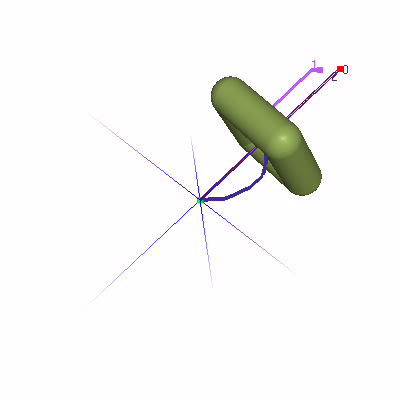}}
    } 
  \end{center}
  \caption[Screenshots from exploration of $3$ homotopy classes in a $X-Y-Z-Time$ configuration space.]{Screenshots from exploration of $3$ homotopy classes in a $X-Y-Z-Time$ configuration space. The loop-shaped obstacle is rotating about an axis. \changedA{The $X,Y$ and $Z$ axes are shown\lo{ in blue}}. Their apparent rotation is due to movement of the \emph{camera} for viewing from different angles. \label{fig:xyzt}} 
\end{figure}

\subsubsection{Exploring Paths in Different Homotopy Classes in a $4$-dimensional Space}

Just as we developed formulae for complete invariants for homology class in the $2$ and $3$ dimensional cases in \cite{planning:AURO:12}, we can now extend the formula to trajectories in higher dimensional spaces using the invariant described in Equation~\eqref{eq:final-psi-eqn}.

In this example we explore homology classes of trajectories in a $3$-dimensional space with moving obstacles. However that makes the configuration space a $4$-dimensional one consisting of the coordinates $X$, $Y$, $Z$ and $Time$.
Thus we present a result in a $X-Y-Z-Time$ configuration space where we find multiple shortest paths in different homology classes in the $4$-dimensional space. Figure~\ref{fig:xyzt} shows the exploration of $3$ homology classes in a $4$-dimensional configuration space consisting of a dynamic obstacle in $3$-dimensions.
The loop-shaped obstacle is rotating about an axis. \changedA{The $X,Y$ and $Z$ axes are shown\lo{ in blue}}. %The apparent rotation of the axes themselves is due to movement of the \emph{camera} for viewing the trajectories from different angles.
As we observe in the sequence, trajectories numbered $0$ and $1$ take off from the start coordinate (green dot) and move towards the ``center'' of the loop. They wait there while $2$ takes a different homotopy class to reach the center later. From there $0$ and $2$ head together towards the goal (red dot), while $1$ wait to take a different trajectory to the goal. Thus the $3$ trajectories are in different homotopy classes.

\lo{Once again, we note that being in different homology classes imply that the trajectories are in different homotopy classes as well.}

% \changedA{
% \subsection{Relation to universal cover}
%
% The graph $\widehat{\mathcal{G}}$ derived from graph $\mathcal{G}$ is strikingly similar to the universal cover of $\mathcal{G}$~\cite{Hatcher:AlgTop}. In fact an alternative
% }

% ----------------------------------------
\changedA{
\lo{\newpage}
\section{Extension to non-Euclidean Ambient Spaces} \label{sec:non-Eu-extension}

Let $L$ be a subspace of $(\mathbb{R}^D - \widetilde{\mathcal{S}})$. In this section we would like to compute complete invariants for homology classes of $(N-1)$-cycles in the quotient space $(\mathbb{R}^D - \widetilde{\mathcal{S}}) / L$.

We write the inclusion map as $\iota: L \hookrightarrow (\mathbb{R}^D - \widetilde{\mathcal{S}})$.
We consider $(N-1)$-chains in $C_{N-1}(\mathbb{R}^D - \widetilde{\mathcal{S}})$, and their images under the quotient map $q_{\#}: C_{\bullet}(\mathbb{R}^D - \widetilde{\mathcal{S}}) \rightarrow C_{\bullet}(\mathbb{R}^D - \widetilde{\mathcal{S}}) / C_{\bullet}(L)$. 
\Ls{
In the following proposition we consider the general pair of spaces $(X,L)$, and for generality we state it for $n$-chains.

Let us consider a $(N-1)$-chain, $\overline{\alpha} \in C_{N-1}(\mathbb{R}^D - \widetilde{\mathcal{S}})$ such that its boundary lies completely in $L$. Let us represent this boundary by $\overline{\beta} \in C_{N-2}(L)$.
This is, in general, extremely difficult to achieve. However we will consider a special condition on the relative cycles

\begin{proposition} \label{prop:relative-invariant}
Let $(X, L)$ be a pair of spaces, $\iota: L \hookrightarrow X$ be the inclusion map, and $q_{\#}: C_{\bullet}(X) \rightarrow C_{\bullet}(X) / C_{\bullet}(L)$ the quotient map for chains.
Consider $\overline{\alpha} \in C_n(X)$ such that its boundary, $\partial \overline{\alpha}$, is either empty or lies completely in $L$. Thus, $q_{\#} (\overline{\alpha})$ is a relative $n$-cycle in $(X, L)$.
Then, $[q_{\#} (\overline{\alpha})] = 0 \in H_n (X, L)$ if and only if there exists some $\overline{\beta} \in C_n(L)$ with $\iota \circ \partial \overline{\beta} = \partial \overline{\alpha}$, such that $[\overline{\alpha} - \iota \circ \overline{\beta}] = 0 \in H_n (X)$ (see Figure~\ref{fig:relative-invariant-illustration}).
\end{proposition}

% RG : NONONONO! THIS IS TRIVIAL...DO NOT SPELL OUT. PUT IT AS A "CLEARLY" IN THE PROOF OF THE COROLLARY TO FOLLOW

\begin{quoteproof}
%  Consider the following part of the long exact sequence for pair $(X , L)$,
% \[
%  H_n(L) \xrightarrow{~\iota_{*}} H_{N-1} (\mathbb{R}^D \!\!-\! \widetilde{\mathcal{S}}) \xrightarrow{~q_{*}} H_{N-1} (\mathbb{R}^D \!\!-\! \widetilde{\mathcal{S}} , L)
% \]
Suppose there exists $\overline{\beta} \in C_{N-1}(L)$ such that $[\overline{\alpha} - \iota \circ \overline{\beta}] = 0$. Using the homomorphism $q_{*}$ induced by $q_{\#}$, and noting that $q_{\#}(\iota\overline{\beta}) = 0$, we have from functoriality of homology $q_{*} ( [\overline{\alpha} - \iota \circ \overline{\beta}] ) =  [q_{\#}(\overline{\alpha})] - [q_{\#}(\iota\overline{\beta})] = [q_{\#}(\overline{\alpha})]$. Thus, $[q_{\#}(\overline{\alpha})] = 0$, concluding the `if' part of the proof.

% The proof follows directly from the functoriality of homology, and the observation that $q_{\#}(\iota\overline{\beta}) = 0$ for $N>1$. There exists a homomorphism $q_{*}: H_{N-1} (\mathbb{R}^D \!\!-\! \widetilde{\mathcal{S}}) \rightarrow H_{N-1} (\mathbb{R}^D \!\!-\! \widetilde{\mathcal{S}} , L)$ such that $q_{*} ( [\overline{\alpha} - \iota \circ \overline{\beta}] ) =  [q_{\#}(\overline{\alpha})] - [q_{\#}(\iota\overline{\beta})] = [q_{\#}(\overline{\alpha})]$.

% Thus, if $[\overline{\alpha} - \iota \circ \overline{\beta}]=0$ for some $\overline{\beta}$, we conclude immediately $[q_{\#}(\overline{\alpha})] = 0$.
\vspace{0.1in}
Next, assume $[q_{\#}(\overline{\alpha})] = 0$. Consider the following diagram with exact rows,
\begin{equation}
\begin{array}{ccccccccc}
 0 & \xrightarrow{\quad} & C_{n+1}(L) & \xrightarrow{~~\iota~~} & C_{n+1}(X) & \xrightarrow{q_{\#}} & C_{n+1}(X,L) & \xrightarrow{\quad} & 0 \\
   &       & \Big\downarrow \tilde{\partial} &        &  \Big\downarrow \partial &        & \Big\downarrow \hat{\partial} &              &  \\
 0 & \xrightarrow{\quad} & C_{n}(L) & \xrightarrow{~~\iota~~} & C_{n}(X) & \xrightarrow{q_{\#}} & C_{n}(X,L) & \xrightarrow{\quad} & 0
\end{array}
\end{equation}
%Here the rows are exact.
The proof follows from the above diagram using the following sequence of arguments:

$[q_{\#}(\overline{\alpha})] = 0$ implies $q_{\#}(\overline{\alpha}) \in C_{n}(X,L)$ is a relative boundary. Thus, there exists some $\overline{\gamma} \in C_{n+1}(X,L)$ such that $q_{\#}(\overline{\alpha}) = \hat{\partial} \overline{\gamma}$.
Due to surjectivity of $q_{\#}$ (since the rows are exact), then there exists a $\overline{A} \in C_{n+1}(X)$ such that $q_{\#}(\overline{A}) = \overline{\gamma}$.
Hence, from the commutativity of the right square, we have $q_{\#}\circ\partial(\overline{A}) = \hat{\partial}\circ q_{\#}(\overline{A}) = \hat{\partial} \overline{\gamma} = q_{\#}(\overline{\alpha})$. Hence, $q_{\#}(\overline{\alpha}-\partial\overline{A}) = 0 \in C_{n}(X,L)$.
Thus, $(\overline{\alpha}-\partial\overline{A}) \in K\!er(q_{\#}) \subseteq C_{n}(X)$.
Finally, using the exactness of the second row, there should hence exist a $\overline{\beta} \in C_{n}(L)$ such that $\iota(\overline{\beta}) = \overline{\alpha}-\partial\overline{A}$. Thus, $\overline{\alpha} - \iota(\overline{\beta}) = \partial\overline{A}$, which is a $n$-boundary. Hence proved.

For simplicity, the following diagram illustrates the quantities introduced in the above argument:
\begin{equation}
\begin{array}{ccccccccc}
 ~ & ~ & ~ & ~ & \overline{A} & \!\!\!\!\xrightarrow{q_{\#}} & \overline{\gamma} & \xrightarrow{\quad} & 0 \\
   &       & ~ &        &  \Big\downarrow \partial &        & \Big\downarrow \hat{\partial} &              &  \\
 \phantom{0} & \phantom{\xrightarrow{\quad}} & \overline{\beta} ~~~ & \xrightarrow{~~\iota~~} & \begin{array}{c} \overline{\alpha},\partial\overline{A} \\ {\scriptstyle{\text{s.t., } (\overline{\alpha}-\partial\overline{A}) \in K\!er(q_{\#})}} \end{array} & \!\!\!\!\xrightarrow{q_{\#}} & q_{\#}(\overline{\alpha}) & \xrightarrow{\quad} & 0
\end{array}
\end{equation}
%the zero in $C_{n}(X,L)$ must have some pre-preimage in $C_{n}(L)$. Thus there exists a $\overline{\beta} \in C_{n}(L)$ such that $q_{\#}\circ\iota(\overline{\beta}) = 0 = q_{\#}(\overline{\alpha}-\partial\overline{A})$.
%Thus, $$
%
%  Again, %if $[q_{\#}(\overline{\alpha})] = 0$, there should be
%  by the property of homomorphisms, there exists
%  some $\overline{\gamma} \in Z_{N-1} (\mathbb{R}^D \!\!-\! \widetilde{\mathcal{S}})$ such that $q_{*}([\overline{\gamma}]) = 0$. Define $\overline{b} = \overline{\alpha} - \overline{\gamma} \in C_{N-1} (\mathbb{R}^D \!\!-\! \widetilde{\mathcal{S}})$.
%
%Thus the `if' part of the statement follows directly. For the `only if' part, we can argue by contradiction: Suppose $[q_{\#}(\overline{\alpha})] = 0$, but there does not exist a
\end{quoteproof}

%RG I DO NOT FIND THIS FIGURE HELPFUL...
\begin{figure}
 \begin{center}
   \subfigure{ \label{fig:relative-invariant-illustration}
     \includegraphics[width=0.39\textwidth, trim=0 0 0 0, clip=true]{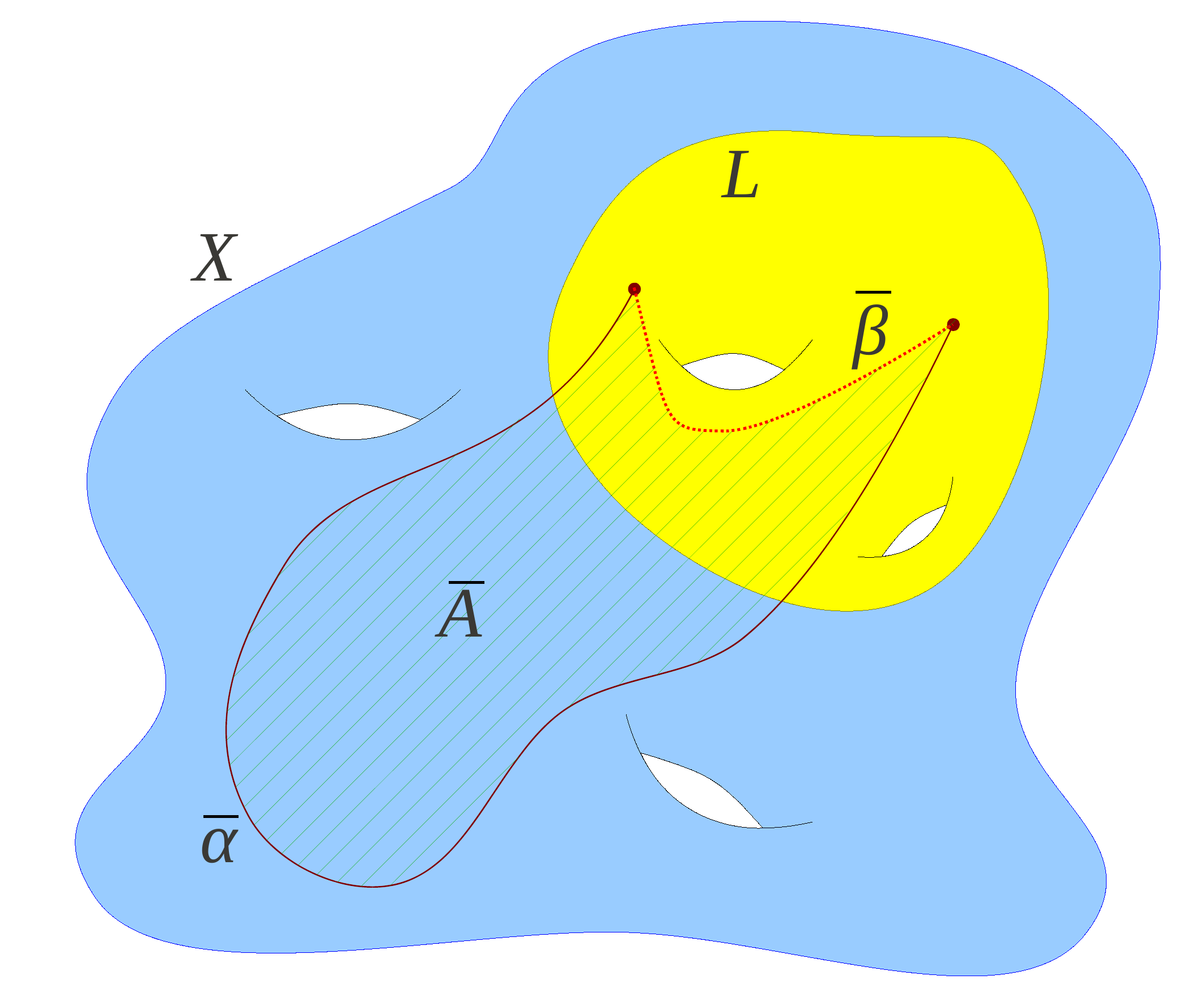}
   } \hspace{0.01in}
   \subfigure{ \label{fig:relative-invariant-traj}
     \includegraphics[width=0.55\textwidth, trim=-100 0 -100 0, clip=true]{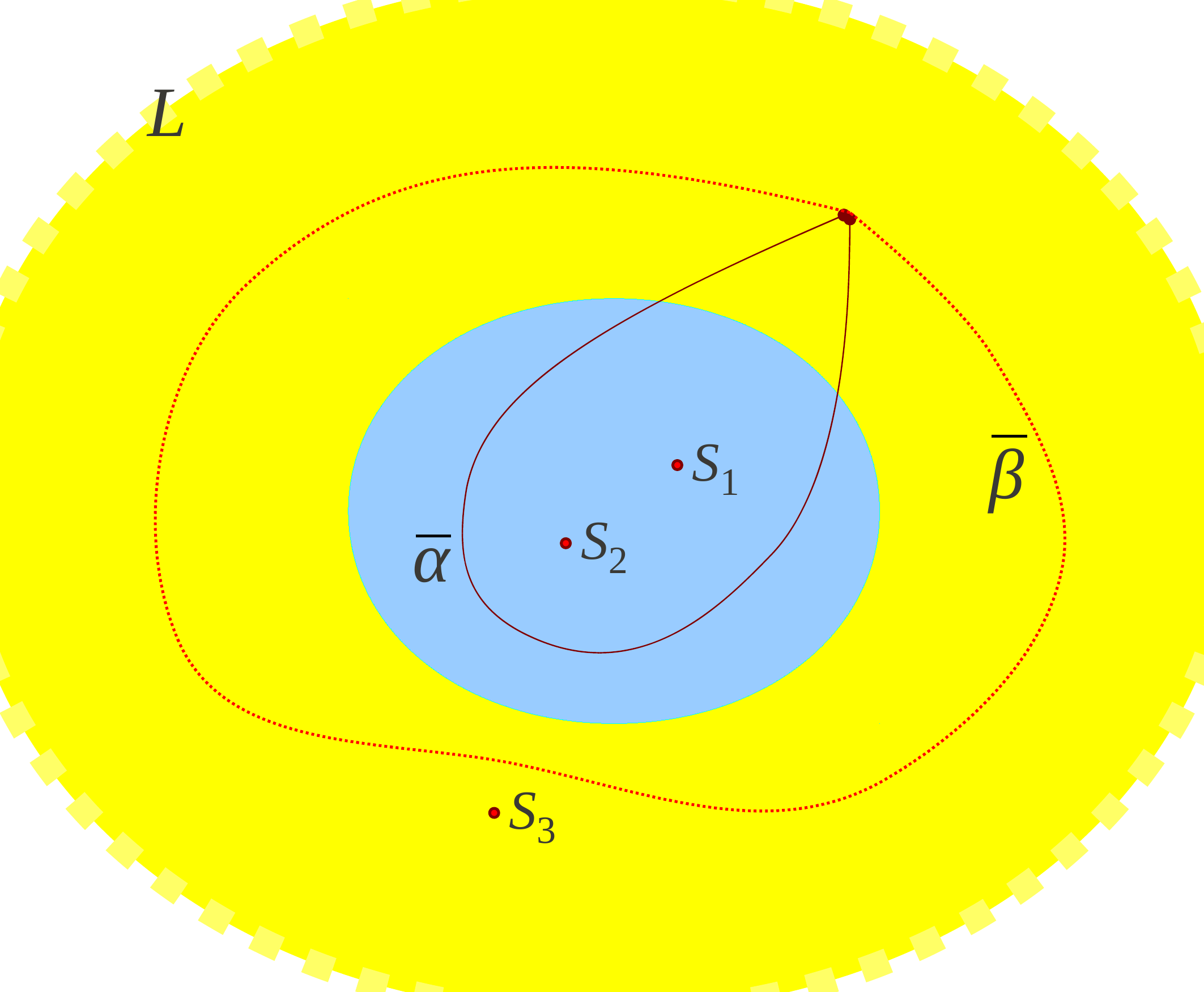}
   }
 \end{center}
 \caption{
 (left) An illustration of a relative cycle. (right). An example with $X=\mathbb{R}^2-(S_1\sqcup S_2 \sqcup S_3),~N=2$. The region, $L$, consists of everything that lies outside the small disk-shaped region, moding which out gives us the $2$-sphere with two punctures (images of $S_1$ and $S_2$). $\overline{\alpha}$ is a non-trivial cycle in $\mathbb{R}^2-(S_1\sqcup S_2 \sqcup S_3)$ since $\phi_{\widetilde{\mathcal{S}}}(\overline{\alpha} ) = \lbrack 1,1,0 \rbrack$. However it is trivial on the punctured sphere. To see this, we observe that in this case $Q = \{\scriptstyle{ [0,0,0], [\pm 1,\pm 1,0], [\pm 2,\pm 2,0], \cdots, [0,0,1], [\pm 1,\pm 1,1], \cdots, [0,0,2], \cdots, \cdots }\}$. Thus we see that $\phi_{\widetilde{\mathcal{S}}}(\overline{\alpha} ) \in Q$. A $\overline{\beta} \in C_{N-1}(L)$ corresponding to the class is shown in the figure. \label{fig:relative-invariant}}
\end{figure}
}
{
% \begin{figure}
%  \begin{center}
%   \includegraphics[width=0.5\textwidth, trim=-100 0 -100 0, clip=true]{figures/relative_example1.pdf}
%  \end{center}
%  \caption{\changedC{An example of computation in quotient space $X/L$. Here $X=\mathbb{R}^2-(S_1\sqcup S_2 \sqcup S_3),~N=2$. The region, $L$, consists of everything that lies outside the small disk-shaped region, moding which out gives us the $2$-sphere with two punctures (images of $S_1$ and $S_2$). $\overline{\alpha}$ is a non-trivial cycle in $\mathbb{R}^2-(S_1\sqcup S_2 \sqcup S_3)$ since $\phi_{\widetilde{\mathcal{S}}}(\overline{\alpha} ) = \lbrack 1,1,0 \rbrack$. However it is trivial on the punctured sphere. To see this, we observe that in this case $Q = \{\scriptstyle{ [0,0,0], [\pm 1,\pm 1,0], [\pm 2,\pm 2,0], \cdots, [0,0,1], [\pm 1,\pm 1,1], \cdots, [0,0,2], \cdots, \cdots }\}$. Thus we see that $\phi_{\widetilde{\mathcal{S}}}(\overline{\alpha} ) \in Q$. A $\overline{\beta} \in C_{N-1}(L)$ corresponding to the class is shown in the figure.} \label{fig:relative-invariant}}
% \end{figure}
}

\begin{corollary} \label{corr:relative-invariant}
 Consider $\overline{\alpha} \in C_{N-1}(\mathbb{R}^D - \widetilde{\mathcal{S}})$ such that its boundary, $\partial \overline{\alpha}$, is either empty or lies completely in $L$. % Thus, $q_{\#} (\overline{\alpha})$ is a relative $n$-cycle in $(X, L)$.
 Consider the set of all the $(N-1)$-chains in $L$ with boundary coinciding with $\partial \overline{\alpha}$ ~(if $\partial \overline{\alpha}=0$, we consider all $(N-1)$-cycles in $L$), and let $Q$ denote the \changedC{set of} $\phi_{\widetilde{\mathcal{S}}}$-image of \changedC{those}. 
 Then, $[q_{\#} (\overline{\alpha})] = 0 \in H_n (X, L)$ if and only if $\phi_{\widetilde{\mathcal{S}}}(\overline{\alpha}) \in Q$.
\end{corollary}
\begin{quoteproof}
 The statement follows directly from 
\Ls{Proposition~\ref{prop:relative-invariant} by setting $X=(\mathbb{R}^D - \widetilde{\mathcal{S}})$ and $n=N-1$ and by noting that}{the definitions of relative homology \changedC{which guarantees the existence of a $\overline{\beta}\in C_{N-1}(L)$ such that}}, $\phi_{\widetilde{\mathcal{S}}}(\overline{\alpha} - \iota \circ \overline{\beta}) = 0$ if and only if $[\overline{\alpha} - \iota \circ \overline{\beta}] = 0$. Moreover, due to the linearity of $\phi_{\widetilde{\mathcal{S}}}$, we have $\phi_{\widetilde{\mathcal{S}}}(\overline{\alpha} - \iota \circ \overline{\beta}) = 0 ~\Rightarrow~ \phi_{\widetilde{\mathcal{S}}}(\overline{\alpha} ) = \phi_{\widetilde{\mathcal{S}}}(\iota \circ \overline{\beta})$. For all computational purpose, $\iota$ becomes the identity map since we use a single coordinate chart on $(\mathbb{R}^D - \widetilde{\mathcal{S}})$.
\end{quoteproof}

One motivation for considering this kind of spaces arise from frontier-based exploration problems in robotics \cite{probRob:Thrun}, where $L$ represents the unexplored/unknown region in a configuration space, and the task at hand is to deploy robots, starting from a point in the known/explored region, to reach $L$ following different topological classes. While we do not discuss a complete exploration problem in this paper, we will describe, with example, how optimal trajectories in the different homology classes for reaching $L$ can be obtained using a graph search-based approach.
As far as implementation for search-based planning for robot trajectories is concerned, we will mostly be interested in $\overline{\alpha}$ that has empty boundary (formed by trajectories sharing the same start and goal points in $(\mathbb{R}^D - \widetilde{\mathcal{O}})$, as shown in \ls{Figure~\ref{fig:relative-invariant-traj}}{Figure~\ref{fig:relative-invariant}}). Thus the $Q$ that will be of our interest is the one for $\partial \overline{\alpha} = \emptyset$.

\subsection{Search-based Implementation}

A graph search-based algorithm, as described earlier, can once again be employed for the case with $N=2$, for finding optimal trajectories in different homology classes on $(\mathbb{R}^D - \widetilde{\mathcal{S}}) / L$. Homology classes of trajectories (which are relative chains in $C_1(\mathbb{R}^D - \widetilde{\mathcal{S}}, L)$) are \changedC{defined informally in a way similar to one in Remark~\ref{def:homology-chains}}.

The complete environment, $\mathbb{R}^D - \widetilde{\mathcal{O}}$, is discretized to create a graph, $\mathcal{G}$, as before. Edges of the graph lying in $L$ are assigned zero costs (a small positive value is used in practice for numerical stability), while for ones in the complement space is assigned the costs induced by a metric of choice (we choose the Euclidean metric of the ambient space for the example in Figure~\ref{fig:result-noneu}).
The construction of the augmented graph is similar to the construction of $\widehat{\mathcal{G}}$ as before, except that now a vertex $\{\mathbf{v},\mathbf{c}\}$ is identified with $\{\mathbf{v},\overline{\mathbf{c}}\}$ if $\mathbf{c}-\overline{\mathbf{c}} \in Q$ (where $Q$ is the set corresponding to $\partial \overline{\alpha} = \emptyset$). We call this derived graph $\widetilde{\mathcal{G}}$.

Figure~\ref{fig:result-noneu} shows an environment that is similar to the one illustrated in Figure~\ref{fig:homology-bump}, except that now everything outside a rectangular region containing the two obstacles is considered to be part of $L$ (the \ls{yellow region}{region near the boundary}, where the metric, and hence the cost of every edge is set to zero). The space under consideration is thus topologically a sphere, with \ls{the yellow region}{$L$} collapsed to a single point. For the search algorithm, we choose the same start coordinate as before (near the bottom of the environment -- \emph{almost} symmetrically placed with respect to the two obstacles), but we place the goal vertex inside $L$ (Exact choice does not matter. Although, if there were multiple path-connected components of $L$, we would have to place one goal vertex in each connected component for exploring all the homology classes).

Figures~\ref{fig:result-noneu}(a)-(e) shows exploration of first $5$ homology classes (in order of path lengths) in $(\mathbb{R}^D - \widetilde{\mathcal{O}}) / L$ by searching in $\widetilde{\mathcal{G}}$. However, we notice that in the classes $3$ and $5$, the parts of the trajectories lying in $(\mathbb{R}^D - \widetilde{\mathcal{O}} - L)$ have disconnected components. Notice that it is not possible to alter such trajectories through small variations to make them fall inside $(\mathbb{R}^D - \widetilde{\mathcal{O}} - L)$, and still remain close to optimal. This is because we use the Euclidean metric on $\mathbb{R}^D$ for length of the trajectories instead of the round metric on $\mathbb{S}^D \approxeq \mathbb{R}^D/L$.

While these solutions are technically optimal in the augmented graph, for exploration problems, where computed trajectories are not desired to have multiple connected components, we can alter the search algorithm slightly in order to obtain trajectories as shown in Figures~\ref{fig:result-noneu}(f)-(j) belonging to the same classes, but connected. Instead of searching in $\widetilde{\mathcal{G}}$, we first perform a pre-computation step where we execute a Dijkstra's search in the subgraph of $\mathcal{G}$ that lies in $L$ starting from the `goal' vertex, and compute the value of $\phi_{\widetilde{\mathcal{S}}}$ up to every other vertex in the subgraph following some path lying inside $L$ (and its boundary, $\partial L$). Let us represent that computed value corresponding to vertex $\mathbf{v}_L \in \mathcal{V}|_L$ by $p(\mathbf{v}_L)$. %Let such a vertex-value pair be $\{\mathbf{v}_L, \mathbf{b}\}$.
The main search is then performed using Dijkstra's algorithm in the subgraph of $\widehat{\mathcal{G}}$ with vertices lying inside $(\mathbb{R}^D - \widetilde{\mathcal{O}} - L)$ (and the boundary, $\partial L$), starting from the `start' vertex, and expanding vertices until the boundary between $L$ and $(\mathbb{R}^D - \widetilde{\mathcal{O}} - L)$ are reached. In addition, a vertex on the boundary, $\{\mathbf{v}'_L, \mathbf{c}\}$, is identified with $\{\mathbf{v}''_L, \mathbf{\overline{c}}\}$ if $\left( (\mathbf{c} - p(\mathbf{v}'_L)) - (\mathbf{\overline{c}} - p(\mathbf{v}''_L)) \right) \in Q$.

\begin{figure*}
  \begin{center}
  \setlength{\fboxsep}{0.01in}
    \subfigure[Class $1$]{
      \fbox{\includegraphics[width=0.17\textwidth, trim=0 0 0 0, clip=true]{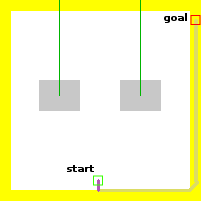}}
    }
    \subfigure[Class $2$]{
      \fbox{\includegraphics[width=0.17\textwidth, trim=0 0 0 0, clip=true]{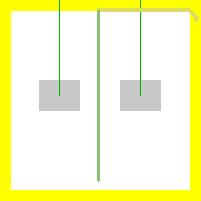}}
    }
    \subfigure[Class $3$]{
      \fbox{\includegraphics[width=0.17\textwidth, trim=0 0 0 0, clip=true]{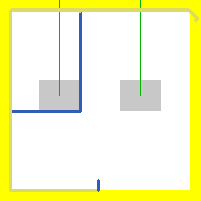}}
    }
    \subfigure[Class $4$]{
      \fbox{\includegraphics[width=0.17\textwidth, trim=0 0 0 0, clip=true]{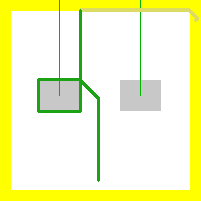}}
    }
    \subfigure[Class $5$]{
      \fbox{\includegraphics[width=0.17\textwidth, trim=0 0 0 0, clip=true]{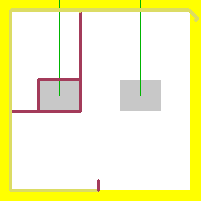}}
    }
    \subfigure[Class $1$]{
      \fbox{\includegraphics[width=0.17\textwidth, trim=0 0 0 0, clip=true]{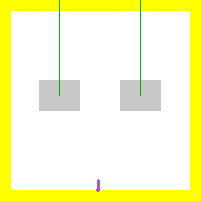}}
    }
    \subfigure[Class $2$]{
      \fbox{\includegraphics[width=0.17\textwidth, trim=0 0 0 0, clip=true]{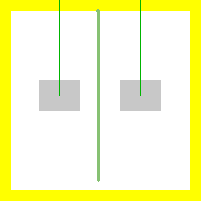}}
    }
    \subfigure[Class $3$]{
      \fbox{\includegraphics[width=0.17\textwidth, trim=0 0 0 0, clip=true]{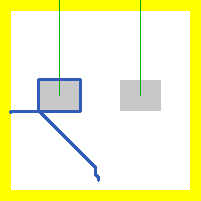}}
    }
    \subfigure[Class $4$]{
      \fbox{\includegraphics[width=0.17\textwidth, trim=0 0 0 0, clip=true]{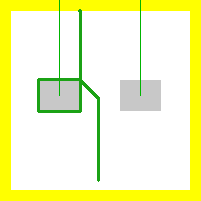}}
    }
    \subfigure[Class $5$ (The trajectory winds around the left obstacle twice.)]{
      \fbox{\includegraphics[width=0.17\textwidth, trim=0 0 0 0, clip=true]{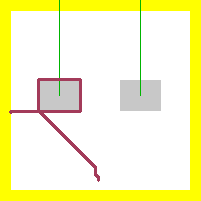}}
    }
    \subfigure[This trajectory belongs to the same class as \emph{Class $3$} (figure (c), (h)) on $(\mathbb{R}^D - \widetilde{\mathcal{O}}) / L$.]{ \label{fig:result-noneu-skipped}
      \fbox{\includegraphics[width=0.23\textwidth, trim=0 0 0 0, clip=true]{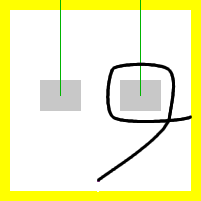}}
    }
  \end{center}
  \caption{The \ls{yellow region}{thin region near the boundary of the rectangular environment, as} shown in the figures is $L$, which we \emph{collapse} to a single point. The gray rectangles are the obstacles. (a)-(e): The first $5$ homology classes of trajectories in $(\mathbb{R}^D - \widetilde{\mathcal{O}}) / L$ connecting a given start point in $(\mathbb{R}^D - \widetilde{\mathcal{O}} - L)$, and an arbitrarily chosen point in $L$ (exact choice does not matter since we mod out $L$, which has a single path connected component) found using graph search algorithm in $\widehat{\mathcal{G}}$. (f)-(j): The solutions obtained using modified algorithm to ensure that the trajectories have single connected components in $(\mathbb{R}^D - \widetilde{\mathcal{O}} - L)$.} \label{fig:result-noneu} \vspace{-0.3in}
\end{figure*}

One interesting observation in the result of Figure~\ref{fig:result-noneu} is that apparently the search does not return any trajectory that winds around the obstacle on the right. This is because on $(\mathbb{R}^D - \widetilde{\mathcal{O}}) / L$ (\emph{i.e.} the sphere punctured by the two obstacles), a trajectory connecting the two chosen points that wind around one obstacle can be deformed over the sphere to make it wind around the other obstacle -- making them homotopic, and hence homologous. This is illustrated in Figure~\ref{fig:result-noneu-skipped}. The reason that the obstacle on the left gets preference in the result of the search algorithm is because the start coordinate is located slightly closer (by $1$ discretization unit) to the obstacle on the left than one on the right.

\section{Conclusion}

% In this paper we have considered the $D$-dimensional Euclidean space, punctured by subspaces $\widetilde{\mathcal{O}}$, and have sought to find homology classes of $(N-1)$-cycles in that space. For that we designed a set of cocycles, represented in terms of integration over differential forms, the evaluation of which over $(N-1)$-cycles gives complete set of invariants for the homology class of the cycles.
% We used the formulae thus obtained for solving optimal path planning problem in robotics, but with topological constraints, and for topological exploration of robot configuration spaces.
% A further generalization of the proposed method allows us to achieve similar objectives in ambient spaces that are not Euclidean, thus extending the potential of the proposed method to the problem of exploration in robotics.
% We have however noted that similar invariants for homology classes is very difficult to design in a way that is computationally efficient as well. At least, such invariants cannot be represented as integrals of differential forms.
% 
% %$[ \{\mathbf{v},\mathbf{c}\}, \{\mathbf{v'},\mathbf{c}'\} ]$ edge $[{}]$
% %The search algorithm, however, in this case involves two steps. In the first step we would like to determine the set $Q$. For this we consider the subgraph of $\mathcal{G}$ that lies inside $L$. We call the connected components of this subgraph $\mathcal{G}_{L_1}$, $\mathcal{G}_{L_2}$, etc. Starting from
% 
% %In construction of the graph, , from discretization of the environment, we assign a single vertex for the entire $L$.

%CONCLUSION FROM ROB'S EMAIL
The problem of optimal path planning (and its higher-dimensional
generalizations to homology $H_N$) has as prerequisite homology cycle
planning. We have addressed this precursor in the context of
obstacle-punctured Euclidean spaces. The novel features of this work
involve (1) the skeletal restructuring of the obstacles
$\widetilde{\mathcal{O}}$ to facilitate (2) the design of a set of
explicit cocycles for a complete set of invariants for the homology
class of the cycles. In this, the language of de Rham cohomology is
the critical technical step, using integration of differential forms
over cycles.
We have demonstrated the use of our methods for solving
homologically-constrained optimal path planning problems in robotics,
and topological exploration of robot configuration spaces. A further
generalization allowed us to achieve similar objectives in ambient
spaces that are not Euclidean, at the expense of an increased
computational complexity. Further work is needed to address this
issue.

}

\section*{Acknowledgements}

We gratefully acknowledge support from the ONR Antidote MURI project, grant no. N00014-09-1-1031.

\bibliography{thesis.bib}

\begin{thebibliography}{10}

\bibitem{William:Clifford}
William~E. Baylis.
\newblock {\em {Clifford (Geometric) Algebras With Applications in Physics,
  Mathematics, and Engineering}}.
\newblock Birkhäuser Boston, 1 edition, 1996.

\bibitem{yagsbpl}
Subhrajit Bhattacharya.
\newblock A template-based c++ library for large-scale graph search and
  planning, 2011.
\newblock See http://subhrajit.net/index.php?WPage=yagsbpl.

\bibitem{planning:AURO:12}
Subhrajit Bhattacharya, Maxim Likhachev, and Vijay Kumar.
\newblock Topological constraints in search-based robot path planning.
\newblock {\em Autonomous Robots}, pages 1--18, June 2012.
\newblock DOI: 10.1007/s10514-012-9304-1.

\bibitem{bott1982differential}
R.~Bott and L.W. Tu.
\newblock {\em Differential Forms in Algebraic Topology}.
\newblock Graduate texts in mathematics. Springer-Verlag, 1982.

\bibitem{Bourgault02informationbased}
Frederic Bourgault, Alexei~A. Makarenko, Stefan~B. Williams, Ben Grocholsky,
  and Hugh~F. Durrant-Whyte.
\newblock Information based adaptive robotic exploration.
\newblock In {\em in Proceedings IEEE/RSJ International Conference on
  Intelligent Robots and Systems (IROS}, pages 540--545, 2002.

\bibitem{cormen2001}
T.~H. Cormen, C.~E. Leiserson, R.~L. Rivest, and C.~Stein.
\newblock {\em Introduction to algorithms}.
\newblock MIT Press, 2nd edition, 2001.

\bibitem{dold1995lectures}
A.~Dold.
\newblock {\em Lectures on algebraic topology}.
\newblock Classics in mathematics. Springer, 2nd edition, 1995.

\bibitem{max:planning:08}
Dave Ferguson, Thomas Howard, and Maxim Likhachev.
\newblock Motion planning in urban environments.
\newblock {\em Journal of Field Robotics}, 25(11-12):939--960, 2008.

\bibitem{Differential:flanders:1989}
Harley Flanders.
\newblock {\em Differential Forms with Applications to the Physical Sciences}.
\newblock Dover Publications, New York, 1989.

\bibitem{GSL:Manual}
Mark Galassi, Jim Davies, James Theiler, Brian Gough, Gerard Jungman, Michael
  Booth, and Fabrice Rossi.
\newblock {\em {Gnu Scientific Library: Reference Manual}}.
\newblock Network Theory Ltd., February 2003.

\bibitem{GL:2006}
R.~Ghrist and S.~LaValle.
\newblock Nonpositive curvature and pareto optimal motion planning.
\newblock {\em SIAM Journal of Control and Optimization}, 45, 2006.

\bibitem{Hart-Astar}
P.~E. Hart, N.~J. Nilsson, and B.~Raphael.
\newblock A formal basis for the heuristic determination of minimum cost paths.
\newblock {\em IEEE Transactions on Systems, Science, and Cybernetics},
  SSC-4(2):100--107, 1968.

\bibitem{Hatcher:AlgTop}
Allen Hatcher.
\newblock {\em Algebraic Topology}.
\newblock Cambridge University Press, 2001.

\bibitem{KoeLik-DLite}
S.~Koenig and M.~Likhachev.
\newblock {D*} {L}ite.
\newblock In {\em Proceedings of the Eighteenth National Conference on
  Artificial Intelligence (AAAI)}, pages 476--483, 2002.

\bibitem{Armadillo}
Conrad Sanderson.
\newblock Armadillo: An open source c++ linear algebra library for fast
  prototyping and computationally intensive experiments.
\newblock Technical report, NICTA, 2010.

\bibitem{seifert1980seifert}
H.~Seifert, W.~Threlfall, J.S. Birman, and J.~Eisner.
\newblock {\em Seifert and Threlfall, A textbook of topology}.
\newblock Pure and applied mathematics. Academic Press, 1980.

\bibitem{Ste95b}
A.~Stentz and M.~Hebert.
\newblock A complete navigation system for goal acquisition in unknown
  environments.
\newblock {\em Autonomous Robots}, 2(2):127--145, 1995.

\bibitem{probRob:Thrun}
Sebastian Thrun, Wolfram Burgard, and Dieter Fox.
\newblock {\em Probabilistic Robotics (Intelligent Robotics and Autonomous
  Agents)}.
\newblock The MIT Press, 2005.

\bibitem{occlusion:yan:06}
Yan Zhou, Bo~Hu, and Jianqiu Zhang.
\newblock Occlusion detection and tracking method based on bayesian decision
  theory.
\newblock In Long-Wen Chang and Wen-Nung Lie, editors, {\em Advances in Image
  and Video Technology}, volume 4319 of {\em Lecture Notes in Computer
  Science}, pages 474--482. Springer Berlin / Heidelberg, 2006.

\end{thebibliography}
\ifx\longversion\undefined
	\bibliographystyle{siam}
\else
	\bibliographystyle{plain}
\fi

\end{document}